\documentclass[11pt]{article}
\usepackage[utf8]{inputenc}

\usepackage{float}
\usepackage{subcaption}
\usepackage{mathtools}
\usepackage{amsmath,amsxtra,amssymb,latexsym,
	epsfig,amscd,amsthm,fancybox,epsfig}
\usepackage[mathscr]{eucal}
\usepackage{graphicx}
\usepackage{epsfig} 
\usepackage{epstopdf}
\usepackage{cases}
\usepackage{color}
\usepackage{tikz}
\usetikzlibrary{positioning,arrows.meta}
\usepackage{hyperref}

\newtheorem{thm}{Theorem}[section]

\newtheorem{deff}{Definition}[section]

\newtheorem{cor}{Corollary}[section]
\newtheorem{lem}{Lemma}[section]
\newtheorem{prop}{Proposition}[section]
\theoremstyle{definition}

\theoremstyle{remark}
\newtheorem{rem}{Remark}[section]

\numberwithin{equation}{section}

\newcommand{\eps}{\varepsilon}

\newcommand{\M}{\mathcal{M}}

\newcommand{\PP}{\mathbb{P}}

\newcommand{\R}{\mathbb{R}}

\newcommand{\Lom}{\mathcal{L}}

\numberwithin{equation}{section}

\newcommand{\TV}{\mathrm{TV}}

\newcommand{\bed}{\begin{displaymath}}
\newcommand{\eed}{\end{displaymath}}
\newcommand{\bea}{\bed\begin{array}{rl}}
	\newcommand{\eea}{\end{array}\eed}

\newcommand{\barray}{\begin{array}{ll}}
	\newcommand{\earray}{\end{array}}

\def\ker{\text{Ker}}

\def\hat{\widehat}
\def\a.s{\text{\;a.s.\;}}

\date{}
\begin{document}
	\title{Identifiability of infinite location-scale Gaussian mixtures}
	\author{Nhu N. Nguyen\thanks{
Department of Mathematics and Applied Mathematical Sciences, University of Rhode Island, Kingston, RI 02881, USA, nhu.nguyen@uri.edu. This author’s research was supported in part by the National Science Foundation under grant DMS-2407669.} \and
    XuanLong Nguyen\thanks{Department of Statistics, University of Michigan, Ann Arbor, MI 48109, USA, xuanlong@umich.edu}}
	\maketitle

\begin{abstract}
The problem of identifiability of the mixing measure for a location-scale mixture of Gaussian distributions is considered. We give sufficient and necessary conditions on the compact support of the mixing measure on unknown location and scale parameters for which the mixing measure's uniqueness can be guaranteed. This is achieved by an analytic and probabilistic framework that characterizes the null space of the Gaussian mixing operator, a linear mapping that transforms a space of mixing measures with compact support into that of Gaussian mixture distributions. This framework exposes a fruitful connection between the infinite location-scale Gaussian mixtures and the exiting behavior of a space-time backward Brownian motion from a compact domain, as well as the nature of the support of the heat potential associated with such a stochastic process. 
As a consequence of this theory we obtain complete characterizations of identifiability for several natural classes of location-scale parameter sets in terms of their topological and geometric properties. Our theory also yields canonical and identifiable representations of infinite Gaussian mixtures with mixing support contained within enclosed boundaries.  

\medskip

\noindent{\bf Keywords.} 
identifiability, Gaussian mixtures, Gaussian cancellation, heat equation, adjoint heat potential,  Brownian motion, exit time, occupation time, canonical representation 

\medskip

\noindent{\bf Subject Classification.} 
60J45, 62E10, 35K05

\medskip
\noindent{\bf Running Title.}  Identifiability of Gaussian mixtures
\end{abstract}

\tableofcontents

\section{Introduction}
In this paper we study the identifiability of mixing measures that arise in a well-known and widely deployed class of statistical models known as Gaussian mixtures. Denote by
$$
\phi(x;\mu,\sigma)
=
\frac{1}{\sqrt{2\pi\sigma}}
\exp\left(
-\frac{(x-\mu)^2}{2\sigma}
\right)
$$
the Gaussian kernel with location $\mu$ and scale parameter $\sigma$, where $x,\mu\in \mathbb R$, and $\sigma>0$.

Let $S$ be a compact subset of $\mathbb R\times (0,\infty)$.
Denote by $\mathcal P(S)$ the space of Borel probability measures having support contained in $S$. 
For $G\in\mathcal P(S)$, the Gaussian mixture density function is given by
$$
p_G(x)
=
\int_S
\phi(x;\mu,\sigma)\,dG(\mu,\sigma),
\qquad x\in\mathbb{R},
$$
and $G$ is referred to as its corresponding mixing measure.
We say that the class $\mathcal P(S)$ is \emph{identifiable for Gaussian mixtures} (or simply, \emph{identifiable}) if the mapping
$
G\mapsto p_G
$
is injective on $\mathcal P(S)$; i.e.,
$$
G_1,G_2\in\mathcal P(S),
\;
p_{G_1}(x)=p_{G_2}(x),\;\;\forall 
x\in\mathbb R
\qquad 
\text{ implies }
\qquad
G_1=G_2.
$$
Equivalently, for each $x\in \mathbb R$, let $\phi_x: \mathbb R \times (0,\infty)$ denote a Gaussian likelihood function: $\phi_x(\mu,\sigma) = \phi(x; \mu,\sigma)$. 
Let $T$ be an operator in $\mathcal{M}(S)$, the space of finite signed measures on $S$, defined by
\begin{equation}\label{eq-deff-T}
    (T\nu)(x)=\int_S\phi_x(\mu,\sigma)d\nu(\mu,\sigma).
\end{equation}
We shall refer to $T$ as the \emph{Gaussian mixing operator}. Identifiability is then equivalent to that the
null space $\mathcal N(S)=\{\nu\in\mathcal M(S):T\nu= 0\}$ contains only the zero element.

This paper addresses the following question: for which parameter sets $S$
is $\mathcal P(S)$ identifiable?  Clearly, identifiability is preserved under
restriction of the parameter set, whereas non-identifiability is preserved
under enlargement.
We generally say that $\mathcal P(S)$ is identifiable;
when no confusion is likely, we may simply say that $S$ is identifiable.

\bigskip

\noindent{\bf Motivation.}
The problem is naturally connected to several areas in analysis and probability, as will be discussed in the sequel (Section \ref{sec:connect-field}). We begin, however, with its most immediate motivation from statistics, where the identifiability question has important roots. Gaussian mixtures are one of the most popular statistical models, going back to Pearson who used it to study the variety of shapes of a crab population from Naples bay \cite{pearson1894contributions}. These days Gaussian mixtures are also deployed as building blocks in large-scale generative AI models for images and texts. The Gaussian mixture model's enduring popularity is due to its flexibility in approximating complex forms of probability density functions, while retaining interpretability. Typically, a Gaussian mixture density $p_G$ is employed to represent a heterogeneous data population, where the mixing measure $G$ represents the nature of the heterogeneity \cite{lindsay1995mixture,titterington1995statistical,mclachlan2000finite,Bishop2006}. In particular, if $G$ is a discrete probability measure with $k < \infty$ supporting atoms, then $p_G$ is called a finite mixture of $k$ Gaussians, each of which provides the distribution for a latent (unobserved) data subpopulations. The supporting atoms for $G$ capture the parameters $\mu$ and $\sigma$ (also known as the mean and the variance) associated with each of the $k$ subpopulations, while the corresponding probability mass describes the proportion. Modern statistics and machine learning literature, particularly in the field of Bayesian nonparametrics, study richer classes for the mixing measure $G$, including the case $G$ has infinite support, resulting in the so-called infinite Gaussian mixtures \cite{hjort2010bayesian,ghosal2017fundamentals}. Therefore, resolving the open question of mixing measure identifiability for infinite Gaussian mixtures has a fundamental implication to the modeler's ability to interpret the model parameters, as well as the possibility and efficiency of parameter recovery from data \cite{nguyen2013convergence,guha2021posterior,nguyen2026optimal}.

\medskip

The identifiability of Gaussian mixtures exhibits a sharp distinction between
finite and unrestricted mixing.  Classical work of Teicher
\cite{teicher1963identifiability} established identifiability of finite
mixtures of Gaussian distributions, whereas the unrestricted class of
infinite location-scale Gaussian mixtures can be non-identifiable without additional structures 
\cite{Teicher1960}.  

An elementary restriction is that an identifiable parameter set must
be topologically thin: necessarily
$\operatorname{Int}(S)=\emptyset$;
see Proposition~\ref{thm-ness-cond} in the Appendix. 
Once the
topologically thick regime is excluded, identifiability is governed by
substantially finer properties of $S$.

A recent paper 
\cite{le2026partial} studies a differential-equation-based mechanism for
non-identifiability of infinite mixture models.  For Gaussian location--scale
mixtures, this addresses a topologically ``thick'' regime.  The genuinely
support-sensitive question begins beyond this regime, when
$\operatorname{Int}(S)=\emptyset$: there, the presence or absence of
identifiability can no longer be detected by interior geometry and instead
depends on the finer properties of $S$.  

It will be shown in this paper that the identifiability for infinite Gaussian mixtures whose mixing measure is supported in $S$ turns out to depend on the geometry of $S$ in a remarkably
delicate way. 
Our results reveal several striking phenomena. 
For example, identifiability cannot be characterized merely by the
``size'' or ``smoothness" of the parameter set. 
Even among topologically
thin sets, seemingly minor geometric changes can reverse the answer: adding
a single point in the ``wrong'' region may destroy identifiability, whereas
adding an entire curve in the ``right'' region may yet preserve it. The
asymmetry between the two parameter directions is equally striking: we prove
that every compact set contained in a graph of the form
$\mu=r(\sigma)$ is identifiable for arbitrary continuous function $r$, whereas there are
graphs of the form $\sigma=s(\mu)$, with $s\in C^\infty$, that are
non-identifiable. Nor is smallness in measure decisive: we exhibit
identifiable parameter sets having arbitrarily large two-dimensional
Lebesgue measure. These phenomena themselves emerge as consequences of
our general theory, which also identifies the topological and geometric mechanisms responsible for
them.

The work most directly related to the identifiability of infinite location-scale Gaussian mixtures that we are aware of is that of Bruni and Koch
\cite{bruni1985identifiability}, who studied continuous mixtures of Gaussian
distributions with both the location and scale parameters unknown.  In
particular, they claimed identifiability under regularity assumptions when
the admissible location--scale parameters are parametrized by a compact $C^1$ curve.  
This claim, however, is false in general.  We construct a counterexample;
see Proposition \ref{prop:three}. In Remark \ref{rem-Bruni}, we also identify a gap in their argument leading to this false assertion.

The limited existing results leave largely open the topological and
geometric classification problem considered here. 
To the best of our knowledge, this
paper is the first to develop a systematic theory of identifiability for
Gaussian location--scale mixtures in terms of the topology and geometry of an arbitrary 
prescribed compact parameter set $S$. Our goal is to present sufficient and necessary conditions for the 
parameter $S$ under which the identifiability property follows. Our approach in doing so is to obtain a complete
characterization of the null space $\mathcal{N}(S)$ of the Gaussian mixing operator, a theory that may be of broad interest. 
In particular, our investigation exposes a deep and fruitful connection between two-parameter Gaussian mixtures to the exiting behavior from a compact domain of a space-time backward Brownian motion, as well as questions related to the support of the adjoint heat potential associated with such stochastic processes.

\subsection{Our approach and main results}
Our approach to resolving identifiability conditions for Gaussian location--scale mixtures is by investigating the null space of the Gaussian mixing operator. Each nonzero member of the null space, which is also called a null measure if it exists, represents a configuration of signed sources whose Gaussian contributions cancel identically. The central question is how the geometry of the mixing measure's support permits or precludes such cancellation.
To this end, we develop an integrated analytic and probabilistic framework that
characterizes the entire null space, and as a consequence we obtain geometric criteria for identifiability and complete characterizations for several natural classes of parameter sets. This framework also yields canonical representations of Gaussian mixtures
with mixing support contained within enclosing boundaries, with uniqueness of the representation guaranteed whenever the chosen boundary is identifiable.

The first hint for the integrated analytic and probabilistic character of our framework can be 
pinpointed to the observation that the Gaussian kernel $\phi(x; \mu, \sigma)$ serves dual
roles: for a fixed location-scale pair $(\mu,\sigma) \in \mathbb{R} \times (0,\infty)$, it provides the probability density function for a normal random variable with mean (or location) $\mu$ and variance (or scale) $\sigma$. However, for a fixed $x \in \mathbb{R}$, it provides the Gaussian likelihood function $\phi_x: (\mu,\sigma) \mapsto \phi_x(\mu, \sigma) = \phi(x; \mu, \sigma)$. 

Taking a functional analytic perspective, the collection of Gaussian likelihood functions $\phi_x$, as $x$ varies in the observation space, is viewed as a class of test functions 
for evaluating the finite measure $\nu \in \mathcal M(S)$. By definition, every null measure $\nu$ annihilates all $\phi_x$, namely $\int \phi_x d\nu = 0$ for all $x$. A key feature of the Gaussian likelihood function comes to the fore, in its role as the fundamental solution of the heat equation: $\Lom f = 0$ for $f\in C^\infty(\mathbb{R}^2)$, where $\Lom=\partial_\sigma - \frac12 \partial^2_{\mu\mu}$ denotes the heat operator. It follows that $\nu$ annihilates the kernel of the heat operator. It can be deduced as a consequence that $\nu$ lies in the range of its adjoint $\check\Lom$. Thus, we may represent every null measure by $\nu = \check\Lom  u$, where the identity holds in the sense of distributions, and $u$ is referred to as the corresponding adjoint heat potential. The challenge lies in a precise characterization of $u$, and in explicating how a restriction on the support set $S$ of the mixing measure $\nu$ determines that of $u$, and precisely when such a restriction reduces $u$ and therefore $\nu$ to be trivial. We will show that it is possible to produce some topological or geometric conditions on $S$ to force such an outcome, which amounts to useful sufficient conditions for Gaussian mixture identifiability. However, a general and more satisfactory solution to this challenge turns out to lie in a probabilistic characterization of the null measure and its adjoint heat potential.  

\medskip

To give an initial intuition for the probabilistic characterization, each location-scale parameter point $(\mu,\sigma)$ is viewed as a source which contributes a Gaussian density profile $x\mapsto \phi(x;\mu,\sigma)$, with signed weight for the source prescribed by the mixing measure $\nu$. The identity $T\nu= 0$, then, expresses exact cancellation of these density profiles. Indeed, writing the Jordan decomposition $\nu=\nu^+-\nu^-$, we have
$$
T\nu=0
\quad\Longleftrightarrow\quad
\int_S \phi(x;\mu,\sigma)\,\nu^+(d\mu,d\sigma)
=
\int_S \phi(x;\mu,\sigma)\,\nu^-(d\mu,d\sigma)
\qquad\text{for every }x\in\mathbb R.
$$
Thus the positive and negative contributions must generate
the same Gaussian output. Equality of their total
masses is necessary but it is far from sufficient for the pointwise
cancellation. The next useful observation is that the Gaussian convolution
identity provides 
an elementary scheme for redistributing location-scale
sources while preserving their Gaussian output. For $s,t>0$,
\begin{align*}
\phi(x;0,s+t)
=
\int_{\mathbb R}
\phi(x;y,s)\phi(y;0,t)\,dy=
\mathbb E\bigl[\phi(x;B_t,s)\bigr], 
\quad\text{for every }x\in\mathbb R,
\end{align*}
where $B$ denotes a standard Brownian motion.
Thus, a source at $(0,s+t)$ can be replaced by a Gaussian
distribution of locations $y$ at a lower scale $s$.
Their difference,
$
\widetilde \nu
=
\delta_{(0,s+t)}
-
\bigl(\phi(y;0,t)\,dy\bigr)\otimes\delta_s,
$
is a nonzero signed measure satisfying $T\widetilde\nu=0$. Although this particular
construction may not satisfy the support restriction for the mixing measure in the underlying problem 
(this $\widetilde\nu$'s support is not compact), a probabilistic construction based on the random source $(B_t,s)$ may be utilized.
This suggests how a null measure $\nu$ and its adjoint heat potential $u$ 
may be suitably represented, and that a space-time Brownian motion will play a natural
role in verifying whether or not a particular Gaussian cancellation scheme is admissible 
subject to the support restriction. 

\bigskip

Next, we shall give a highlight of our main results
in more detail. In addition, Section \ref{sec:connect-field} describes several connections between the Gaussian mixture identifiability problem considered here and several interesting problems arising in other fields. We explain the similarities as well as the key differences, and in some cases, the implications that our results may entail in such settings.
The reader in a hurry may choose to proceed directly to Section \ref{sec:prelim} and beyond, and return to Section \ref{sec:connect-field} after the first reading.

\medskip
\noindent{\bf Scale separation.}
As described above, an exact (non-trivial) cancellation requires the positive and negative
contributions to agree throughout the observation space of $x$,
including their tails. We shall state our first and relatively
simple result which disallows such a possibility, by exploiting a 
particular structure of the Gaussian tails and its Fourier transform. 
This result shows that, for
compact $S$, source groups separated by a positive scale
gap cannot produce the same nonzero Gaussian profile.
Consequently, in this case, if their combined output vanishes, each group
must already have zero output. The proof uses Hardy's
uncertainty principle: a common nonzero profile would require
incompatible Gaussian decay in observational and Fourier space  (e.g., \cite[pg. 158]{stein2003fourier}).
Therefore, we show that every $\nu$ in the null space $\mathcal N(S)$
can be localized across every gap of $\pi_\sigma(S)$, the set of scale levels occurring
in $S$.
Precisely, $\nu$ annihilates all $\phi_x\mathbf 1_{\{\sigma<c\}}$, for every $c\notin \pi_\sigma(S)$. Therefore, if
$\operatorname{Int}(\pi_\sigma(S))$, the interior of $\pi_\sigma(S)$, is empty then it provides enough scale gaps to separate all scale levels and then
the Stone--Weierstrass theorem yields
that $\nu$ annihilates $\phi_xh(\sigma)$, for all continuous function $h(\sigma)$ on $\pi_\sigma(S)$.
From this, we can recover the full two-dimensional Fourier transform of $\nu$, which will be deduced to be identical to zero.
Thus, we obtain Theorem \ref{thm:thin-scale-proj} which reads that
$$
\operatorname{Int}\pi_\sigma(S)=\emptyset
\quad\Longrightarrow\quad
\mathcal P(S)\text{ is identifiable}.
$$
An interesting consequence of 
this result is that 
$S$ actually can have arbitrarily large two-dimensional Lebesgue measure while $\mathcal P(S)$ is still identifiable; see Remark \ref{rem-S-large}. It will also be shown that the identifiability of $\mathcal P(S)$ can be reduced to that of each restriction of $S$ in each connected component in the projection $\pi_\sigma(S)$; see Theorem \ref{thm:var-reduction}. Moreover,
it also shows on an arbitrary compact and countable set, countably infinite Gaussian mixtures are identifiable. This covers the infinite location-scale Gaussian mixtures corresponding to mixing measures supported by countable and compact subsets, which are typically studied in standard Bayesian nonparametric literature \cite{hjort2010bayesian,ghosal2017fundamentals}.

\medskip

\noindent{\bf Heat potentials and support localization.}
The scale-gap argument in Theorem \ref{thm:thin-scale-proj} identifies a stage at which cancellation
must already be complete if a scale gap exists. Viewing Gaussian mixing as diffusion
towards decreased scale levels, with $T\nu$ as the final profile
at level zero, for a null measure, the sources above a scale
gap have zero final output independently of those below it.
Since Gaussian diffusion cannot turn a nonzero integrable
profile into zero, their response must already vanish upon
reaching any level inside the gap.

When the permissible scales contain an interval, 
i.e., within it there is no gap for isolating the cancellation processes,
a nonzero response from sources at a given scale level may persist and subsequently be cancelled by
sources at lower scale levels. To track this possibility, we 
consider
the adjoint heat potential alluded to earlier, which
records the residual response after diffusion down to
lower scale levels, while 
$T\nu$ is the final output.
Particularly, for compact $S$, it is established that every null measure $\nu\in\mathcal N(S)$
has a potential $u$ satisfying
$$
u\in L^2(\mathbb R^2),
\qquad
\check{\Lom}u=\nu,
\qquad
\operatorname{supp}u\subset\mathcal K(S).
$$
Where, as mentioned above, $\check{\Lom}$ is the adjoint heat operator, the equation
holds in the sense of distributions, and $\mathcal K(S)$ (see Definition \ref{def-mu-hull}) is
the compact horizontal (or $\mu$-directional) hull obtained by filling the bounded
gaps between source locations at each scale level.
Thus, the (global) Gaussian cancellation condition, namely, $T\nu = 0$, has the effect of confining the intermediate response to a region determined by the spatial arrangement of $S$.
If this hull has zero planar Lebesgue measure, integrability
forces $u=0$, and the source equation gives $\nu=0$.
This yields that (cf. Theorem \ref{thm-sufficient-iden-K(S)})
$$
\mathcal K(S)\text{ has
zero area}\quad\Longrightarrow\quad\mathcal P(S)\text{ is identifiabile}.
$$
Based on the adjoint heat potential $u$, from a different angle arising in the PDE approximation literature \cite{cias2026runge,treves1967locally}, Runge approximation provides a further geometric criterion:
$\mathcal P(S)$ is identifiable if $S$ has zero planar measure
and its horizontal sections are connected or empty at a dense
set of scale levels. These criteria reveal a fundamental
asymmetry between location and scale. On compact domains,
every continuous graph $\mu=r(\sigma)$ is identifiable,
while for graphs $\sigma=s(\mu)$, if
$s$ satisfies some conditions
 as in Corollary \ref{prop:iden-graph-over-location}, then its graph is identifiable. However, we also construct
a non-identifiable graph $\sigma=s(\mu)$ with $s\in C^\infty$
(Proposition~\ref{prop:three}). 
Thus, it is noted that smoothness alone does
not determine whether or not a graph permits Gaussian cancellation.
\medskip

\noindent{\bf Brownian exit behavior and Gaussian cancellation.}
The adjoint heat potential $u$ and its support localization results determine where an intermediate response
may persist as the Gaussian cancellation progresses towards decreased scale levels. To understand whether or not the Gaussian density profiles can be generated and
cancelled using only sources in $S$, we endow the
potential with a probabilistic interpretation as follows. 
Set $K=\mathcal K(S)$
and consider the backward space--time Brownian motion
$Y_t^q=(\mu+B_t,\sigma-t)$ started at $q=(\mu,\sigma)$,
where $B$ is standard Brownian motion. We will stop the process
at its first exit from the compact $K$, allowing immediate exit.
For a null source $\nu$, it is shown that its corresponding adjoint heat potential $u$ is the density
of the signed expected occupation accumulated before exit,
with starting points weighted by $\nu$. Thus the potential
records occupation within the hull $K$. The backward space-time Brownian motion is so constructed that the associated exit
operator denoted by $E_K$ to record the combined signed exit measure
is shown to preserve the Gaussian output, i.e.,
$T(E_K\lambda)=T\lambda$, for all $\lambda\in\mathcal M(K)$.

Combining the occupation representation for $u$ with its support
localization, we obtain the exact equivalence (cf. Theorem \ref{prop:pre-annihator-rep-via-exit}):
$$
\text{For }\nu\in\mathcal M(S),\qquad
T\nu=0
\quad\Longleftrightarrow\quad
E_K\nu=0.
$$
Global Gaussian cancellation is therefore equivalent to
cancellation of the combined exit measure.
Although it appears formally that we have merely replaced one
integral condition for the null measure by another, the latter
characterization under the exit operator driven by the Brownian
motion opens up the door to identifying sources responsible for
Gaussian cancellation. Verification for Gaussian mixture identifiability 
becomes a rather simple matter as a consequence.

Indeed, call a point active if the Brownian motion $Y_t^q$ can remain
in $K$ for a positive time with positive probability. Precisely,
we introduce the active set (also known as the Brownian interior of $K$)
$$
O(K)
:=
\{q\in K:\mathbb P(\tau_K^q>0)>0\},\quad\text{where}\quad
\tau_K^q=\inf\{t>0:Y_t^q\notin K\}.
$$
Let $A=S\cap O(K)$ be the set of active points in $S$, and let
$F=\partial K\setminus S$ be the forbidden exit boundary, where $\partial K$ stands for the boundary of $K$.
Then, it is shown that the null space must admit the following representation (cf. Theorem \ref{thm-null-space-characterization}):
\begin{equation}\label{intro-eq-null}
\mathcal N(S)
=
\bigl\{
\lambda-E_K\lambda:
\lambda\in\mathcal M(A),\ (E_K\lambda)|_F=0
\bigr\}.
\end{equation}
The restriction on $F$ concerns signed contributions
collectively: individual exit laws may charge $F$, provided
their contributions cancel there. Identifiability holds
precisely when every nonzero signed active source leaves
a nonzero net exit contribution on $F$.
Then, as an application, we obtain a geometric sufficient condition as (cf. Theorem \ref{thm:complete-iden-chacterization}):
$$
S\cap O(K)=\emptyset\quad\Longrightarrow\quad\mathcal P(S)\text{ is identifiabile}.
$$
When $\partial K\subset S$, this verifiable condition is necessary and sufficient:
$$
\text{Whenever }\partial K\subset S,\qquad S\cap O(K)=\emptyset\quad\Longleftrightarrow\quad\mathcal P(S)\text{ is identifiabile}.
$$


\bigskip

\noindent{\bf Several complete characterizations.}
In Section \ref{sec:5} we give examples of classes of curved supports for which a complete characterization of identifiability can be obtained by applying the general theory described above. The application of the general theorems requires verifying some conditions about the exit behavior of the Brownian motion. Such verifications are relatively simple in most cases. Yet, all identifiability/ non-identifiability results are probably new and and in many cases the contrast is quite striking. We now turn to these results.

\medskip

\noindent{\bf Single graphs.}
Graphs illustrate the different roles of location and scale.
On compact interval domains, every continuous graph
$\mu=r(\sigma)$ is identifiable (see Corollary \ref{cor-continuous-scale-graph}). For graphs of the opposite
orientation, $\sigma=s(\mu)$, our results cover quasiconvex
profiles, concave profiles with a unique maximizer, and
single-peaked profiles satisfying the derivative conditions, among others, as in Corollary \ref{prop:iden-graph-over-location}.
However, we also construct a
non-identifiable graph $(\mu,s(\mu))$ with $s\in C^\infty$
(cf. Proposition \ref{prop:three}).
Thus smoothness alone does not result in identifiability,
even within the class of single graphs.
See Figures \ref{fig:single graph} and \ref{fig:three-branch}.

\medskip

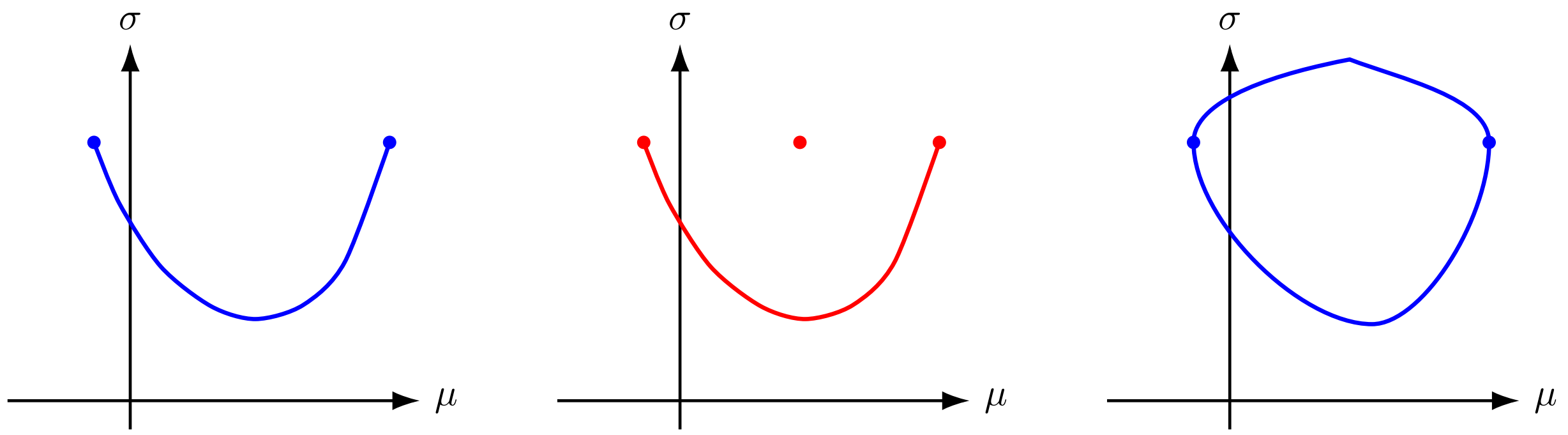
\begin{figure}[H]
\centering


\newcommand{\xmin}{-2.2}
\newcommand{\xmax}{5.2}
\newcommand{\ymin}{-0.45}
\newcommand{\ymax}{5.6}

\tikzset{
    axis/.style={
        black,
        line width=0.85pt,
        -{Latex[length=2.8mm,width=1.9mm]}
    },
    bluecurve/.style={
        blue,
        line width=1.2pt
    },
    redcurve/.style={
        red,
        line width=1.2pt
    },
    blackcurve/.style={
        black,
        line width=1.05pt
    },
    bluepoint/.style={
        circle,
        fill=blue,
        inner sep=1.35pt
    },
    redpoint/.style={
        circle,
        fill=red,
        inner sep=1.35pt
    },
    blackpoint/.style={
        circle,
        fill=black,
        inner sep=1.35pt
    }
}


\begin{subfigure}[t]{0.31\textwidth}
\centering

\resizebox{0.94\linewidth}{!}{%
\begin{tikzpicture}[x=0.56cm,y=0.64cm]

    \draw[axis]
        (0,\ymin) -- (0,\ymax)
        node[above] {$\sigma$};

    \draw[axis]
        (\xmin,0) -- (\xmax,0)
        node[right] {$\mu$};

    \draw[bluecurve]
        plot[smooth] coordinates {
            (-0.65,4.05)
            (-0.20,3.10)
            (0.55,2.10)
            (1.45,1.48)
            (2.25,1.28)
            (3.10,1.50)
            (3.85,2.18)
            (4.65,4.05)
        };

    \node[bluepoint] at (-0.65,4.05) {};
    \node[bluepoint] at (4.65,4.05) {};

\end{tikzpicture}%
}

\caption{U-shape satisfies Corollary \ref{cor:U-shape} is identifiable.}

\end{subfigure}
\hfill
\begin{subfigure}[t]{0.31\textwidth}
\centering

\resizebox{0.94\linewidth}{!}{%
\begin{tikzpicture}[x=0.56cm,y=0.64cm]

    \draw[axis]
        (0,\ymin) -- (0,\ymax)
        node[above] {$\sigma$};

    \draw[axis]
        (\xmin,0) -- (\xmax,0)
        node[right] {$\mu$};

    \draw[redcurve]
        plot[smooth] coordinates {
            (-0.65,4.05)
            (-0.20,3.10)
            (0.55,2.10)
            (1.45,1.48)
            (2.25,1.28)
            (3.10,1.50)
            (3.85,2.18)
            (4.65,4.05)
        };

    \node[redpoint] at (-0.65,4.05) {};
    \node[redpoint] at (4.65,4.05) {};

    \node[redpoint] at (2.15,4.05) {};

\end{tikzpicture}%
}

\caption{Adding an arbitrary point or points in the top cap makes the support set non-identifiable. See Corollary \ref{cor:domain-exit-balayage}.}

\end{subfigure}
\hfill
\begin{subfigure}[t]{0.31\textwidth}
\centering

\resizebox{0.94\linewidth}{!}{%
\begin{tikzpicture}[x=0.56cm,y=0.64cm]

    \draw[axis]
        (0,\ymin) -- (0,\ymax)
        node[above] {$\sigma$};

    \draw[axis]
        (\xmin,0) -- (\xmax,0)
        node[right] {$\mu$};

    \draw[bluecurve]
        (-0.65,4.05)
        .. controls (-0.65,4.75) and (0.95,5.15) ..
            (2.15,5.35)
        .. controls (3.05,5.05) and (4.65,4.75) ..
            (4.65,4.05)
        .. controls (4.65,2.85) and (3.45,1.20) ..
            (2.55,1.20)
        .. controls (1.25,1.20) and (-0.65,2.85) ..
            (-0.65,4.05);

    \node[bluepoint] at (-0.65,4.05) {};
    \node[bluepoint] at (4.65,4.05) {};

\end{tikzpicture}%
}

\caption{Completing the curve to obtain a closed curve may result in either  identifiability or non-identifiability character depending the flatness of the top point. 
See the proof of Corollary \ref{thm:pq-threshold}.
}

\end{subfigure}

\caption{Examples of U-shaped support sets.}
\label{fig:U-shape}

\end{figure}

\noindent{\bf U-shape supports.}
A support formed by two-branched curves joined by a bottom
segment is identifiable under the one-sided regularity
condition of Corollary \ref{cor:U-shape}.
This conclusion is sensitive
to the missing upper cap: adding any point in its relative
interior destroys identifiability, since that point's
Gaussian output can be reproduced by its Brownian exit
measure on the original U-shape (see Corollary \ref{cor:domain-exit-balayage}).
Suitable curve extensions can nevertheless preserve
identifiability when every point of the enlarged support
remains an immediate-exit point of its horizontal hull.
Thus the effect of enlarging a support depends on how
the added parameters change its exit behavior: a single
cap point can enable cancellation, while an entire
additional curve need not.
See Figure \ref{fig:U-shape}.

\medskip

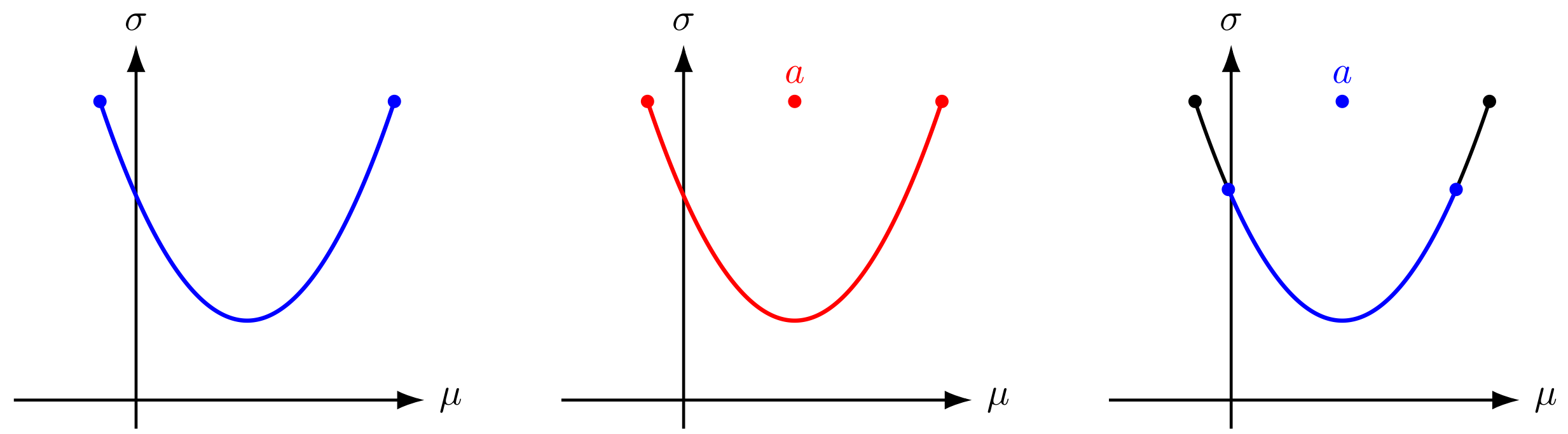
\begin{figure}[H]
\centering


\newcommand{\xmin}{-2.2}
\newcommand{\xmax}{5.2}
\newcommand{\ymin}{-0.45}
\newcommand{\ymax}{5.6}

\newcommand{\muzero}{2.0}
\newcommand{\sigmazero}{1.25}
\newcommand{\epslevel}{4.7}

\newcommand{\xleft}{-0.65}
\newcommand{\xright}{4.65}

\newcommand{\parabcoef}{0.49128}

\newcommand{\xblueleft}{-0.05}
\newcommand{\xblueright}{4.05}

\newcommand{\f}[1]{%
  \parabcoef*((#1)-\muzero)^2+\sigmazero%
}

\tikzset{
    axis/.style={
        black,
        line width=0.85pt,
        -{Latex[length=2.8mm,width=1.9mm]}
    },
    bluecurve/.style={
        blue,
        line width=1.2pt
    },
    redcurve/.style={
        red,
        line width=1.2pt
    },
    blackcurve/.style={
        black,
        line width=1.05pt
    },
    bluepoint/.style={
        circle,
        fill=blue,
        inner sep=1.35pt
    },
    redpoint/.style={
        circle,
        fill=red,
        inner sep=1.35pt
    },
    blackpoint/.style={
        circle,
        fill=black,
        inner sep=1.35pt
    }
}


\begin{subfigure}[t]{0.31\textwidth}
\centering

\resizebox{0.94\linewidth}{!}{%
\begin{tikzpicture}[x=0.56cm,y=0.64cm]

    \draw[axis]
        (0,\ymin) -- (0,\ymax)
        node[above] {$\sigma$};

    \draw[axis]
        (\xmin,0) -- (\xmax,0)
        node[right] {$\mu$};

    \draw[
        bluecurve,
        domain=\xleft:\xright,
        samples=150,
        smooth
    ]
        plot (\x,{\f{\x}});

    \node[bluepoint]
        at (\xleft,{\f{\xleft}}) {};

    \node[bluepoint]
        at (\xright,{\f{\xright}}) {};

\end{tikzpicture}%
}

\caption{Identifibility of a parabola because parabola is convex. 
}

\end{subfigure}
\hfill
\begin{subfigure}[t]{0.31\textwidth}
\centering

\resizebox{0.94\linewidth}{!}{%
\begin{tikzpicture}[x=0.56cm,y=0.64cm]

    \draw[axis]
        (0,\ymin) -- (0,\ymax)
        node[above] {$\sigma$};

    \draw[axis]
        (\xmin,0) -- (\xmax,0)
        node[right] {$\mu$};

    \draw[
        redcurve,
        domain=\xleft:\xright,
        samples=150,
        smooth
    ]
        plot (\x,{\f{\x}});

    \node[redpoint]
        at (\xleft,{\f{\xleft}}) {};

    \node[redpoint]
        at (\xright,{\f{\xright}}) {};

    \node[redpoint]
        at (\muzero,\epslevel) {};

    \node[
        red,
        above=1.5pt
    ]
        at (\muzero,\epslevel)
        {$a$};

\end{tikzpicture}%
}

\caption{Non-identifibility because starting from the source point $a$, space-time backward Brownian motion can hit the parabola in positive and finite time.}

\end{subfigure}
\hfill
\begin{subfigure}[t]{0.31\textwidth}
\centering

\resizebox{0.94\linewidth}{!}{%
\begin{tikzpicture}[x=0.56cm,y=0.64cm]

    \draw[axis]
        (0,\ymin) -- (0,\ymax)
        node[above] {$\sigma$};

    \draw[axis]
        (\xmin,0) -- (\xmax,0)
        node[right] {$\mu$};

    \draw[
        blackcurve,
        domain=\xleft:\xblueleft,
        samples=70,
        smooth
    ]
        plot (\x,{\f{\x}});

    \draw[
        bluecurve,
        domain=\xblueleft:\xblueright,
        samples=130,
        smooth
    ]
        plot (\x,{\f{\x}});

    \draw[
        blackcurve,
        domain=\xblueright:\xright,
        samples=70,
        smooth
    ]
        plot (\x,{\f{\x}});

    \node[blackpoint]
        at (\xleft,{\f{\xleft}}) {};

    \node[blackpoint]
        at (\xright,{\f{\xright}}) {};

    \node[bluepoint]
        at (\xblueleft,{\f{\xblueleft}}) {};

    \node[bluepoint]
        at (\xblueright,{\f{\xblueright}}) {};

    \node[bluepoint]
        at (\muzero,\epslevel) {};

    \node[
        blue,
        above=1.5pt
    ]
        at (\muzero,\epslevel)
        {$a$};

\end{tikzpicture}%
}

\caption{Identifiability because there is no active source in $S$.}

\end{subfigure}

\caption{A complete characterization of identifiability of parabola with a source point.}
\label{fig:parabola}

\end{figure}

\noindent{\bf Parabola with a source point.}
The parabolic U-shape exhibits a sharper support threshold.
For the parabola arc $\Gamma$ as in Corollary \ref{cor:U-shape-sharp-threshold},
adding any point $a$
in the relative interior of its upper cap produces a
non-identifiable support $S=\Gamma\cup\{a\}$, while every
proper compact subset of $S$ is identifiable; see Figure \ref{fig:parabola}..
The Gaussian output of $a$ is reproduced by a Brownian
exit measure with full support on $\Gamma$.
Identifiability of $\Gamma$ makes this representation unique,
so the entire arc is required for cancellation.
More generally, the same minimality holds whenever a source
point outside an identifiable compact set admits an equivalent
mixture representation with full support on that set (see Corollary \ref{cor:one-point-full-support-target}).


\medskip

\begin{figure}[H]
\centering


\newcommand{\xmin}{-1.55}
\newcommand{\xmax}{1.55}
\newcommand{\ymin}{-0.45}
\newcommand{\ymax}{3.75}

\newcommand{\szero}{2.0}

\tikzset{
    axis/.style={
        black,
        line width=0.85pt,
        -{Latex[length=2.8mm,width=1.9mm]}
    },
    bluecurve/.style={
        blue,
        line width=1.2pt,
        line join=round
    },
    redcurve/.style={
        red,
        line width=1.2pt,
        line join=round
    },
    blackcurve/.style={
        black,
        line width=1.05pt
    }
}


\begin{subfigure}[t]{0.46\textwidth}
\centering

\resizebox{0.75\linewidth}{!}{%
\begin{tikzpicture}[x=1.55cm,y=1.18cm]

    \draw[axis]
        (0,\ymin) -- (0,\ymax)
        node[above] {$\sigma$};

    \draw[axis]
        (\xmin,0) -- (\xmax,0)
        node[right] {$\mu$};

    \draw[redcurve]
        plot[domain=-1:1, samples=260, smooth]
            (\x,{2 + pow(max(0,1 - pow(abs(\x),2.1)),1.0/3.0)})
        --
        plot[domain=1:-1, samples=260, smooth]
            (\x,{2 - pow(max(0,1 - pow(abs(\x),2.1)),1.0/3.0)})
        -- cycle;

\end{tikzpicture}%
}

\caption{The case $p>2$: $S_{p,q}$ is non-identifiable.}

\end{subfigure}
\hfill
\begin{subfigure}[t]{0.46\textwidth}
\centering

\resizebox{0.75\linewidth}{!}{%
\begin{tikzpicture}[x=1.55cm,y=1.18cm]

    \draw[axis]
        (0,\ymin) -- (0,\ymax)
        node[above] {$\sigma$};

    \draw[axis]
        (\xmin,0) -- (\xmax,0)
        node[right] {$\mu$};

    \draw[bluecurve]
        plot[domain=-1:1, samples=260, smooth]
            (\x,{2 + pow(max(0,1 - pow(abs(\x),1.1)),1.0/3.0)})
        --
        plot[domain=1:-1, samples=260, smooth]
            (\x,{2 - pow(max(0,1 - pow(abs(\x),1.1)),1.0/3.0)})
        -- cycle;

\end{tikzpicture}%
}

\caption{The case $1\leq p\leq 2$: $S_{p,q}$ is identifiable.
}
\end{subfigure}

\caption{Identifiability characterization for the convex curves of the form
$S_{p,q}
=
\{(\mu,\sigma)\in\mathbb R\times(0,\infty): |\mu|^p+|\sigma-s_0|^q=1\}$.
}
\label{fig:spq}

\end{figure}
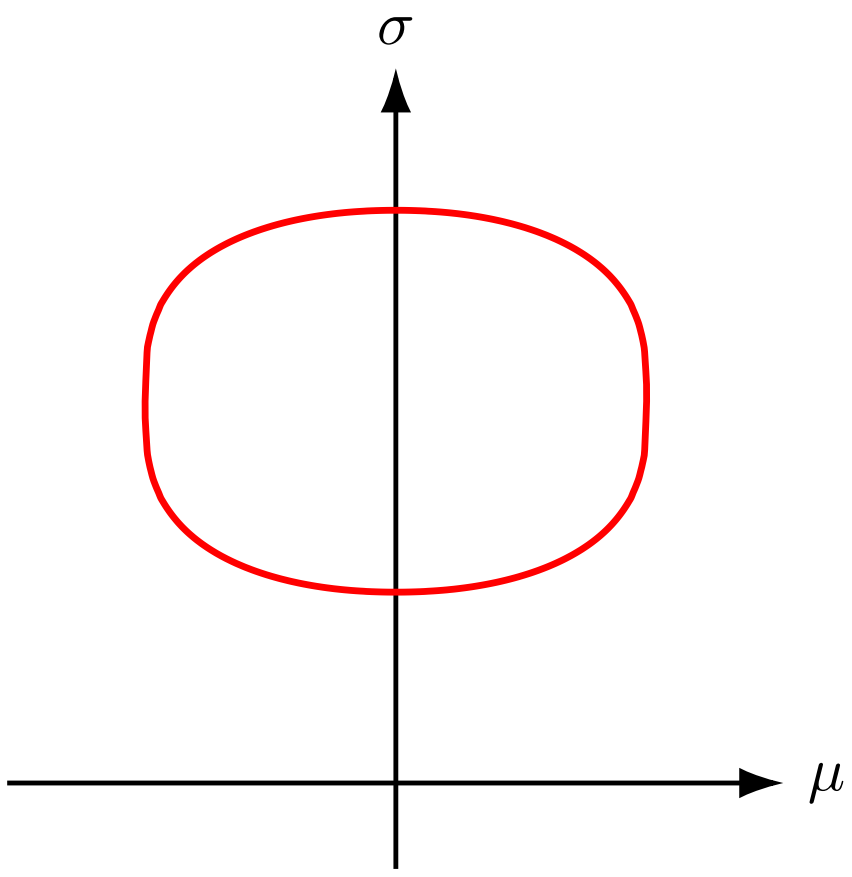

\noindent{\bf Convex curves.}
We also exhibit a sharp threshold for a class of convex curves. This shows that
even among analytic, strictly convex curves,
flatness at the uppermost point decides identifiability.
	Let
	$
	p,q\ge1,
	\;
	s_0>1,
	$
	and define
	$$
	S_{p,q}
	:=
	\left\{
	(\mu,\sigma)\in\mathbb R\times(0,\infty):
	|\mu|^p+|\sigma-s_0|^q=1
	\right\}.
	$$
Corollary \ref{thm:pq-threshold} shows if $p>2$, $\mathcal P(S_{p,q})$ is non-identifiable while if $p\leq 2$ then $\mathcal P(S_{p,q})$ is identifiable.
See Figure \ref{fig:spq}.

\bigskip

\noindent{\bf Identifiability for support in parameterized curves.}
A natural question is whether identifiability can be transported through
 transformations of the parameter space.  
We show that identifiability itself exhibits a nontrivial invariance under a
class of M\"obius transformations.  
Precisely, the M\"obius transformation
	\begin{equation*}
		\mathcal T(\mu,\sigma)
		=
		\left(
		\frac{\lambda\mu+\eta\sigma+\xi}
		{\gamma\sigma+\delta},
		\,
		\frac{\alpha\sigma+\beta}
		{\gamma\sigma+\delta}
		\right),
	\end{equation*}
	where
	$
	\alpha,\beta,\gamma,\delta,
	\lambda,\eta,\xi\in\mathbb R,
	$
    $
		\alpha\delta-\beta\gamma
		=
		\lambda^2
		>
		0,
        $
	$	\gamma\sigma+\delta\ne0,
		\;
		\frac{\alpha\sigma+\beta}
		{\gamma\sigma+\delta}>0,\;
			\forall
		\sigma\in\pi_\sigma(S),
        $
	satisfies
	\begin{equation*}
		\mathcal P(S)\text{ is identifiable}
		\quad\Longleftrightarrow\quad
		\mathcal P(\mathcal T(S))\text{ is identifiable}.
	\end{equation*}

\bigskip

\noindent{\bf Canonical boundary representations.}
Our theory also yields a representation principle
for Gaussian mixtures. For a rectangle or closed disk
$K\subset\mathbb R\times(0,\infty)$, every
$G\in\mathcal P(K)$ has a unique representative
$G_\Gamma=E_KG$ on a fixed boundary curve $\Gamma$ which 
preserves the Gaussian mixture output, i.e., $p_{G_\Gamma}=p_G$. 
For a rectangle, the bottom and
two vertical sides suffice; for a disk, every boundary
point is an immediate-exit point and the entire circle
is needed. In both cases, $\Gamma$ is the smallest compact
subset of $K$ generating its full mixture class; see Theorem \ref{prop:rep-mixing-to-boundary}.

In other words, every Gaussian mixture with compact parameter support can be represented \emph{uniquely} by a Gaussian mixture whose mixing measure is supported on a fixed boundary curve of a rectangular or disk-shaped enclosure, cf. Corollary \ref{cor:circle-representation}. This gives curved supports a direct modeling role: a two-dimensional mixing distribution can be replaced by a uniquely determined
measure on a chosen curve without changing the observed
mixture distribution. This uniqueness concerns the boundary representative;
the original mixing measure need not be identifiable. Such canonical representations of mixing measures supported on identifiable curves provide
an attractive and justifiable option for a probabilistic
modeler working with infinite Gaussian mixture models.

\subsection{Connections to other fields}\label{sec:connect-field}
Although this work is motivated by the identifiability question for Gaussian
location--scale mixtures, the same mathematical structure appears in several
different settings.  In particular, the problem admits reformulations
and interpretations in the studies of Fredholm integral equations, restricted Fourier--Laplace transforms,
parabolic inverse-source problems, caloric approximation, generalized moments,
and characteristic kernels.  We briefly describe these connections below.

Recall that for a finite signed measure $\nu\in\mathcal (S)$, 
define
$
T\nu(x)
:=
\int_S \phi_x(\mu,\sigma)\,d\nu(\mu,\sigma)
$. 
Thus, $S$ is identifiable precisely if for every $\nu\in\mathcal M(S)$,
$
T\nu= 0$
implies
$
\nu=0
$.

\medskip
\noindent
\textbf{A geometric view of first-kind integral equations.}
Recovering a mixing measure $G$ supported on $S$ from a given mixture density $p_G$ amounts to solving the measure-valued Fredholm integral equation of the first kind
$$
TG(x)=\int_S\phi_x(\mu,\sigma)\,dG(\mu,\sigma)=p_G(x),
$$
see e.g.,
\cite{Kress2014,bredies2013inverse} and the references therein.
Our results provide geometric criteria for injectivity and
describe the entire null space through the exit operator.
This characterization also determines all alternative
mixing measures compatible with a given output: for any
$G_0\in\mathcal P(S)$,
$
\{G\in\mathcal P(S):p_G=p_{G_0}\}
=
(G_0+\mathcal N(S))\cap\mathcal P(S).
$
Thus, support geometry governs both uniqueness and the
structure of the ambiguity when uniqueness fails.


\medskip
\noindent
\textbf{Fourier--Laplace uniqueness and Heisenberg uniqueness pairs.}
Taking the Fourier transform in the observation variable $x$ gives
$$
\widehat{T\nu}(t)
=
\int_S e^{it\mu-\frac12\sigma t^2}\,d\nu(\mu,\sigma),
\qquad t\in\mathbb R.
$$
Thus identifiability amounts to determining a measure supported
on $S$ from the restriction of its two-variable Fourier--Laplace
transform to the parabolic curve
$
\Lambda_{\mathrm{par}}
=
\{(it,-t^2/2):t\in\mathbb R\}
\subset\mathbb C^2.
$
The transform has Fourier dependence on location and Laplace
dependence on scale.

This is a Fourier--Laplace analogue of the geometric uniqueness
problems studied through Heisenberg uniqueness pairs
(see e.g., \cite{hedenmalm2011heisenberg,hedenmalm2020klein}),
where a support constraint permits recovery of a measure from
its Fourier transform on a thin observation set.
These uniqueness problems connect support geometry with completeness of exponential systems and uniqueness for PDEs, including the Klein–Gordon equation \cite{hedenmalm2011heisenberg}
and the Helmholtz and Laplace equations
\cite{fernandez2019heisenberg}.
In our setting, the observation set is a fixed complex parabolic curve, and our results give geometric criteria on the support $S$ under which the restricted Fourier–Laplace transform determines the mixing measure. 


Viewed in this way, our results give geometric criteria for
a fixed parabolic observation curve to determine measures
on a prescribed support $S$.
The exit-operator characterization also describes the entire
null space of this restriction through signed sources and
their Brownian exit measures, identifying all cancellations
that remain when uniqueness fails.

\medskip
\noindent\textbf{Inverse heat-source problems with final-time data.}
From its definition, observe that
$T\nu$ admits a direct inverse-source interpretation: regarding $\sigma$
as elapsed diffusion time, $T\nu$ is the terminal temperature profile
generated by a measure-valued space--time source $\nu$.  Thus,
identifiability is equivalent to injectivity of the final-state map on
sources supported on $S$. 
Inverse problems for parabolic equations with final-time overdetermination
form a classical subject; see, in particular,
\cite{isakov1991inverse,bushuyev1995global,yamamoto2009carleman} and the monograph \cite{isakov1990inverse}.
Without structural restrictions on the source, however, a single terminal
profile is far from sufficient for uniqueness: for instance, any nonzero
compactly supported function $w$ vanishing near the initial and terminal
times produces an invisible source
$f=\left(\partial_\sigma-\frac12\partial_{\mu\mu}\right)w$,
whose terminal field vanishes identically.
The question arising here is therefore whether \emph{support geometry alone}
can restore uniqueness.  Our identifiability problem asks precisely which
space--time (i.e., location--scale) sets $S$ exclude such invisible sources, thereby turning an
intrinsically underdetermined inverse heat-source problem into a geometric
uniqueness problem.
Our results, therefore, provide a complete answer for such a question.

\medskip
\noindent
\textbf{Caloric polynomial approximation.}
Expanding the caloric exponentials gives
$$
e^{z\mu+\frac12 z^2\sigma}
=
\sum_{n=0}^{\infty}
H_n(\mu,\sigma)\frac{z^n}{n!},
$$
where $H_n$ is the $n$-th heat polynomial.  Heat polynomials form a
classical family of polynomial solutions of the heat equation and have long
played a role in expansions and approximation of caloric functions; see e.g.,
\cite{rosenbloom1959expansions,widder1962analytic}.

Identifiability is equivalent to the uniform approximation property
$
\overline{\operatorname{span}\{H_n|_S:n\ge0\}}
=
C(S),
$
the space of continuous functions on (the compact) $S$; see Theorem \ref{thm:iden-characterization}.
Thus every identifiability theorem in this paper is simultaneously a
uniform approximation theorem for heat polynomials on $S$, while every
non-identifiability result exhibits a compact set on which this family fails
to be dense.

This places the problem within the broader theory of approximation by
solutions of partial differential equations, including Runge-type
approximation for the heat equation
\cite{treves1967locally,jones1975approximation,enciso2019approximation}.
There is, however, an important distinction: classical Runge theorems
approximate local caloric solutions by global ones, whereas here one asks
whether the much more rigid family of caloric polynomials can approximate
\emph{arbitrary} continuous data on $S$.  Unlike ordinary polynomials, the
heat polynomials do not form an algebra, so Stone--Weierstrass theorem does not
apply; the answer is instead governed by the geometry of $S$.


\medskip
\noindent
\textbf{A generalized moment problem.}
The heat-polynomial expansion also gives a moment-theoretic formulation.
For a finite measure $G$ supported on $S$, define its parabolic moments by
$m_n(G)
:=
\int_S H_n(\mu,\sigma)\,dG(\mu,\sigma),
\; n\ge0.$
Then identifiability of $\mathcal P(S)$ is equivalent to determinacy of the
measure $G$ by the sequence $\{m_n(G)\}_{n\ge0}$.  Thus the problem fits
naturally into the theory of generalized moment problems, where a measure is
to be recovered from its integrals against a prescribed family of test
functions
\cite{kreuin1977markov,KarlinStudden1966}.

The structure of the available moments is, however, highly degenerate from
the viewpoint of the classical multivariate moment problem
\cite{ShohatTamarkin1943,PutinarVasilescu1999}.  At parabolic degree
$j+2k=n$, the full collection of monomial moments is replaced by the single
caloric combination
$m_n(G)=\int_S H_n(\mu,\sigma)\,dG(\mu,\sigma)$;
for instance,
$m_2(G)=\int_S(\mu^2+\sigma)\,dG$
contains only the combination of $\int\mu^2\,dG$ and $\int\sigma\,dG$,
rather than the two moments separately.  Our results therefore identify
geometries of $S$ for which this highly compressed moment system is
nevertheless determining, and geometries for which determinacy fails.

\medskip
\noindent
\textbf{Characteristic kernels and kernel mean embeddings.}
There is also a natural reproducing-kernel formulation. Define the symmetric and positive definite kernel, for 
$\theta=(\mu,\sigma)$ and $\theta'=(\mu',\sigma')$,
$$
K(\theta,\theta')
:=
\int_{\mathbb R}\phi_x(\theta)\phi_x(\theta')\,dx
=
\frac{1}{\sqrt{2\pi(\sigma+\sigma')}}
\exp\!\left(
-\frac{(\mu-\mu')^2}{2(\sigma+\sigma')}
\right).
$$
This is the $L^2$ probability-product kernel between the corresponding
Gaussian distributions \cite{jebara2004probability}.  Therefore, for any finite signed
measure $\nu$ supported on $S$,
$$
\|T \nu\|_{L^2(\mathbb R)}^2 =
\iint_{S\times S}
K(\theta,\theta')\,d\nu(\theta)\,d\nu(\theta').
$$
Consequently, if $G$ and $H$ are probability measures on $S$, then
the Maximum Mean Discrepancy 
associated with $K$ between $G$ and $H$
is given by
$\|p_G-p_H\|_{L^2(\mathbb R)} = \|T(G-H)\|_{L^2(\mathbb{R})}.$
Thus identifiability of Gaussian mixture with mixing distributions concentrated on $S$ is equivalent to the restriction of $K$ to $S\times S$
being characteristic, in the terminology of reproducing kernel mean embeddings
\cite{sriperumbudur2010hilbert,sriperumbudur2011universality};
equivalently, in the present setting, it is integrally strictly positive
definite on finite signed measures supported on $S$.
Our results may therefore be read as geometric criteria describing precisely
when this nonstationary Gaussian probability-product kernel is characteristic
on a prescribed subset of the location--scale parameter space.

\bigskip
 
The rest of the paper is organized as follows. 
Section \ref{sec:prelim} provides some technical preliminaries.
Section \ref{sec:topo} presents sufficient topological conditions for identifiability, 
which are easily verifiable.
Section \ref{sec:geometric} presents a characterization of the null space $\mathcal N(S)$ based on the backward space-time Brownian motion and its exit behaviors. This theory leads to sufficient and necessary conditions for identifiability.
Section \ref{sec:5} applies the theory to shed light on identifiability of mixing measures supported on curves. It also contains a canonical mixing representation theorem on curved supports.
Appendix \ref{app:proofs} is devoted to proofs of some Lemmas and Propositions. 

A word of caution: in our exposition we use terms such as (location-scale) source, 
(Gaussian) output, residual output, and so on to provide some intuition and 
interpretation for the objects in discussion. Such terms do not appear 
in the rigorous definitions and theorems that follow.

\section{Preliminaries }\label{sec:prelim}

The Gaussian kernel $\phi(x; \mu, \sigma)$ has two basic functions: for a fixed location-scale pair $(\mu,\sigma) \in \mathbb{R} \times (0,\infty)$, it provides the probability density function for a normal random variable \footnote{For mathematical convenience, in this paper the normal variance is denoted by $\sigma$ as opposed to $\sigma^2$. } with mean (or location) $\mu$ and variance (or scale) $\sigma$. For a fixed $x \in \mathbb{R}$, it provides the Gaussian likelihood function
$\phi_x: (\mu,\sigma) \mapsto \phi_x(\mu, \sigma) = \phi(x; \mu, \sigma)$. Accordingly, in our study of identifiability for location-scale Gaussian mixtures, there are two distinct and deeply connected perspectives; one is probabilistic and the other functional analytic. They are both crucial in our formulation of necessary and sufficient conditions for identifiability, culminating in a precise characterization of Gaussian-null measures, namely, all elements that make up the null space $\mathcal{N}(S)$ of the Gaussian mixing operator. In the following we shall discuss these two perspectives along with several elementary facts, which will be useful for subsequent developments.

\subsection{Gaussian convolution and cancellation}
\label{sec:gaussian-cancellation}
Given a finite signed measure $\nu$ on a compact set $S \subset \mathbb{R} \times (0,\infty)$, the output of the Gaussian mixing operator, namely, $T\nu(\cdot) = \int \phi(\cdot;\mu,\sigma) \nu(d\mu,d\sigma)$, can be viewed as a weighted combination of the Gaussian output $\phi(\cdot; \mu,\sigma)$, with $\nu$ providing a description of the sources
of location-scale point $(\mu,\sigma) \in S$. The failure of identifiability, that is, the existence of a non-trivial $\nu$ such that $T\nu(x) = 0$ for all $x\in \mathbb{R}$, corresponds to the possibility of cancellation of the Gaussian outputs under a suitable signed weighting of the sources. Convolution provides the basic mechanism for such cancellation phenomena.

At a fixed scale $s>0$, the Gaussian output of a
finite signed location measure $\eta$ is a convolution:
$$
T(\eta\otimes\delta_s)(x)
=\int_{\mathbb R}\phi(x;\mu,s)\,\eta(d\mu)
=\bigl(\phi(\,\cdot\,;0,s)*\eta\bigr)(x),
$$
and its Fourier transform
$$
\widehat{T(\eta\otimes\delta_s)}(\xi)
=e^{-s\xi^2/2}\widehat\eta(\xi),
\text{ where }\widehat\eta(\xi):=\int_{\mathbb R}e^{-i\mu\xi}\,\eta(d\mu).
$$
Since the Gaussian multiplier never vanishes, uniqueness
of the Fourier transform gives
$
T(\eta\otimes\delta_s)=0 \Longrightarrow \eta=0.
$
Thus, for each fixed scale, Gaussian location mixtures
are identifiable without any restriction on the location
support. In particular, a nonzero Gaussian-null measure, i.e.,
a nontrivial member of $\mathcal{N}(S)$ if it exists,
cannot be supported on a single scale level.

Across different scale levels, the Gaussian convolution
identity again provides 
an elementary mechanism for redistributing
sources while preserving their Gaussian output. For $s,t>0$,
\begin{align*}
\phi(x;0,s+t)
=
\int_{\mathbb R}
\phi(x;\mu,s)\phi(\mu;0,t)\,d\mu=
\mathbb E\bigl[\phi(x;B_t,s)\bigr],
\end{align*}
where $B$ denotes a standard Brownian motion.
Thus a source at $(0,s+t)$ can be replaced by a Gaussian
distribution of locations at a lower scale $s$.
Their difference,
$
\widetilde\nu
=
\delta_{(0,s+t)}
-
\bigl(\phi(\mu;0,t)\,d\mu\bigr)\otimes\delta_s,
$
is a nonzero signed measure satisfying $T\widetilde\nu=0$.
Probabilistically, this replacement evolves the location
along Brownian paths for time $t$, while decreasing the
scale by the same amount, $(s+t)-t=s$.
The replacement measure here has unbounded location support,
so this particular construction lies outside our compact
setting. It nevertheless raises the central 
question: which sets $S$ admit distinct source configurations
supported on $S$ that produce the same Gaussian output?
As we will see in
Section \ref{sec:geometric}, the basic mechanism expressed by
Gaussian convolution leads to instances of Gaussian output 
cancellation, which amounts to conditions for Gaussian mixtures that fail to
be identifiable. To develop sufficient conditions for identifiability,
we need another perspective of Gaussian mixtures, where
the Gaussian likelihood functions play the role of test functions
for the source distribution.



\subsection{Gaussian likelihood as test functions}
\label{sec:gaussian-likelihood}


The collection of Gaussian likelihood functions $\{\phi_x: S \subset \mathbb R\times (0,\infty) \rightarrow \mathbb{R} \;|\; x \in \mathbb R\}$ may
be viewed as a class of test functions for evaluating the finite measure $\nu \in \mathcal{M}(S)$.
In particular, $\nu \in \mathcal{N}(S)$ means that $\langle v, \phi_x \rangle := \int \phi_x d \nu = 0$ for 
all test functions $\phi_x$ as $x$ varies in $\mathbb R$. As $S \subset \mathbb R\times (0,\infty)$, the finite measure $\nu$ can be uniquely identified by $C_b(S)$, the collection of continuous and bounded functions on $S$. As implied by the Gaussian cancellation phenomena described above, without restriction on the set $S$ the 
class of Gaussian 
likelihoods $\phi_x$ would be too small to ensure the identification of $\nu$. Thus, as we will see
shortly in Proposition \ref{thm:iden-characterization}, an identifiability condition on the set $S$ may be seen equivalently as a restriction according
to which the linear span of the class of Gaussian likelihood functions under the uniform norm is dense in $C_b(S)$.
Since we work with compact $S$, $C_b(S) = C(S)$, the space of continuous functions on $S$.

A key feature about $\phi_x$ is that it provides the fundamental solution for the heat equation
$\Lom f = 0$ for $f\in C^\infty(\mathbb{R}^2)$, where $\Lom=\partial_\sigma - \frac12 \partial^2_{\mu\mu}$
denotes the heat operator. Another useful fact is that $\phi_x$ may be approximated by the collection
of two-variable Hermite polynomials, which we now define.
For $(\mu,\sigma)\in \mathbb R\times(0,\infty)$, 
expanding in powers of $t$
the Fourier
transform $
\widehat{\phi}(t;\mu,\sigma)
=
e^{it\mu-\frac12\sigma t^2}
$ of Gaussian
kernel in the variable $x$, by analytic continuation we obtain the identity
\begin{equation}\label{eq-heat-poly}
	e^{z\mu+\frac12\sigma z^2}
	=
	\sum_{n=0}^{\infty}
	H_n(\mu,\sigma)\frac{z^n}{n!},\quad z\in \mathbb C,
\end{equation}
where $H_n$ are two-variable Hermite polynomials given by
$$
	H_n(\mu,\sigma)
	=
	n!
	\sum_{j=0}^{\lfloor n/2\rfloor}
	\frac{\mu^{\,n-2j}\sigma^j}
	{2^j j!(n-2j)!}.
	$$
It can be verified that $\Lom H_n = 0$ for all $n \geq 0$, and henceforth the $H_n$ are also referred to as heat polynomials. The approximation property mentioned earlier is formalized as follows. 

\begin{prop}\label{prop:equi-Gauss-Hermite}
Let $$\mathcal V_\phi:=\overline{
			\operatorname{span}
			\{\phi_x:x\in\mathbb{R}\}}\qquad\text{and}\qquad
    \mathcal V_H:=\overline{
		\operatorname{span}\{H_n|_S:n\geq0\}
	},
    $$
where the closure is taken under sup-norm in $C(S)$. Then, $\mathcal V_\phi=\mathcal V_H$.
\end{prop}

The following proposition makes precise the equivalence of identifiability for infinite location-scale Gaussian mixtures and an approximation property of Gaussian likelihood functions $\{(\mu,\sigma) \mapsto \phi_x(\mu,\sigma) | x\in\mathbb R\}$ and heat polynomials $\{H_n(\mu,\sigma)|_S:n\geq 0\}$ in the space of continuous functions on $S$.

\begin{prop}[]\label{thm:iden-characterization}
	Let
	$
	S\subset \mathbb{R}\times(0,\infty)
	$
	be compact. 
	The following statements are equivalent:
	\begin{itemize}
		\item[(i)] 
		$\mathcal P(S)$ is identifiable.
            \item[(ii)] The null space $\mathcal N(S)=\{0\}.$
		\item[(iii)] 
        $\mathcal V_\phi=
		C(S)
		$, i.e.,
       the linear span of the class of Gaussian
		likehood functions $\{\phi_x:x\in\mathbb{R}\}$ is dense in $C(S)$
		with respect to the uniform norm.
    \item[(iv)]
    $\mathcal V_H
	= C(S)$, i.e.,
   the linear span of the class of heat polynomials restricted on $S$, $\{H_n|_{S}: n\geq 0\}$ is dense in $C(S)$
		with respect to the uniform norm.
	\end{itemize}
\end{prop}



\begin{rem}
The proofs of Proposition \ref{prop:equi-Gauss-Hermite} and Proposition \ref{thm:iden-characterization} 
can be found in Appendix \ref{sec:proof-Sec2-Appendix}. The equivalences of (i), (ii) and (iii) in Proposition \ref{thm:iden-characterization} hold generally for any likelihood family arising from a valid probability kernel, not just the Gaussian kernel.
    The connection between Gaussian kernel and heat polynomials stated in Proposition \ref{prop:equi-Gauss-Hermite} leads to the equivalence of (iv).
\end{rem}

We end this preparation by noting that both functional analytic and probabilistic viewpoints
can be utilized to formulate separately sufficient or necessary conditions for identifiability, which 
amount to topological or geometric conditions on the set $S$ that make the null space $\mathcal{N}(S)$ trivial. Such results will be developed in
Section \ref{sec:topo}. To obtain sufficient \emph{and} necessary conditions, we will establish a
complete characterization of the null space $\mathcal{N}(S)$, which includes a precise representation 
of the members of the null space. It is remarkable that both probabilistic and functional analytic viewpoints 
are deeply integrated in this theory, which is presented in Section \ref{sec:geometric}.

\section{Sufficient conditions of identifiability}\label{sec:topo}


We shall first present a sufficient condition for the null space $\mathcal{N}(S)$ of the Gaussian mixing operator to be trivial. This condition is formulated by exposing the impossibility of cancellation of Gaussian outputs $\{\phi(\cdot| \mu,\sigma)|(\mu,\sigma) \in S\}$, provided that the scale projection of the set of sources $S$ has an empty interior. The benefit of this rather simple condition is that it has an intuitive and strong connection to the identifiability of convolution Gaussian mixtures, i.e., for each fixed scale level, where Fourier analysis plays an useful role (cf. Section \ref{sec:gaussian-cancellation}). This result is presented in Section \ref{sec:topo-thin}. 
Section \ref{sec:heat adjoint} embarks on the general regime where such a condition on the scale projection may be violated. One needs to anticipate cancellation across the possibly connected location-scale domain $S$, and
at this point make use of the fact that the two-variable Gaussian likelihood functions $\phi_x$ are the fundamental solution to the heat equation (cf. Section \ref{sec:gaussian-likelihood}). Building on the preliminaries presented in the previous section, it will be established that any element in the null space $\mathcal{N}(S)$ for a compact $S$ can be represented in terms of a compactly supported and square integrable adjoint heat potential. Not only does this representation allow us to arrive at relatively simple sufficient topological conditions on the support $S$ for which the null set is trivial, i.e., the Gaussian mixtures on $\mathcal{S}$ is identifiable (see Section \ref{sec:sufficient-adjoint}), it also provides the foundation for subsequent developments in Section \ref{sec:geometric}.





 \subsection{A condition on topological thinness of $\sigma$-projection}\label{sec:topo-thin}

We begin by examining whether sources on opposite sides of a positive scale gap can compensate for one another. The following basic lemma shows that, when two sets of location-scale sources are separated by a positive scale gap, the cancellation of their combined Gaussian outputs forces the output of each set of sources to vanish separately. 
Thus, cancellation must occur within each set, without compensation from sources on the other side of the scale gap.



\begin{lem}[Scale-gap decoupling]
	\label{lem:scale-gap}
	Let $S_-,S_+\subset \mathbb R\times(0,\infty)$ be compact. Suppose that, for some
	$
	0<a<b,
	$
	one has
	$$
	S_-\subset\{(\mu,\sigma):\sigma\le a\},
	\qquad
	S_+\subset\{(\mu,\sigma):\sigma\ge b\}.
	$$
	Let $\nu_-$ and $\nu_+$ be finite signed measures supported on $S_-$ and $S_+$, respectively. If
	$$
	(T\nu_-)(x)+(T\nu_+)(x)=0,
	\qquad
	\text{for every }x\in\mathbb R,
	$$
	then
	$$
	T\nu_-=T\nu_+=0.
	$$
\end{lem}

\begin{proof}
	Set
	$$
	f(x):=(T\nu_-)(x)=-(T\nu_+)(x).
	$$
	Since $S_-\cup S_+\subset \mathbb R\times(0,\infty)$ is compact,
	there exist constants
	$
	M<\infty,
	\;
	\underline\sigma>0,
	$
	such that 	
	$$
	|\mu|\le M,
	\;\;
\sigma\ge\underline\sigma\quad\text{for all}\quad(\mu,\sigma)\in S_-\cup S_+.
	$$

	Using the representation of $f$ through $\nu_-$, with $\|\cdot\|_{TV}$ denoting the total variation,
    we obtain
	\begin{align*}
		|f(x)|
		&\le
		\int_{S_-}
		\frac{1}{\sqrt{2\pi\sigma}}
		\exp\left(
		-\frac{(x-\mu)^2}{2\sigma}
		\right)
		\,d|\nu_-|(\mu,\sigma)
		\le
		\frac{\|\nu_-\|_{\mathrm{TV}}}
		{\sqrt{2\pi\underline\sigma}}
		\exp\left(
		-\frac{(|x|-M)^2}{2a}
		\right).
	\end{align*}
	Choose $a' \in (a,b)$.
	Since
	$$
	\frac{x^2}{2a'}
	-
	\frac{(|x|-M)^2}{2a}
	$$
	is bounded above on $\mathbb R$, there exists $C_1<\infty$ such that
	\begin{equation}
		\label{eq:gap-x-decay}
		|f(x)|
		\le
		C_1
		\exp\left(
		-\frac{x^2}{2a'}
		\right),
		\qquad
		\forall x\in\mathbb R.
	\end{equation}
	
	We next use the representation of $f$ through $\nu_+$. 
We have the following bound for its Fourier transform $\hat f(t)$
	\begin{equation}
		\label{eq:gap-f-decay}
		\widehat f(t)
	=
	-\int_{S_+}
	e^{it\mu-\sigma t^2/2}
	\,d\nu_+(\mu,\sigma)\quad\text{and }\quad
        |\widehat f(t)|
		\le
		\|\nu_+\|_{\mathrm{TV}}
		\exp\left(
		-\frac{bt^2}{2}
		\right),\quad\forall t\in \mathbb R.
	\end{equation}
By Hardy's uncertainty principle, see, e.g.,
\cite[pg. 158]{stein2003fourier}
	it follows from \eqref{eq:gap-x-decay} and \eqref{eq:gap-f-decay} that
	$
	f\equiv0.
	$
	Consequently,
	$
	T\nu_-=T\nu_+=0.
	$
\end{proof}

From 
this lemma, we have the following result.
\begin{thm}
	\label{thm:thin-scale-proj}
	Let
	$
	S\subset \mathbb R\times(0,\infty)
	$
	be compact, and denote
	$$
	\pi_\sigma(S)
	:=
	\bigl\{
	\sigma>0:\ (\mu,\sigma)\in S
	\text{ for some }\mu\in\mathbb R
	\bigr\}.
	$$
	If
	$
	\operatorname{Int}\pi_\sigma(S)=\emptyset,
	$
then $\mathcal P(S)$ is identifiable.
\end{thm}

\begin{proof}
	Since $S$ is compact, $\pi_\sigma(S)$ is a compact subset of $(0,\infty)$.
Let $\nu\in\mathcal N(S)$ be arbitrary but fixed. Then $(T\nu)(x) = 0$ for all $x\in \mathbb R$. We need to show that $\nu = 0$. The proof proceeds in several steps as follows.

\medskip

\noindent{\bf Step 1.} (Localization of cancellation across scale gaps).
	Let
	$
	c\in\mathbb R\setminus \pi_\sigma(S),
	$
	and define
	$$
	\nu_c^-
	:=
	\mathbf 1_{\{\sigma<c\}}\nu,
	\qquad
	\nu_c^+
	:=
	\mathbf 1_{\{\sigma>c\}}\nu.
	$$
	We first show that $\nu_c^+,\nu_c^-\in\mathcal N(S)$, i.e., 
		\begin{equation}
		\label{eq:cut-null}
		(T\nu_c^-)(x)
		=
		(T\nu_c^+)(x)
		=
		0,\qquad\forall x\in\mathbb R.
	\end{equation}
    Indeed,	since $c\notin \pi_\sigma(S)$,
	$
	\nu=\nu_c^-+\nu_c^+
$. 	 
If $\nu_c^+$ or $\nu_c^-$ vanishes, \eqref{eq:cut-null} follows directly. Otherwise, since $\pi_\sigma(S)$ is compact and $c\notin \pi_\sigma(S)$, there exist
	$
	a<c<b
	$
	such that
	$
	\operatorname{supp}\nu_c^-
	\subset
	\{\sigma\le a\},
	\;
	\operatorname{supp}\nu_c^+
	\subset
	\{\sigma\ge b\}.
	$
	Since $T\nu= 0$,
    $
	(T\nu_c^-)(x)+(T\nu_c^+)(x)=0,\;\forall x\in\mathbb R.
	$
	Lemma \ref{lem:scale-gap} therefore yields
	$
		(T\nu_c^-)(x)
		=
		(T\nu_c^+)(x)
		=
		0,\;\forall x\in\mathbb R.
$
	

\medskip
    
\noindent{\bf Step 2.} (Abundance of gaps and application of Stone-Weierstrass theorem). Since $c\notin \pi_\sigma(S)$, the function $\mathbf 1_{\{\sigma<c\}}$ is continuous on $\pi_\sigma(S)$.
	Consider the class of functions
	$$
	\mathcal A
	:=
	\operatorname{span}
	\left(
	\{1\}
	\cup
	\{\mathbf 1_{\{\sigma<c\}}:c\notin \pi_\sigma(S)\}
	\right)
	\subset C(\pi_\sigma(S)).
	$$
	Since
	$
	\mathbf 1_{\{\sigma<c\}}\mathbf 1_{\{\sigma<d\}}=\mathbf 1_{\{\sigma<\min\{c,d\}\}}
	$,
	$\mathcal A$ is an algebra containing constant functions.
	We claim that $\mathcal A$ separates the points of $\pi_\sigma(S)$. Let
	$
	\sigma_1,\sigma_2\in \pi_\sigma(S),
	\;
	\sigma_1<\sigma_2.
	$
	Since
	$
	\operatorname{Int}\pi_\sigma(S)=\emptyset,
	$
	the interval $(\sigma_1,\sigma_2)$ cannot be contained in $\pi_\sigma(S)$. Hence there exists
	$
	c\in(\sigma_1,\sigma_2)\setminus \pi_\sigma(S).
	$
	For this $c$,
	$
	\mathbf 1_{\{\sigma_1<c\}}=1
	$
    while 
    $
	\mathbf 1_{\{\sigma_2<c\}}=0.
	$
	Thus $\mathcal A$ separates points. By the Stone--Weierstrass theorem,
	\begin{equation}
		\label{eq:sw}
\mathcal A	\text{ is dense in }
		C(\pi_\sigma(S)).
	\end{equation}
	
	For $h\in C(\pi_\sigma(S))$, define
	$$
	T(\nu,h)(x)
	:=
	\int_S
	h(\sigma)\phi(x;\mu,\sigma)
	\,d\nu(\mu,\sigma).
	$$
We have
	\begin{align*}
		\|T(\nu,h)\|_\infty
		&\le
		\|h\|_\infty
		\|\nu\|_{\mathrm{TV}}
		\sup_{\substack{x\in\mathbb R\;(\mu,\sigma)\in S}}
		\phi(x;\mu,\sigma)\le
		\frac{\|\nu\|_{\mathrm{TV}}}
		{\sqrt{2\pi\underline\sigma}}
		\|h\|_\infty.
	\end{align*}
	Hence, for fixed $\nu$,
	$
	T(\nu,h):C(\pi_\sigma(S))\rightarrow C_b(\mathbb R)
	$
	is continuous in $h$.
	Since $T\nu = 0$, and from \eqref{eq:cut-null},
	$$
	T(\nu,1)=0\quad\text{ and }\quad
	T(\nu,\mathbf 1_{\{\sigma<c\}})=0
	\;
	\text{for every }c\notin \pi_\sigma(S).
	$$
	Consequently,
	$
	T(\nu,h)=0 
	$ 
    for all $h\in\mathcal A$.
	Using \eqref{eq:sw} and the continuity of $T(\nu,h)$ in $h$, we obtain
	\begin{equation}
		\label{eq:weighted-null}
		\int_S
		h(\sigma)
		\phi(x;\mu,\sigma)
		\,d\nu(\mu,\sigma)
		=
		0,\quad\forall h\in C(\pi_\sigma(S)),\;	\forall x\in\mathbb R.
	\end{equation}

\medskip

\noindent{\bf Step 3.} (Application of Fourier inversion theorem).	We now recover from \eqref{eq:weighted-null} the full two-dimensional Fourier transform of $\nu$.
	Taking the Fourier transform of \eqref{eq:weighted-null} with respect to $x$ gives
	\begin{equation}
		\label{eq:weighted-f-null}
		\int_S
		h(\sigma)
		e^{it\mu-\sigma t^2/2}
		\,d\nu(\mu,\sigma)
		=
		0,\quad\text{for every }
		h\in C(\pi_\sigma(S)),
		\;
		t\in\mathbb R.
	\end{equation}
	The interchange of integrals is justified since
	$
	\int_S
	\int_{\mathbb R}
	|h(\sigma)|
	\phi(x;\mu,\sigma)
	\,dx\,
	d|\nu|(\mu,\sigma)
	\le
	\|h\|_\infty\|\nu\|_{\mathrm{TV}}
	<
	\infty.
	$
	
	Fix arbitrary
	$
	t_1,t_2\in\mathbb R
	$, applying \eqref{eq:weighted-f-null} with
	$$
	h_1(\sigma)
	=
	e^{\sigma t_1^2/2}\cos(t_2\sigma)\qquad 
	\text
	{and}
	\qquad
	h_2(\sigma)
	=
	e^{\sigma t_1^2/2}\sin(t_2\sigma)
	$$
gives
	$$
	\int_S
	e^{it_1\mu}\cos(t_2\sigma)
	\,d\nu(\mu,\sigma)
	=
	0\qquad 
	\text{and}
	\qquad
	\int_S
	e^{it_1\mu}\sin(t_2\sigma)
	\,d\nu(\mu,\sigma)
	=
	0.
	$$
	Combining the two identities yields
	\begin{equation}
		\label{eq:full-f-zero}
		\int_S
		e^{i(t_1\mu+t_2\sigma)}
		\,d\nu(\mu,\sigma)
		=
		0
		\qquad
		\text{for every }(t_1,t_2)\in\mathbb R^2.
	\end{equation}
	
	The left-hand side of \eqref{eq:full-f-zero} is the Fourier transform of the finite signed measure $\nu$ on $\mathbb R^2$. By uniqueness of the Fourier transform of finite measures,
	$
	\nu=0.
	$
    This concludes the proof.
\end{proof}

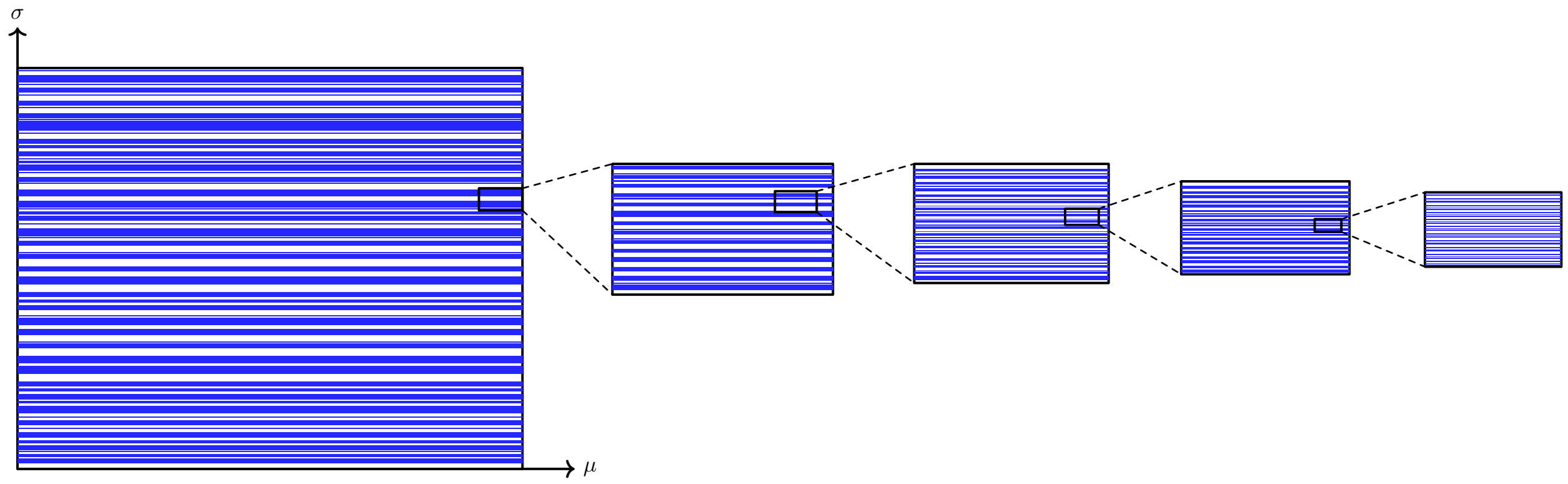
\begin{figure}[H]
\centering
\resizebox{\textwidth}{!}{%
\begin{tikzpicture}[x=1cm,y=1cm,line cap=round,line join=round]

\tikzset{
    axis/.style={->, very thick},
    frame/.style={draw=black, very thick},
    zoombox/.style={draw=black, very thick},
    zoomline/.style={dashed, thick},
    blueband/.style={draw=none, fill=blue!85}
}

%
\newcommand{\BandPattern}[4]{%

    \foreach \yy in {0.10,0.32,0.54,0.76,0.98,1.20,1.42,1.64,1.86,2.08,2.30,2.52,2.74,2.96,3.18,3.40,3.62,3.84,4.06,4.28,4.50,4.72,4.94,5.16,5.38,5.60,5.82,6.04,6.26,6.48,6.70,6.92} {
        \fill[blueband] (#1,#2+\yy) rectangle (#1+#3,#2+\yy+0.08);
    }

    \foreach \yy in {0.20,0.58,0.96,1.34,1.72,2.10,2.48,2.86,3.24,3.62,4.00,4.38,4.76,5.14,5.52,5.90,6.28,6.66} {
        \fill[blueband] (#1,#2+\yy) rectangle (#1+#3,#2+\yy+0.05);
    }

    \foreach \yy in {0.44,1.13,1.82,2.51,3.20,3.89,4.58,5.27,5.96,6.65} {
        \fill[blueband] (#1,#2+\yy) rectangle (#1+#3,#2+\yy+0.035);
    }

    \foreach \yy in {0.14,0.29,0.47,0.69,0.88,1.07,1.27,1.48,1.68,1.89,2.17,2.39,2.63,2.79,3.01,3.29,3.47,3.71,3.97,4.21,4.47,4.69,4.91,5.09,5.33,5.55,5.77,6.01,6.21,6.43,6.61,6.85} {
        \fill[blueband] (#1,#2+\yy) rectangle (#1+#3,#2+\yy+0.018);
    }
}

\draw[axis] (0,0) -- (0,7.6) node[above] {$\sigma$};
\draw[axis] (0,0) -- (9.6,0) node[right] {$\mu$};

\def\xA{0.0}
\def\yA{0.0}
\def\wA{8.7}
\def\hA{6.9}

\draw[frame, fill=white] (\xA,\yA) rectangle (\xA+\wA,\yA+\hA);
\begin{scope}
    \clip (\xA,\yA) rectangle (\xA+\wA,\yA+\hA);
    \BandPattern{\xA}{\yA}{\wA}{\hA}
\end{scope}

\def\zxA{7.95}
\def\zyA{4.45}
\def\zwA{0.75}
\def\zhA{0.38}
\draw[zoombox] (\zxA,\zyA) rectangle (\zxA+\zwA,\zyA+\zhA);

\def\xB{10.25}
\def\yB{3.00}
\def\wB{3.80}
\def\hB{2.25}

\draw[frame, fill=white] (\xB,\yB) rectangle (\xB+\wB,\yB+\hB);
\begin{scope}
    \clip (\xB,\yB) rectangle (\xB+\wB,\yB+\hB);

    \foreach \yy in {0.08,0.24,0.40,0.56,0.72,0.88,1.04,1.20,1.36,1.52,1.68,1.84,2.00,2.16} {
        \fill[blueband] (\xB,\yB+\yy) rectangle (\xB+\wB,\yB+\yy+0.060);
    }
    \foreach \yy in {0.14,0.44,0.74,1.04,1.34,1.64,1.94} {
        \fill[blueband] (\xB,\yB+\yy) rectangle (\xB+\wB,\yB+\yy+0.030);
    }
    \foreach \yy in {0.19,0.31,0.63,0.95,1.11,1.43,1.71,1.87,2.07} {
        \fill[blueband] (\xB,\yB+\yy) rectangle (\xB+\wB,\yB+\yy+0.015);
    }
\end{scope}

\draw[zoomline] (\zxA+\zwA,\zyA+\zhA) -- (\xB,\yB+\hB);
\draw[zoomline] (\zxA+\zwA,\zyA)       -- (\xB,\yB);

\def\zxB{13.05}
\def\zyB{4.42}
\def\zwB{0.72}
\def\zhB{0.36}
\draw[zoombox] (\zxB,\zyB) rectangle (\zxB+\zwB,\zyB+\zhB);

\def\xC{15.45}
\def\yC{3.20}
\def\wC{3.35}
\def\hC{2.05}

\draw[frame, fill=white] (\xC,\yC) rectangle (\xC+\wC,\yC+\hC);
\begin{scope}
    \clip (\xC,\yC) rectangle (\xC+\wC,\yC+\hC);

    \foreach \yy in {0.05,0.16,0.27,0.38,0.49,0.60,0.71,0.82,0.93,1.04,1.15,1.26,1.37,1.48,1.59,1.70,1.81,1.92} {
        \fill[blueband] (\xC,\yC+\yy) rectangle (\xC+\wC,\yC+\yy+0.040);
    }
    \foreach \yy in {0.11,0.33,0.55,0.77,0.99,1.21,1.43,1.65,1.87} {
        \fill[blueband] (\xC,\yC+\yy) rectangle (\xC+\wC,\yC+\yy+0.018);
    }
    \foreach \yy in {0.09,0.21,0.42,0.67,0.88,1.09,1.32,1.58,1.79} {
        \fill[blueband] (\xC,\yC+\yy) rectangle (\xC+\wC,\yC+\yy+0.010);
    }
\end{scope}

\draw[zoomline] (\zxB+\zwB,\zyB+\zhB) -- (\xC,\yC+\hC);
\draw[zoomline] (\zxB+\zwB,\zyB)       -- (\xC,\yC);

\def\zxC{18.05}
\def\zyC{4.20}
\def\zwC{0.58}
\def\zhC{0.28}
\draw[zoombox] (\zxC,\zyC) rectangle (\zxC+\zwC,\zyC+\zhC);

\def\xD{20.05}
\def\yD{3.35}
\def\wD{2.90}
\def\hD{1.60}

\draw[frame, fill=white] (\xD,\yD) rectangle (\xD+\wD,\yD+\hD);
\begin{scope}
    \clip (\xD,\yD) rectangle (\xD+\wD,\yD+\hD);

    \foreach \yy in {0.03,0.11,0.19,0.27,0.35,0.43,0.51,0.59,0.67,0.75,0.83,0.91,0.99,1.07,1.15,1.23,1.31,1.39,1.47} {
        \fill[blueband] (\xD,\yD+\yy) rectangle (\xD+\wD,\yD+\yy+0.028);
    }
    \foreach \yy in {0.07,0.23,0.39,0.55,0.71,0.87,1.03,1.19,1.35,1.51} {
        \fill[blueband] (\xD,\yD+\yy) rectangle (\xD+\wD,\yD+\yy+0.012);
    }
    \foreach \yy in {0.14,0.30,0.46,0.62,0.78,0.94,1.10,1.26,1.42} {
        \fill[blueband] (\xD,\yD+\yy) rectangle (\xD+\wD,\yD+\yy+0.007);
    }
\end{scope}

\draw[zoomline] (\zxC+\zwC,\zyC+\zhC) -- (\xD,\yD+\hD);
\draw[zoomline] (\zxC+\zwC,\zyC)       -- (\xD,\yD);

\def\zxD{22.35}
\def\zyD{4.08}
\def\zwD{0.46}
\def\zhD{0.22}
\draw[zoombox] (\zxD,\zyD) rectangle (\zxD+\zwD,\zyD+\zhD);

\def\xE{24.25}
\def\yE{3.48}
\def\wE{2.35}
\def\hE{1.28}

\draw[frame, fill=white] (\xE,\yE) rectangle (\xE+\wE,\yE+\hE);
\begin{scope}
    \clip (\xE,\yE) rectangle (\xE+\wE,\yE+\hE);

    \foreach \yy in {0.02,0.08,0.14,0.20,0.26,0.32,0.38,0.44,0.50,0.56,0.62,0.68,0.74,0.80,0.86,0.92,0.98,1.04,1.10,1.16,1.22} {
        \fill[blueband] (\xE,\yE+\yy) rectangle (\xE+\wE,\yE+\yy+0.020);
    }
    \foreach \yy in {0.05,0.17,0.29,0.41,0.53,0.65,0.77,0.89,1.01,1.13} {
        \fill[blueband] (\xE,\yE+\yy) rectangle (\xE+\wE,\yE+\yy+0.009);
    }
\end{scope}

\draw[zoomline] (\zxD+\zwD,\zyD+\zhD) -- (\xE,\yE+\hE);
\draw[zoomline] (\zxD+\zwD,\zyD)       -- (\xE,\yE);

\end{tikzpicture}%
}
\caption{Schematic illustration of $S=[0,1]\times K_N$: although $S$ has large measure but $\mathcal P(S)$ is still identifiable because the projection in the $\sigma$-direction has empty interior, suggested by the increasingly fine interlacing of blue bands and white gaps under successive zooms.}
\label{fig:SN-zoom}
\end{figure}

\begin{rem}\label{rem-S-large}
It is shown in Proposition \ref{thm-ness-cond} that $S$ must be topological thin to be identifiable.
However, by Theorem \ref{thm:thin-scale-proj}, we illustrate that $S$ can have arbitrarily large Lebesgue measure while $\mathcal P(S)$ is still identifiable.

We first claim that:
for every integer $N \ge 1$, there exists a compact set
$K_N \subset [1,2]$ such that
$$
\operatorname{Int}(K_N)=\emptyset
\qquad\text{and}\qquad
\text{one-dimensional Lebesgue measure of } K_N>1-1/N.
$$
Indeed, enumerate the rational numbers in $[1,2]$ as
$\mathbb{Q}\cap[1,2]=\{q_1,q_2,\ldots\}.$
Define
$$
U_N
=
\bigcup_{n=1}^{\infty}
\left(
q_n-\frac{1}{N2^{n+2}},
q_n+\frac{1}{N2^{n+2}}
\right),
\;
K_N=[1,2]\setminus U_N.
$$
Since $U_N$ is open, $K_N$ is a closed subset of the compact
interval $[1,2]$, and hence is compact. Moreover,
$K_N\cap\mathbb{Q}=\emptyset$. Thus,
$\operatorname{Int}(K_N)=\emptyset.$
By countable subadditivity,
 Lebesgue measure of
$\bigl(U_N\cap[1,2]\bigr)
\le
\sum_{n=1}^{\infty}\frac{2}{N2^{n+2}}
=
\frac{1}{2N}.
$
Consequently,
 Lebesgue measure of
$(K_N)
\ge
1-\frac{1}{2N}
>
1-\frac{1}{N}.
$

By this idea, 
we can easily construct a compact set $S$ having arbitrarily large Lebesgue measure and satisfying conditions in Theorem \ref{thm:thin-scale-proj} so that $\mathcal P(S)$ is identifiable. 
See Figure \ref{fig:SN-zoom}.
\end{rem}

A simple consequence of Theorem \ref{thm:thin-scale-proj} is the following result.

\begin{cor}\label{thm-countable-iden}
    Let $S\subset\mathbb R\times(0,\infty)$ be compact and countable.
    Then the space of infinite countably Gaussian mixing distribution in $S$ is identifiable.
\end{cor}

Another application of scale-gap argument is the following reduction result.
\begin{prop}
	\label{cor:var-gap}
	Let
	$
	S=S_-\cup S_+
	$
	where $S_-$ and $S_+$ are compact subsets of
	$\mathbb R\times(0,\infty)$ satisfying
	$
	\max\pi_\sigma(S_-)
	<
	\min\pi_\sigma(S_+).
	$
	Then
	$$
	\mathcal P(S)\text{ is identifiable}
	\quad
    \Longleftrightarrow
	\quad
    \text{both }
	\mathcal P(S_-)
	\text{ and }
	\mathcal P(S_+)
	\text{ are identifiable}.
    $$
\end{prop}

\begin{proof}
	If $\mathcal P(S)$ is identifiable, then every subclass obtained by
	restricting the support is identifiable, so the forward implication
	is immediate.
    Conversely, suppose that $\mathcal P(S_-)$ and $\mathcal P(S_+)$
	are identifiable. 
    Let
 $\nu\in\mathcal N(S)$.
	Write
	$
	\nu=\nu_-+\nu_+,
	$
	where $\nu_\pm$ denotes the restriction of $\nu$ to $S_\pm$.
	Then
	$
	(T\nu_-)(x)+(T\nu_+)(x)=0,\;\forall x\in\mathbb R.
	$
	By Lemma~\ref{lem:scale-gap},
	$
	(T\nu_-)(x)=(T\nu_+)(x)=0\;\forall x\in\mathbb R.
	$
	The assumed identifiability on $S_-$ and $S_+$ gives
	$
	\nu_-=\nu_+=0.
	$
	Thus $\nu=0$, and the conclusion follows.
\end{proof}

\begin{thm}[Reduction to nondegenerate scale components]
	\label{thm:var-reduction}
	Let
	$
	S\subset\mathbb R\times(0,\infty)
	$
	be compact.
	Let $\mathscr I(\pi_\sigma(S))$ denote the collection of all nondegenerate
	connected components of $\pi_\sigma(S)$. 
	For each
	$
	I\in\mathscr I(\pi_\sigma(S)),
	$
	define
	$$
	S_I
	:=
	S\cap(\mathbb R\times I).
	$$
	Then
	\begin{equation}
		\label{eq:-reduction-equi}
		\mathcal P(S)\text{ is identifiable}
		\quad\Longleftrightarrow\quad
		\mathcal P(S_I)\text{ is identifiable for every }
		I\in\mathscr I(\pi_\sigma(S)).
	\end{equation}
\end{thm}
\begin{proof}
	The forward implication is immediate, since every probability measure
	supported on $S_I$ is also supported on $S$.
	
	For the converse, suppose that
	$
	\mathcal P(S_I)
	$
	is identifiable for every
	$
	I\in\mathscr I(\pi_\sigma(S))
	$, 
    and that $\nu \in \mathcal{N}(S)$.
	Fix
	$
	I\in\mathscr I(\pi_\sigma(S))
	$, which
	has the form
	$
	I=[a,b],\;
	a< b.
	$
    	There exist monotone
	sequences
	$
	\ell_n\uparrow a,
	\;
	r_n\downarrow b,
	$
	such that
	$
	\ell_n,r_n\notin \pi_\sigma(S).
	$
    	For every
	$
	c\notin \pi_\sigma(S),
	$
	the scale-gap argument \eqref{eq:cut-null} in the proof of
	Lemma \ref{lem:scale-gap} shows that
	$
	\mathbf 1_{\{\sigma<c\}}\nu\in\mathcal N(S).
	$
Let
	$
	\nu_n
	:=
	\mathbf 1_{\{\ell_n<\sigma<r_n\}}\nu.
	$
  Since
	$
	\nu_n
	=
	\mathbf 1_{\{\sigma<r_n\}}\nu
	-
	\mathbf 1_{\{\sigma<\ell_n\}}\nu
	$,
    $\nu_n\in\mathcal N(S).$

    	The sets
	$
	S\cap\mathbb R\times (\ell_n,r_n)
	$
	decrease to $S_I$. Since $|\nu|$ is finite,
	continuity from above gives
	\begin{equation}
		\label{eq:comp-tv-conv}
		\left\|
		\nu_n-\nu|_{S_I}
		\right\|_{\mathrm{TV}}
		\longrightarrow0\text{ as }n\to\infty.
	\end{equation}
Let
	$
	\underline\sigma:=\min\pi_\sigma(S)>0.
	$
Then,
	\begin{align*}
		\sup_{x\in\mathbb R}
		\left|
		(T\nu_n)(x)-(T(\nu|_{S_I}))(x)
		\right|
		\le
		\sup_{
        \mathbb R\times S
        }
		\phi(x;\mu,\sigma)
		\left\|
		\nu_n-\nu|_{S_I}
		\right\|_{\mathrm{TV}}\le
		\frac{1}{\sqrt{2\pi\underline\sigma}}
		\left\|
		\nu_n-\nu|_{S_I}
		\right\|_{\mathrm{TV}}.
	\end{align*}
	By \eqref{eq:comp-tv-conv},
	the right-hand side tends to zero as $n\to\infty$. Since
	$
	T\nu_n=0
	$
	for every $n$, it follows that
	$
	T\nu|_{S_I}=0.
	$
	By identifiability of $\mathcal P(S_I)$,
	$
	\nu|_{S_I}=0.
	$
	We have therefore shown that $\nu$ vanishes on every nondegenerate
	scale component.
	
	Set
	$
	K:=\operatorname{supp}|\nu|.
	$
	Then $K$ is a compact subset of $S$. We claim that
$\operatorname{Int}\pi_\sigma(K)=\emptyset.
	$	
	Suppose, to the contrary, that
$\operatorname{Int}\pi_\sigma(K)\ne\emptyset.$
	Then $\pi_\sigma(K)$ contains a nontrivial open interval $J$.
	Since
	$
	J\subset \pi_\sigma(S)
	$
	and $J$ is connected, it is contained in a connected component
	$
	I=[a,b]
	$
	of $\pi_\sigma(S)$. Because $J$ is nontrivial, this component is nondegenerate and
	therefore belongs to $\mathscr I(\pi_\sigma(S))$.
	Choose
	$
	\sigma_0\in J\cap(a,b).
	$
	Since
	$
	\sigma_0\in\pi_\sigma(K),
	$
	there exists $\mu_0$ such that
	$
	(\mu_0,\sigma_0)\in K.
	$
	Since $\sigma_0$ lies in the interior of $I$, there exists
	$\varepsilon>0$ such that
	$
	(\sigma_0-\varepsilon,\sigma_0+\varepsilon)
	\subset I.
	$
    But
	$
	\nu|_{S_I}=0
	$
	and so,
	$
	|\nu|
	\left(
	S\cap
	\left(
	\mathbb R\times
	(\sigma_0-\varepsilon,\sigma_0+\varepsilon)
	\right)
	\right)
	=0.
	$
	This contradicts
$(\mu_0,\sigma_0)\in\operatorname{supp}|\nu|.
	$
	Therefore
	$
\operatorname{Int}\pi_\sigma(K)=\emptyset.
	$
	
	By
	Theorem \ref{thm:thin-scale-proj},
	$\mathcal P(K)$ is identifiable. Since $\nu$ is supported on $K$ and
	$
	T\nu=0,
	$
	it follows that
	$
	\nu=0.
	$
	Hence $\mathcal P(S)$ is identifiable.
\end{proof}

\subsection{Null space characterization via adjoint heat potentials}\label{sec:heat adjoint}

The scale-gap results (Lemma \ref{lem:scale-gap} and Theorem \ref{thm:var-reduction}) show that for Gaussian-null measure $\nu$, which produces the Gaussian cancellation as $T\nu(x) =0$ for all $x\in \mathbb R$, source components separated by a positive scale gap must cancel separately. To understand this phenomenon more deeply, and to prepare for the general setting in which positive scale gaps are not necessarily available, we shall examine how the distribution $\nu$ of the Gaussian sources asserts its impacts on spatial locations at intermediate scale levels. 

Noting that a Gaussian output at a scale level may be cancelled by sources at lower scale levels, this suggests a notion of ``residual" contribution of the Gaussian output that remains at each intermediate scale level. Specifically,
at the scale level $r$, a source at $(\mu,\sigma)$ with $\sigma>r$ has ``diffused" through $\sigma-r$ units of scale. Summing these signed contributions leads to setting
\begin{equation}\label{eq-deff-u-explicit}
u(x,r)
:=
\int_{\{\sigma>r\}}
\phi(x;\mu,\sigma-r)\,\nu(d\mu,d\sigma).
\end{equation}
Thus $u(\cdot, r)$ records the residual profile at the scale level $r$, while at level zero
$u(\cdot,0) = T\nu$ collects the total Gaussian output attributed to the sources prescribed by the signed measure $\nu$. 

To assess precisely what the scale-gap argument detects at a cut when there
exists a scale gap in $\pi_{\sigma}(S)$,
consider 
 $r\in(a,b)$ in the scale gap.
By definition \eqref{eq-deff-u-explicit}, only the upper restriction
$\nu_{b_+}:=\nu|_{\{\sigma\ge b\}}$ contributes to the residual profile $u$ at levels $r$.
The Gaussian semigroup identity
$
\phi(x;\mu,\sigma) = 
\int_{\mathbb R}
\phi(x;y,r)\,\phi(y;\mu,\sigma-r)\,dy
$
and Fubini's theorem therefore give
$$
T\nu_{b_+}(x)
= \int \phi(x;\mu,\sigma) \nu_{b+}(d\mu, d\sigma) = 
\int_{\mathbb R}\phi(x;y,r)\,u(y,r)\,dy.
$$
Since
$u(\,\cdot\,,r)\in L^1(\mathbb R)$, injectivity of Gaussian
convolution yields
$$
T\nu_{b_+}=0
\quad\Longleftrightarrow\quad
u(\,\cdot\,,r)=0.
$$
The scale-gap argument of Lemma \ref{lem:scale-gap} asserts that if 
$T\nu =0$, then $T\nu_{b_+} = 0$. It follows that the upper scale level 
sources annihilate all Gaussian test functions exactly when their 
residual profile vanishes at the intermediate level $r$.
In other words, the upper scale level sources have already cancelled their entire 
response upon reaching the gap. 

Now, consider the case that there is no scale gap. In general, the residual profile $u$ records what remains at each intermediate
level, rather than only the final output $T\nu$. Even when the output $T\nu$ equals zero everywhere, the net contribution from sources at larger scale levels may remain nonzero at an intermediate level. 
Stating in the language
of residual profiles, under identifiability, whenever the final output
$u(\cdot, 0)$ vanishes the residual profiles must also
vanish at all intermediate scale levels as well.

It will be seen that the residual profile $u$ plays a crucial
role in our investigation of the Gaussian source $\nu$.
It may be verified from \eqref{eq-deff-u-explicit} that $u$ satisfies the partial differential
equation $-\partial_r u - \frac 12 \partial^2_{xx} u = \nu$, where the
identity is understood in the sense of distributions. In other words,
the residual profile $u$ may be identified as a solution to the adjoint
heat equation, given the source $\nu$, and will be referred henceforth as
the corresponding adjoint heat potential for $\nu$. More importantly, we shall
use this relation as the fundamental representation for the
members of the null space $\mathcal{N}(S)$.




Although $u$ admits the explicit representation
\eqref{eq-deff-u-explicit}, we shall use instead a structural argument to establish
that, since $\nu$ has compact support, the null condition $T\nu=0$ ensures there is a distribution $u\in L^2(\R^2)$ compactly
supported and satisfying the relation
$$
\check \Lom u=\nu.
$$
where
$\check \Lom$ is the adjoint of the heat operator 
$\Lom=\partial_\sigma-\frac 12\partial^2_{\mu\mu}$.

\medskip

This result is given in the following Proposition, which, in summary, is motivated from the 
observation that (i) the Gaussian likelihood functions are the fundamental
solution to the heat equation, i.e.,  $\Lom \phi_x = 0$, and (ii)
by definition any member of the null space $\mathcal{N}(S)$ annihilates
all Gaussian likelihood functions, i.e., $\int \phi_x d\nu = 0$ for all $x\in \mathbb R$.

\begin{prop}\label{prop-potential}
	Denote $\ker \;\Lom
	:=
	\{f\in C^\infty(\mathbb R^2):
	\Lom f=0\}.$
	Assume $\nu$ is a finite signed measure with compact support such that $\langle\nu,f\rangle=0,
	\;\forall\;f\in\ker\; \Lom,$ where $\langle\nu,f\rangle:=\int fd\nu$.
	Then, there is
	$$
	u\in L^2(\R^2)\cap\mathcal E'(\R^2)\quad\text{such that}\quad\check \Lom u=\nu,
	$$
    where $\mathcal E'(\R^2)$ is the space of compactly supported distributions.
\end{prop}

\begin{proof}
	By the solvability of heat equation in $\R^2$,  $\Lom$
	is surjective. Hence the induced quotient map
	$$
	\widetilde \Lom:
	C^\infty(\mathbb R^2)/\ker\; \Lom
	\longrightarrow
	C^\infty(\mathbb R^2),
	\qquad
	[f]\longmapsto \Lom f,
	$$
	is a continuous bijection between Fr\'echet spaces. By the open mapping
	theorem, $\widetilde \Lom^{-1}$ is continuous.
	
	Since $\langle\nu,f\rangle=0,
	\;\forall\;f\in\ker\; \Lom$,
    the continuous function $f\mapsto \langle \nu,f\rangle$ descends
    to a continuous linear functional
	$
	\widetilde\nu([f]):=\langle\nu,f\rangle
	$
    on
	$C^\infty(\mathbb R^2)/\ker \;\Lom$. Moreover, since $\nu$ is a finite measure with compact support,
	$
	u:=\widetilde\nu\circ\widetilde \Lom^{-1}
	$
    is a continuous linear functional on $C^\infty(\mathbb R^2)$.
    It follows from the standard characterization (see, e.g., \cite[Theorem 2.3.1]{Hormander1990}) $(C^\infty(\mathbb R^2))'=\mathcal E'(\mathbb R^2)$ that $u$ is a compactly supported distribution.
Moreover, by construction,
\begin{equation}
		\label{eq-transpose-equation}
	u(\Lom f)=\langle\nu,f\rangle,
	\quad \forall f\in C^\infty(\mathbb R^2),\qquad
	\text{or equivalently}\qquad
		\check \Lom u=\nu.
	\end{equation}
    
	Next, we prove that $u$ (regarding as a distribution) belongs to
	$L^2(\mathbb R^2)$.
	Let \((\xi,\tau)\) denote the Fourier variables dual to
	\((\mu,\sigma)\), respectively. Taking Fourier transforms in
	\eqref{eq-transpose-equation}
	gives
	\begin{equation}
		\label{eq-Fourier-equation}
		\Big(\frac12\xi^2-i\tau\Big)\widehat u(\xi,\tau)
		=
		\widehat\nu(\xi,\tau).
	\end{equation}
	Since \(\nu\) is a finite measure,
	$
	|\widehat\nu(\xi,\tau)|
	\le \|\nu\|_{\mathrm{TV}}
	=:M.
	$
Moreover,
		$
		\left|
		\frac12\xi^2-i\tau
		\right|^2
		=
		\frac14\xi^4+\tau^2.
		$
	Therefore,
	\begin{equation}
		\label{eq-Fourier-bound}
		|\widehat u(\xi,\tau)|^2
		\le
		\frac{M^2}{\tau^2+\xi^4/4}\; \text{ for } (\xi,\tau)\neq(0,0).
	\end{equation}
	Next, we claim that the right-hand side of
	\eqref{eq-Fourier-bound} is integrable away from the origin.
	First, for \(|\xi|\ge1\),
	$
	\int_{\mathbb R}
	\frac{d\tau}{\tau^2+\xi^4/4}
	=
	\frac{2\pi}{|\xi|^2},
	$
	and hence
	$	\int_{|\xi|\ge1}
		\int_{\mathbb R}
		\frac{d\tau\,d\xi}{\tau^2+\xi^4/4}
		=
		2\pi
		\int_{|\xi|\ge1}\frac{d\xi}{\xi^2}
		<\infty.$
	On the remaining region
	$
	|\xi|\le1,\; \xi^2+\tau^2\ge4,
	$
	we have \(|\tau|\ge\sqrt3\), and hence
	$
	\frac{1}{\tau^2+\xi^4/4}
	\le
	\frac{1}{\tau^2}.
	$
	It follows that
	\begin{equation}
		\label{eq-integrable}
		\frac{1}{\tau^2+\xi^4/4}
		\in
		L^1\bigl(\mathbb R^2\setminus B(0,2)\bigr).
	\end{equation}
	
	On the other hand, since $u$ is a compactly supported distribution, by Paley–Wiener-Schwartz theorem,
	$\widehat u$ is a smooth function on \(\mathbb R^2\)
	(in fact, the restriction of an entire function of exponential type)
    (e.g., \cite[Theorem 7.3.1]{Hormander1990}).
	In particular,
	$
	\widehat u
	$
    is 
    square-integrable in $B(0,2)$.
	This together with \eqref{eq-Fourier-bound} and \eqref{eq-integrable} leads to that
	\begin{equation}
		\label{eq-u-hat-L2}
		\widehat u\in L^2(\mathbb R^2).
	\end{equation}
	Therefore, Plancherel's theorem,
	together with the injectivity of the Fourier transform on tempered
	distributions, implies that $u$ is represented by an
	$L^2(\mathbb R^2)$-function. 
\end{proof}

We proceed with the following approximation result, which follows from
Malgrange's approximation theorem
\cite[Chapter I, Theorems 2 and $2'$]{Malgrange1956},
together with the elementary fact that every polynomial solution of
the heat equation is a finite linear combination of the heat polynomials
$H_n(\mu,\sigma)$, $n\geq 0$, defined as in \eqref{eq-heat-poly}.
\begin{lem}\label{lem:heat-poly-approx}
Finite linear combinations of the heat polynomials $H_n(\mu,\sigma)$,
$n\geq 0$, are dense in 
$\ker \;\Lom
	=
	\{f\in C^\infty(\mathbb R^2):
	\Lom f=0\}$
with respect to uniform convergence on compact subsets.
\end{lem}

We arrive at the fundamental representation for elements of the null space $\mathcal{N}(S)$.
\begin{prop}\label{thm-exist-pre-ani}
	Let $S\subset\mathbb R\times(0,\infty)$ be compact.
    For any $\nu\in\mathcal N(S)$, there is 
    \begin{equation}\label{eq-transpose-equation-2}
      u\in L^2(\R^2)\cap\mathcal E'(\R^2)\quad\text{such that}\quad \nu = \check \Lom u.
    \end{equation}
\end{prop}

\begin{proof}
	
	
Let $\nu\in\mathcal N(S)$ be arbitrary but fixed.
Because $\nu$ has compact support, the functional
$
	f\mapsto \langle\nu,f\rangle:=\int_S f\,d\nu
$
	is continuous on $C^\infty(\mathbb R^2)$. 
    Since $\langle \nu,\phi_x\rangle=0$ for all $x$, $\langle\nu,f\rangle=0$ for all $f\in \mathcal V_\phi$.
  It follows from Proposition \ref{prop:equi-Gauss-Hermite} that $\langle\nu,f\rangle=0$ for all $f\in \mathcal V_H$.
This combining with Lemma \ref{lem:heat-poly-approx} yields that
	\begin{equation*}
		\langle\nu,f\rangle=0,
		\qquad \forall f\in\ker\; \Lom.
	\end{equation*}
%
By Proposition \ref{prop-potential}, there is $u\in L^2(\R^2)\cap\mathcal E'(\R^2)$ such that
	$	\check \Lom u=\nu.$
The proof is complete.
\end{proof}

Next, we describe the support of the adjoint heat potential $u$ that arises in Proposition \ref{thm-exist-pre-ani}.
Intuitively, compact support of $u$ provides a region of vanishing at spatial infinity. Away from the source set, spatial analyticity makes this vanishing rigid: at each fixed scale level, along any horizontal interval containing no source, the potential cannot vanish on an open subinterval without vanishing throughout. The zero field therefore extends inward from spatial infinity until the source set is encountered. Therefore, at each fixed scale level, this eliminates the two exterior regions and leaves only the source points and the bounded gaps between them. These remaining regions form the $\mu$-directional hull $\mathcal K(S)$ of $S$, which is defined as follows.

\begin{deff}[$\mu$-directional hull]
	\label{def-mu-hull}
	Let $S\subset\mathbb R\times(0,\infty)$ be compact. 
	The $\mu$-directional hull of $S$ is defined by
	$$
	\mathcal K(S)
	:=
	S
	\cup
	\bigcup_{c>0}
	\left\{
	J\times\{c\}:
	J\text{ is a bounded connected component of }
	\mathbb R\setminus S_c
	\right\},
	$$
where for a set $K\subset\mathbb R\times (0,\infty)$, $K_c$ is the horizontal fiber at a scale level $c>0$ of $K$ defined by
$$
K_c:=\{\mu\in\mathbb R:(\mu,c)\in K\}.
$$
\end{deff}

The $\mu$-directional hull has a precise and useful description, as
given in the following lemma.
\begin{lem}\label{lem:K(S)-compact}
Let $S\subset\mathbb R\times(0,\infty)$ be compact. Then, 
$
(\mathcal K(S))_c=
[\min S_c,\max S_c]$ if $S_c\neq\emptyset,$
and
$(\mathcal K(S))_c=\emptyset$ if $S_c=\emptyset.$
As a result,
$\mathcal K(S)$ is compact.
\end{lem}

The proof of Lemma \ref{lem:K(S)-compact} is given in Appendix \ref{sec:proof34-Appendix}.
Next, we establish the following support result for $u$, which is similar to \cite[Theorems 5 and 16]{cias2026runge}.
However, we include a direct proof based on spatial analyticity of solutions of the heat equation.
\begin{prop}[Support localization]
	\label{lem-mu-support}
	Let $S\subset\mathbb R\times(0,\infty)$ be compact and let
$
	u\in\mathcal E'(\mathbb R^2)
$
	satisfy
$
	\operatorname{supp}(\check \Lom u)\subset S.
$
	Then
$
	\operatorname{supp}u\subset\mathcal K(S).
$
\end{prop}

\begin{proof}
	Let
$
(\mu_0,\sigma_0)\notin\mathcal K(S).
$
	Then $\mu_0$ belongs to an unbounded connected component $J'$ of
	$
	\mathbb R\setminus S_{\sigma_0}.
	$
	Assume, without loss of generality, that component $J'$ is unbounded to the
	left.
	
	Since $u$ is compactly supported, we may choose
	$
	\mu_1<\mu_0,
	\;
	\mu_1\in J',
	$
	so that $u$ vanishes in a neighborhood of
	$(\mu_1,\sigma_0)$. Since $\mu_1$ and $\mu_0$ lie in the same
	connected component $J'$, the horizontal segment
$
	K:=[\mu_1,\mu_0]\times\{\sigma_0\}
$
	is disjoint from $S$.
	Since $K$ is compact and $S$ is closed, there exists an open
	neighborhood $V$ of $K$ such that
$
	V\cap S=\emptyset.
$
	Therefore
$$
	\check\Lom u=0
	\qquad\text{in }V.
$$
Since
$
K
=
[\mu_1,\mu_0]\times\{\sigma_0\}
$
is a compact subset of the open set $V$, there exists small
$\eps>0$ such that the rectangle
$
(\mu_1-\eps,\mu_0+\eps)
\times
(\sigma_0-\eps,\sigma_0+\eps)
$ is contained in $V$. So,
$$
\check \Lom u=0
\qquad\text{in }(\mu_1-\eps,\mu_0+\eps)
\times
(\sigma_0-\eps,\sigma_0+\eps).
$$
Because of the choice of $\mu_1$, the distribution $u$ vanishes in a
neighborhood of
$
(\mu_1,\sigma_0).
$ We may therefore choose
$0<\eps'<\eps$ such that
$$
u=0
\quad\text{on}\quad
(\mu_1-\eps',\mu_1+\eps')
\times
(\sigma_0-\eps',\sigma_0+\eps')
\subset
(\mu_1-\eps,\mu_0+\eps)
\times
(\sigma_0-\eps,\sigma_0+\eps).
$$

	Next, by the hypoellipticity of the heat operator, $u\in C^\infty(V)$.
	The classical spatial-analyticity theorem for parabolic equations
	with analytic coefficients
	\cite{LandisOleinik1974} then implies that, for every fixed
	$\sigma$, the function $\mu\mapsto u(\mu,\sigma)$ is real
	analytic on each connected component of the corresponding spatial
	slice of $V$.

Fix
$
\sigma\in(\sigma_0-\eps',\sigma_0+\eps').
$
The real-analytic function
$
\mu\mapsto u(\mu,\sigma)
$
vanishes on the nonempty open interval
$
(\mu_1-\eps',\mu_1+\eps').
$
The identity theorem for real-analytic functions therefore gives
$$
u(\mu,\sigma)=0
\qquad
\text{for every }
\mu\in(\mu_1-\eps',\mu_0+\eps').
$$
Since $\sigma$ was arbitrary in $(\sigma_0-\eps',\sigma_0+\eps')$, it follows that
$$
u=0
\quad\text{on}\quad
(\mu_1-\eps',\mu_0+\eps')
\times
(\sigma_0-\eps',\sigma_0+\eps').
$$
In particular, $u$ vanishes in a neighborhood of
$
(\mu_0,\sigma_0).
$
Therefore $(\mu_0,\sigma_0)\notin\operatorname{supp}(u)$. The proof is complete.
\end{proof}

\begin{rem}
Due to the support localization result, from now on we will also refer to $\mathcal{K}(S)$ as the support hull for the adjoint heat potential $u$. 
\end{rem}
\subsubsection{Sufficient conditions impacting support of adjoint heat potentials}
\label{sec:sufficient-adjoint}

\begin{thm}
	\label{thm-sufficient-iden-K(S)}
	Let $S\subset\mathbb R\times(0,\infty)$ be compact such that its $\mu$-directional hull $\mathcal K(S)$ has two-dimensional Lebesgue measure zero.
	Then,
	$\mathcal P(S)$ is identifiable.
\end{thm}

\begin{proof}
Let $\nu\in \mathcal N(S)$.
By Proposition \ref{prop-potential}, there is $u\in L^2(\R^2)\cap\mathcal E'(\R^2)$ such that
	\begin{equation*}
		\check \Lom u=\nu.
	\end{equation*}
Therefore, because $u$ is compactly supported distributions and, by
	\eqref{eq-transpose-equation-2},
	$
	\operatorname{supp}\check \Lom u
	=
	\operatorname{supp}\nu
	\subset S,
	$
	the support localization, Proposition \ref{lem-mu-support} yields
	\begin{equation}
		\label{eq-support-u-1}
		\operatorname{supp}u\subset \mathcal K(S).
	\end{equation}
Thus, $\mathcal K(S)$ has two-dimensional Lebesgue measure zero. 
	Together with \eqref{eq-support-u-1}, this implies
	$
	u=0\text{ in }L^2(\mathbb R^2).
	$
	Equation
	\eqref{eq-transpose-equation-2} then gives
	$
	\nu=0.
	$
Therefore, the identifiability 
follows. 
\end{proof}

It is interesting to note that the theory of Runge approximation for linear partial differential operators is concerned precisely with conditions according to which the support of the adjoint potential $u$ is contained in $S$ itself.  See \cite{cias2026runge} and the references therein. In particular,
under the geometric condition \eqref{eq-sufficient-cond} below, provided that $\nu = \check \Lom u$, and $\nu$ is the supported in the compact $S$, then the adjoint heat operator $u$ is supported in $S$. The same vanishing argument therefore applies with the measure-zero assumption imposed directly on $S$.

\begin{thm}
	\label{thm-sufficient-iden}
	Let $S\subset\mathbb R\times(0,\infty)$ be compact and has two-dimensional Lebesgue measure zero.
	Suppose that there is a dense set $C$ of $\mathbb R$ such that
\begin{equation}\label{eq-sufficient-cond}
(\mathbb R^2\setminus S)\cap \R\times\{c\}\text{ has no bounded connected component for every }0<c\in C.
\end{equation}
	Then,
	$\mathcal P(S)$ is identifiable.
\end{thm}

\begin{proof}
Let $\nu\in\mathcal N(S)$, then
there is $u\in L^2(\R^2)\cap\mathcal E'(\R^2)$ such that
$		\check \Lom u=\nu$
as in Proposition \ref{prop-potential}. 
Since \eqref{eq-sufficient-cond} and \cite[Corollary 18]{cias2026runge},
	$(S,\mathbb R^2)
		\;\text{is an $\Lom$-Runge pair}.$
Therefore, because $u$ is a compactly supported distribution and, 
	$
	\operatorname{supp}\check \Lom u
	=
	\operatorname{supp}\nu
	\subset S,
	$
	the support characterization of Runge pairs
	\cite[Theorem 5]{cias2026runge} yields
		$\operatorname{supp}u\subset S.$
	As $S$ has Lebesgue measure zero, the same $L^2$
vanishing argument gives $u=0$ and hence $\nu=0$.
Therefore, the identifiability 
follows. 
\end{proof}

An immediate consequence of either Theorem \ref{thm-sufficient-iden-K(S)} or Theorem \ref{thm-sufficient-iden} is the identifiability of location-scale Gaussian mixtures when the support set $S$ is a graph over scale variables:
$S = \{(r(\sigma), \sigma): J\}$ where $J$ is a compact interval in $(0,\infty)$ and
$r$ any continuous function on $J$. See Corollary \ref{cor-continuous-scale-graph}. Theorem \ref{thm-sufficient-iden} is stronger than Theorem \ref{thm-sufficient-iden-K(S)}, as the former covers, for instance, the case $S$ is a graph over location variables satisfying a certain characteristic fiber condition, which is not applicable to the latter; see Corollary \ref{thm-fiber-condition}. Section \ref{sec:parameterized-curve} provides a catalog of identifiability results for set $S$ parameterized by curves. 

Although Theorem \ref{thm:thin-scale-proj}, Theorem \ref{thm-sufficient-iden-K(S)} and Theorem \ref{thm-sufficient-iden} combine to provide a considerably broad coverage of easy-to-verify topological conditions for identifiability, we will see that such sufficient conditions are not necessary. Moreover, the approach developed thus far does not tell us when a given set $S$ is not identifiable. 
(We already restricted $S$ to be topologically thin as if  
$\operatorname{Int}(S)\neq \emptyset$ then $\mathcal P(S)$ is non-identifiable;
see Proposition \ref{thm-ness-cond} in the Appendix.)
Obtaining necessary \emph{and} sufficient conditions of identifiability requires additional information about the support set $S$ and a probabilistic 
approach to be developed in the following section.

\section{Probabilistic characterization of the null space of Gaussian mixing operator} \label{sec:geometric}


In this section we will show that the Gaussian-null measures and their corresponding adjoint heat potentials admit a precise probabilistic representation based on the exit behavior of the space-time Brownian motion. This leads to a complete characterization of sufficient and necessary conditions for identifiability of location-scale Gaussian mixtures.

It is worth noting from the outset that the probabilistic approach to be described in the sequel is integral to the functional analytic framework developed in the prior section. Although this framework gives rise to the fundamental role of adjoint heat potentials via the representation $\nu = \check \Lom u$, what has been lacking so far is a suitable mechanism that allows for the possibility that the null space is non-empty. As discussed in Section \ref{sec:prelim}, the probabilistic viewpoint of Gaussian cancellation via convolution provides a natural mechanism based on which a stochastic process will be suitably constructed to trace the evolution of Gaussian sources and to quantify the cancellation phenomena associated with their contributions. 

To provide an intuition for our probabilistic construction, recall from Section \ref{sec:heat adjoint} that the adjoint-heat potential $u$ records the residual response generated by a Gaussian-null source $\nu\in\mathcal N(S)$. The representation $\nu=\check{\mathcal L}u$, together with the localization result of Proposition \ref{lem-mu-support}, namely,
$\operatorname{supp}u\subset\mathcal K(S)$, asserts that sources
carried by $S$ can leave a residual response inside the $\mu$-directional hull while their contributions must cancel outside it. This confinement
is a crucial property of the residual response, even as the individual contributions
continue to spread downward in scale levels. We now seek to understand how cancellation
at the allowed source locations can induce such a confined
response. A key question is whether and how $u$ can be nonzero while its finite signed source $\nu=\check \Lom u$ remains supported on $S$, because this corresponds precisely to the violation of identifiability.



For this purpose, it is natural to follow the spreading from
each source and identify the distribution of compensating
sources that would cancel its subsequent evolution.
Backward space--time Brownian motion provides this description.
Starting from $q=(\mu,\sigma)$, let
$
Y_t^q=(\mu+B_t,\sigma-t),
$
where $B_t$ is a standard Brownian motion.
For any $t\in (0,\sigma)$, the spatial coordinate of $Y_t^q$ is distributed by a normal distribution $N(\mu, t)$, its scale coordinate is fixed at $\sigma-t$. As a result, the Gaussian output produced by $Y_t^q$ is random, and in expectation at a spatial location $x$ it is precisely $\int \phi(x;y,\sigma-t) \phi(y; \mu, t) dy = \phi(x;\mu, \sigma)$. At the terminal time $t=\sigma$, the spatial position is $\mu+B_\sigma$, whose density is $\phi(x;\mu,\sigma)$. Thus, starting from any point $q=(\mu,\sigma)$, as one simulates the evolution of location-scale sources via the process $(Y_t^q)_{t\in [0,\sigma]}$, the Gaussian output remains preserved in expectation at any time. The identity $T\nu=0$ therefore expresses cancellation of the terminal distributions generated by the signed sources.

Now, let us simulate the process $Y_t^q$ with a starting source measure for the starting point $q$ supported in $S$. Since the support of the potential $u$ corresponding to a Gaussian-null measure is confined in the hull $\mathcal{K}(S)$, the simulated process can be stopped upon exiting $\mathcal{K}(S)$. Positive survival time identifies the points of $S$ that can generate occupation inside the hull. These are the sources that contribute to the residual response $u$. The exit law then supplies the boundary term in the source equation for $u$: the source of the stopped occupation equals the starting measure minus its exit distribution. A remarkable consequence of the Feynman-Kac formula is if a $\nu \in \mathcal{N}(S)$ is taken to be the starting source measure, then the occupation measure is precisely its adjoint potential $u$, which then implies that the exit law for the starting source $\nu$ must vanish. 

From here, one obtains a precise characterization for the null space: it is shown in Theorem \ref{thm-null-space-characterization} that every Gaussian-null measure arises as the difference of an admissible starting source measure and its exit distribution. Since the exit distribution is supported on the boundary of the hull, the requirement that $\nu = \check \Lom u$ remain supported on $S$ forces the exit distribution for the admissible starting source measure to vanish on the part of the boundary outside of $S$, where no compensating source is allowed. Taking into account the specific topological or geometric structure of the set $S$ and the support of the exit law introduced by the Brownian motion, now it is possible to determine whether the null space is trivial or not. Theorem \ref{thm:complete-iden-chacterization} formulates the necessary and sufficient conditions for identifiability in terms of the structure of set $S$ and the exit behavior of the Brownian motion.

\subsection{Backward space-time Brownian motion}
Let $B_t$ be a standard Brownian motion.
\begin{deff}[Space-time Brownian motion]
	For $q=(\mu,\sigma)\in\mathbb R\times (0,\infty)$, denote the space--time Brownian motion
	$$
	Y_t^q=Y_t^{(\mu,\sigma)}
	:=
	(\mu+B_t,\sigma-t),
	\qquad 0\le t\le \sigma.
	$$
\end{deff}

We start with the basic mechanism that is to replace an initial source by the
distribution of its locations at an intermediate stage, retaining
the remaining variance at each location. For a fixed time
$0<t<\sigma$, the Gaussian semigroup identity gives
$
\phi(x;\mu,\sigma)
=
\int_{\mathbb R}
\phi(x;y,\sigma-t)\,
\phi(y;\mu,t)\,dy.
$
Thus the original source and the redistributed sources produce
exactly the same Gaussian output. Lemma \ref{lem:Brownian-stop-rep} below allows the intermediate
time to depend on the path. 

\begin{lem}[Brownian stopping representation]
	\label{lem:Brownian-stop-rep}
	Let
	$
	\mu\in\mathbb R,
	\;
	\sigma>0.
	$
	Suppose that $\tau$ is a stopping time satisfying there is $\underline{\sigma}>0$ and
	\begin{equation}\label{eq-bounded-tau}
	0\le\tau\le \sigma-\underline{\sigma}
	\quad
	\text{a.s.}
	\end{equation}
	Define
	$
	Q_\tau
	:=
	\operatorname{Law}
	\left(
	\mu+B_\tau,\,
	\sigma-\tau
	\right).
	$
	Then
	\begin{equation}
		\label{eq:Brownian-mix-iden}
		p_{Q_\tau}(x)
		=
		\phi(x;\mu,\sigma)
		\qquad
		\text{for every }x\in\mathbb R.
	\end{equation}
	Equivalently,
	$
	p_{Q_\tau}(x)
	=
	p_{\delta_{(\mu,\sigma)}}(x),\;\forall x\in\mathbb R.
	$
\end{lem}

\begin{proof}
	Fix
	$
	x\in\mathbb R
	$
	and define
$X_t^x:=\phi(x;\mu+B_{t\wedge\tau},\sigma-t\wedge\tau).
	$
Combing It\^{o}'s formula, and the fact $\Lom\phi_x=0$, which
	the drift term vanishes, we obtain
	\begin{align*}
		X_t^x
		=\phi(x;\mu,\sigma)+
\int_0^{t\wedge\tau}\partial_\mu\phi
		\left(
		x;
		\mu+B_s,
		\sigma-s
		\right)dB_s.
	\end{align*}
    Hence
	$
	X_t^x
	$
	is a local martingale.
	By \eqref{eq-bounded-tau}, we have
	$
	0<X_t^x
	\le
	\frac{1}
	{\sqrt{2\pi\underline{\sigma}}}.
	$
	Thus
	$
	X_t^x
	$
	is a bounded martingale.
Therefore, 
	the optional stopping theorem yields
$$
		\phi(x;\mu,\sigma)
		=
		X_0^x
		=
		\mathbb E[X_{\tau}^x]
		=
		\mathbb E
		\left[
		\phi
		(
		x;
		\mu+B_{\tau},
		\sigma-\tau
		)
		\right]=p_{Q_\tau}(x).
$$
The proof is complete.    
\end{proof}
\begin{rem}
\label{rem:Fourier-stop}
	
Taking Fourier transforms in identity \eqref{eq:Brownian-mix-iden} gives
$
\int
e^{i\xi m-v\xi^2/2}\,Q_\tau(dm,dv)
=
e^{i\xi\mu-\sigma\xi^2/2},
$
for all
$
\xi\in\mathbb R.
$
Actually, this identity also follows by optional stopping of the complex
exponential martingale $\exp(i\xi B_t+\xi^2t/2)$.
\end{rem}

\begin{rem}\label{rem:non-identiable}
The above construction yields a recipe for non-identifiability: we may specify a starting distribution for the location-scale sources and allow the Brownian paths to determine where compensating sources are allocated. Suppose that the starting distribution and the stopped distribution are supported by
$A_0$ and $A_1$, respectively, both of which are two compact subsets of $\mathbb R\times(0,\infty)$. 
Start the process from $(M,\Sigma)\sim G_{A_0}\in\mathcal P(A_0)$,
independently of the Brownian motion $B_t$. If a stopping time $\tau\ge0$ satisfies that
$(M+B_\tau,\Sigma-\tau)\in A_1$ almost surely, then
optional stopping applied to the Gaussian kernel gives
$$
p_{G_\tau}=p_{G_{A_0}}, \qquad
\textrm{where }
G_\tau:=\operatorname{Law}(M+B_\tau,\Sigma-\tau)
\in\mathcal P(A_1).
$$
Thus stopping redistributes the sources from $A_0$ to $A_1$
while preserving their Gaussian output.

Whenever $\mathbb P(\tau>0)>0$, this redistribution is
nontrivial.
Hence $G_\tau\ne G_{A_0}$, so 
$\mathcal P(A_0\cup A_1)$
is not identifiable. Choosing the stopping time through
an exit rule from a suitable compact set connects this construction directly to the structure of the support. This is formalized by the exit operator that we will now describe.
\end{rem}

\subsubsection{Brownian interior, occupation measure, and exit operator}

The preceding construction produces equal Gaussian
mixture outputs by replacing an initial location-scale mixing distribution with its stopped law. We now ask whether arbitrary mixing distributions
with the same mixture density admit a common description of
this kind. More precisely, can the geometry of $S$ determine
a stopping rule under which equality of Gaussian mixtures is
equivalent to equality of the stopped distributions?

To address this question, we first define the following.
\begin{deff}
Let $K$ be a closed set of $\mathbb R\times (0,\infty)$.
Define the first exit time from $K$ by
$$
\tau_K^q=\inf\{t>0:Y_t^q\notin K\}.
$$
The backward Brownian interior of $K$ is
$$
O(K)
=
\{q\in K:\PP(\tau_K^q>0)>0\}.
$$
Define the exit law
$$
L_q^K=\operatorname{Law}(Y_{\tau_K^q}^q).
$$
For a finite signed Borel measure $\lambda\in\mathcal M(K)$, define the occupation measure $M_K\lambda$ by
$$
\langle M_K\lambda,f\rangle
=
\int_K
\mathbb E\left[\int_0^{\tau_K^q}f(Y_t^q)\,dt\right]\lambda(dq),
$$
for bounded Borel functions $f$.
Define the exit operator $E_K$ acting on $\lambda\in\mathcal M(K)$ by
\begin{equation}
    E_K\lambda(\cdot)= \int_KL_q^K(\cdot)\lambda(dq).
\end{equation}
\end{deff}

This stopped evolution produces two complementary objects.
The occupation measure $M_K\nu$ records the signed source
before exit, while the exit measure $E_K\nu$ records
the contribution at the exit locations.

Denote
$$
\underline{\sigma}_K=\min\{\sigma:(\mu,\sigma)\in K\},\qquad
\overline{\sigma}_K=\max\{\sigma:(\mu,\sigma)\in K\}.
$$
We first obtain some basic properties.
\begin{lem}\label{lem:exit-law}
    For every nonempty and compact set  $K\subset\R\times(0,\infty)$, we have the followings:
\begin{enumerate}
\item[(i)]
The set $O(K)$ is Borel, and $q\mapsto L_q^K$ is a
Borel probability kernel. Moreover, $\operatorname{Int} K \subset O_K$. 

\item[(ii)]
For every $q=(\mu,\sigma)\in K$,
\begin{equation}\label{eq-bounded-stopping}
0\le\tau_K^q\le \sigma-\underline{\sigma}_K\le \overline{\sigma}_K-\underline{\sigma}_K.
\end{equation}
Moreover,
$$
q\in O(K)
\quad\Longleftrightarrow\quad
\PP(\tau_K^q>0)=1.
$$

\item[(iii)]
For every $q \in K$,
$
L_q^K\bigl(\partial K\setminus O(K)\bigr)=1
$. If $q\notin O(K)$, which entails $q\in \partial K \setminus O(K)$, then $L_q^K=\delta_q$.
\end{enumerate}
\end{lem}
\begin{proof}
Let $(\mathcal F_t)_{t\ge0}$ be the
	completed natural filtration of $B_t$.
    On continuous Brownian path space and by closeness of $K$, for $r>0$,
$
	\{\tau_K^q\ge r\}
	=
	\bigcap_{t\in\mathbb Q\cap(0,r)}
	\{Y_t^q\in K\}.
	$
Evaluation of the path at $\tau_K^q$ is measurable, so the exit
laws form a Borel kernel. The last claim regarding the inclusion 
property of $O(K)$ is proved below.

Now, we move to part (ii).
It is clear that the time coordinate implies
$\tau_K^q\le \sigma-\underline{\sigma}_K$.
Next, for any $\varepsilon>0$, choose an integer $m_0$ such that
	$1/m_0<\varepsilon$. Then
	$$
	\{\tau_K^q>0\}
	=
	\bigcup_{m\ge m_0}
	\{\tau_K^q\ge 1/m\}
	\in\mathcal F_\varepsilon.
	$$
	Consequently,
	$$
	\{\tau_K^q>0\}
	\in
	\mathcal F_{0+}
	:=
	\bigcap_{\varepsilon>0}\mathcal F_\varepsilon.
	$$
	By Blumenthal's zero--one law,
	\begin{equation}\label{eq-0-1-tau}
\mathbb P(\tau_K^q>0)\in\{0,1\}.	    
	\end{equation}
Therefore, $q\in O(K)$ is equivalent to $\mathbb P(\tau_K^q>0)=1$.
This also entails that $O(K) \supset \operatorname{Int} K$,
since a continuous path
starting from any $q\in \operatorname{Int} K$ remains in $K$ for a positive time.

Next, we prove part (iii). 
Note that if $\tau_K^q>0$, continuity gives
$Y_{\tau_K^q}^q\in\partial K$. If $\tau_K^q=0$, then as argued
above, $q \in \partial K$. In either cases, $Y_{\tau_K^q}^q\in\partial K$ almost surely.

By the definition of $\tau_K^q$, after time $\tau_K^q$ the path leaves
$K$ arbitrarily soon; that is, for every $n\ge1$,
$$
\mathbb P_q\!\left(
Y_{\tau_K^q+s}^q\notin K \text{ for some } s\in(0,1/n)
\,\middle|\,\mathcal F_{\tau_K^q}
\right)=1
\qquad\text{a.s.}
$$
Here $\PP_q$ indicates that the starting point is $q$.
By the strong Markov property at $\tau_K^q$, 
$$
\mathbb P_{Y_{\tau_K^q}^q}(\tau_K<1/n)=1
\quad\text{for every }n,
$$
and hence
$
\mathbb P_{Y_{\tau_K^q}^q}(\tau_K=0)=1
\;\text{a.s.}
$
That means $Y_{\tau_K^q}^q\notin O(K)$ a.s.
Therefore the exit law $L_q^K$ is concentrated on $\partial K\setminus O(K)$.
Moreover, at an starting point $q\notin O(K)$,
$\tau_K^q=0$ almost surely and $L_q^K=\delta_q$.
\end{proof}

\begin{lem}
\label{lem:exit-preserve}
The operator $E_K$ preserves total mass and satisfies
$$
\|E_K\lambda\|_{\TV}\le\|\lambda\|_{\TV},
\qquad
E_K^2=E_K,
$$
and
\begin{equation}\label{eq-preserve-mixture-O}
T(E_K\lambda)=T\lambda.
\end{equation}
\end{lem}
\begin{proof}
Mass preservation and the total variation contraction follow
from the probability kernel representation of $E_K$.
By Lemma \ref{lem:exit-law}(iii), $L_q^K$ is concentrated on 
$\partial K \setminus O(K)$. Moreover, if $q\notin O(K)$, then $L_q^K=\delta_q$.
Therefore, $E_K$ acts
as the identity on measures concentrated on $K\setminus O(K)$, so
$E_K^2=E_K$.

	Fix
	$
	x\in\mathbb R,
	$
by Lemma \ref{lem:Brownian-stop-rep}, we have
	$
		\phi(x;\mu,\sigma)
		=
		\mathbb E
		[
		\phi
		(
		x;
		\mu+B_{\tau_K^q},
		\sigma-\tau_K^q
		)
		].
$
Integrating against $\lambda$ proves
	\eqref{eq-preserve-mixture-O}.
\end{proof}

\begin{lem}\label{prop-W_K}
 Let $K\subset\mathbb R\times(0,\infty)$ be nonempty and compact.
    For $\lambda\in\mathcal M(K)$, we have the following 
    \begin{itemize}
        \item[(i)] $\operatorname{supp}(M_K\lambda)\subset K$ and $\|M_K\lambda\|_{\TV}
\le(\overline{\sigma}_K-\underline{\sigma}_K)\|\lambda\|_{\TV}$.
        \item[(ii)] $\check {\mathcal L} (M_k\lambda)=\lambda-E_K\lambda$ and $M_K\lambda=M_K(\lambda|_{O(K)})$.
    \end{itemize}
\end{lem}

\begin{proof}
The support and norm assertions follow directly from the
definition and the bound on $\tau_K^q$.
We next prove part (ii). Since the generator of the backward Brownian motion is $-\Lom$, Dynkin's formula gives
for any $\varphi\in C_c^\infty(\R^2)$
$$
\mathbb E\varphi(Y_{\tau_K^q}^q)-\varphi(q)
=
-\mathbb E\int_0^{\tau_K^q}\mathcal L\varphi(Y_t^q)\,dt.
$$
Integration against $\lambda$ yields
$$
\langle M_K\lambda,\mathcal L\varphi\rangle
=
\langle\lambda-E_K\lambda,\varphi\rangle.
$$
It shows that $\check {\mathcal L} (M_k\lambda)=\lambda-E_K\lambda$.
Finally, if $q\notin O(K)$ then $\tau_K^q=0$ a.s., so
the occupation time is zero almost
surely. This proves the final assertion.
\end{proof}

The next proposition provides the fundamental linkage between occupation measure and adjoint heat potential associated with a source measure.

\begin{prop}\label{prop:pre-anni-rep-exit}
    Let $K\subset\mathbb R\times(0,\infty)$ be nonempty and compact.
    	Suppose that $u\in L^1(\mathbb R^2)$ satisfies
	$$
	\operatorname{supp}u\subset K,
	\qquad
	\check{\mathcal L}u=\nu,
	$$
    where $\nu$ is a finite signed Borel measure on $\mathbb{R}^2$.
    Then, the continuous linear functional associated with the distribution
    $u$ is given by,
    for every $f\in C_c^\infty(\mathbb R^2)$,
	\begin{equation}\label{eq-F-K-rep}
		\langle u,f\rangle
		=
		\int_K
		\mathbb E\left[
		\int_0^{\tau_K^q}f(Y_t^q)\,dt
		\right]\nu(dq).
	\end{equation}
Equivalently,
\begin{equation}
u=M_K\nu.
\end{equation}
\end{prop}
\begin{proof}
Since differential operators do not enlarge distributional supports,
$
\operatorname{supp}\nu
\subseteq \operatorname{supp}u \subseteq K.
$
Choose a nonnegative function $w\in C_b^\infty(\mathbb R^2)$
whose zero set is exactly $K$. Then $wu=0$ almost everywhere, so
$$
(\check{\mathcal L}+nw)u=\nu,
\qquad n\ge1.
$$

Take a $\sigma_0<\underline{\sigma}_K$, and consider an arbitrary function 
$f\in C_c^\infty(\mathbb R^2)$ such that $f$ vanishes
in a neighborhood of $\{\sigma\le\sigma_0\}$.
For $\sigma>\sigma_0$, define
$$
\psi_n(\mu,\sigma)
=
\mathbb E\!\left[
\int_0^{\sigma-\sigma_0}
\exp\!\left(-n\int_0^t w(Y_r^{(\mu,\sigma)})\,dr\right)
f(Y_t^{(\mu,\sigma)})\,dt
\right],
$$
and extend $\psi_n$ by zero to $\{\sigma\le\sigma_0\}$.
In addition, on each half-space $\{(\mu,\sigma)| \sigma\le\Sigma\}$, with $\Sigma>\sigma_0$,
the upper integration limit in the above definition may be replaced by $\Sigma-\sigma_0$, since 
$f(Y_t^{(\mu,\sigma)}) = f(\mu+B_t, \sigma-t) = 0$ when $t \geq \sigma-\sigma_0$.
Differentiation under
the integral and expectation is justified by the boundedness of
all derivatives of $w$ and $f$. Thus $\psi_n\in C^\infty(\mathbb R^2)$,
including across $\sigma=\sigma_0$.

The time-dependent, inhomogeneous Feynman--Kac formula
(see \cite[Section~5.7]{shreve1991brownian}) gives
$$
(\mathcal L+nw)\psi_n=f.
$$
Moreover, since $w\ge0$,
$$
|\psi_n(q)|
\le (\overline{\sigma}_K-\sigma_0)\|f\|_\infty,
\qquad \forall q=(\mu,\sigma)\in K,\quad n\ge1.
$$

Choose $\chi\in C_c^\infty(\mathbb R\times(0,\infty))$
equal to one on a neighborhood of $K$.
We obtain for every $n\geq 1$
$$
\begin{aligned}
\langle u,f\rangle
=\langle u,(\mathcal L+nw)\psi_n\rangle 
=\langle u, \mathcal L \psi_n \rangle 
= \langle u,\mathcal L (\chi\psi_n)\rangle
=\langle\check{\mathcal L} u,\chi\psi_n\rangle
 =\int_K\psi_n(q)\,\nu(dq).
\end{aligned}
$$
In the above display, the second equality is due to $wu = 0$ almost everywhere, the third equality holds because all derivatives of $\chi$ vanish on 
a neighborhood of the support of $u$, and the fourth equality is due to distributional 
duality and $\chi \psi_n \in C_{c}^{\infty}(\mathbb R^2)$. 

To arrive at the desired claim requires taking the obtained expression in the limit of $n\rightarrow \infty$. For each fixed $q=(\mu,\sigma)\in K$, put
$$
A_t :=\int_0^t w(Y_r^q)\,dr,
\qquad \tau:=\tau_K^q.
$$
On every continuous sample path, $Y_r^q \in K$ as long as $r\leq \tau$. It follows that $w(Y_r^q)= 0$ and so $A_t=0$ for $t\leq \tau$. On the other hand, for $t>\tau$, the sample path visits
$K^c$ before $t$, and continuity gives a non-empty interval on which
$w(Y_r^q)>0$, resulting in $A_t>0$. It follows that
$$
e^{-nA_t}\longrightarrow \mathbf 1_{\{t\leq\tau\}}
\quad\text{for Lebesgue-almost every }t.
$$
Since $\tau_K^q\le\sigma-\underline{\sigma}_K<\sigma-\sigma_0$,
dominated convergence in time and probability yields
$$\psi_n(q)\longrightarrow
\mathbb E\!\left[\int_0^{\tau_K^q}f(Y_t^q)\,dt\right].$$
The uniform bound above and $|\nu|(K)<\infty$ now allow
dominated convergence against $\nu$, giving
$$\langle u,f\rangle
=
\int_K
\mathbb E\!\left[\int_0^{\tau_K^q}f(Y_t^q)\,dt\right]
\,\nu(dq).$$
This proves $u=M_K\nu$.
\end{proof}

An immediate consequence of Proposition \ref{prop:pre-anni-rep-exit} is the first main theorem of this section. It says that every Gaussian cancellation supported on $S$ is detected as a cancellation by the exit operator.
This is the first characterization of the null space of the Gaussian mixing operator.

\begin{thm}[Null space characterization]
\label{prop:pre-annihator-rep-via-exit}
 Let $S\subset\mathbb R\times (0,\infty)$ be 
 nonempty and 
 compact.
    Then,
    \begin{equation}
     \mathcal N(S)=\{\nu\in\mathcal M(S): E_{\mathcal K(S)}\nu=0\}.    \end{equation}
\end{thm}
\begin{proof}
Let $\nu\in\mathcal N(S)$.
Applying Proposition \ref{thm-exist-pre-ani}, there is 
    $
      u\in L^2(\R^2)\cap\mathcal E'(\R^2)\text{ such that }\check \Lom u=\nu,
  $
  and by Proposition \ref{lem-mu-support}, $\operatorname{supp} u\subset \mathcal K(S)$.
  By Proposition \ref{prop:pre-anni-rep-exit},
  $u=M_{\mathcal K(S)}\nu$.
Applying Lemma \ref{prop-W_K}(ii)
now gives
$$
\nu=\check{\mathcal L}u=\check{\mathcal L}(M_{\mathcal K(S)}\nu)=\nu-E_{\mathcal K(S)}\nu,
$$
hence $E_{\mathcal K(S)}\nu=0$.

Conversely, assume $\nu\in \mathcal M(S)$ satisfying $E_{\mathcal K(S)}\nu=0$.
By \eqref{eq-preserve-mixture-O} in Lemma \ref{lem:exit-preserve}, we have $T\nu=T(E_{\mathcal K(S)}\nu)=0.$
Thus $\nu\in\mathcal N(S)$.
\end{proof}

\subsection{Sufficient and necessary conditions for identifiability}

We are now in position to formulate sufficient and necessary conditions for Gaussian mixture identifiability in terms of the topological and geometric structure of the source measure $S$. To this end, we present a representation theorem for the members of the null space, which follows from Theorem \ref{prop:pre-annihator-rep-via-exit}.

\begin{thm}[Null-space representation]\label{thm-null-space-characterization}
    Let $S\subset\mathbb R\times (0,\infty)$ be nonempty and compact.
    Denote
    \begin{equation}
       A=S\cap O(\mathcal K(S)),\qquad F=\partial (\mathcal K(S))\setminus S.
    \end{equation}
    Then,
    \begin{equation}
    \label{eq:null-measure-rep}
        \mathcal N(S)=\{\lambda-E_{\mathcal K(S)}\lambda: \lambda\in\mathcal M(A), \;(E_{\mathcal K(S)}\lambda)\vert_{F}=0\}
    \end{equation}
\end{thm}
\begin{proof}
    For this proof, let $K$ stand for $\mathcal K(S)$.
    Let $\nu\in \mathcal N(S)$ and be composed as
    \begin{equation}
        \nu=\lambda+\eta,\text{ where }\lambda=\nu\vert_A,\;\eta=\nu\vert_{S\setminus A}.
    \end{equation}
  Since $E_{K}$ is the exit operator with respect to $K$, it acts as identity operator on measure concentrated on $K\setminus O(K)\supset S\setminus A$ (cf. Lemma \ref{lem:exit-preserve}). Thus, by Theorem \ref{prop:pre-annihator-rep-via-exit}
  \begin{equation}
      0=E_{K}\nu=E_{K}\lambda +\eta.
  \end{equation}
  As a result, 
  \begin{equation}
  \label{eq:null-measure-rep-2}
\nu=\lambda-E_{K}\lambda.
  \end{equation}
Moreover, $E_{K}\lambda=-\eta$, which is concentrated on $S\setminus O(K)$. We also know that $E_{K} \lambda$ is concentrated on $\partial K \setminus O(K)$, so it must concentrate all its mass on $(\partial K \cap S) \setminus O(K)$. Hence,
$(E_{K}\lambda)\vert_F=0$.

Conversely, suppose that a measure $\nu$ admits the representation \eqref{eq:null-measure-rep-2}, where $\lambda\in\M(A)$ satisfies $(E_{K}\lambda)|_F=0$.
By Lemma \ref{lem:exit-law}(iii), $E_{K}\lambda$ is concentrated on $\partial K\setminus O(K) \subset \partial K$. It follows that $E_{K}\lambda$ concentrates all its mass on $S$. Since $S$ is closed, this entails
$\operatorname{supp}(E_{K}\lambda)\subset S$.
Clearly, $\operatorname{supp}(\lambda) \subset S$.
Hence
$$
\nu=\lambda-E_{K}\lambda\in\M(S).
$$
By \eqref{eq-preserve-mixture-O}, $T(E_{K}\lambda)=T\lambda$. Which shows that $\nu\in\mathcal N(S)$.
\end{proof}

\begin{rem}
    The representation \eqref{eq:null-measure-rep} establishes that every Gaussian mixture cancellation arises as a redistributing scheme of location-scale sources guided by the backward Brownian diffusion process. This process starts with an admissible source measure $\lambda$, which concentrates all its mass on the set $A=S \cap O(\mathcal K(S))$, and after a positive time it ends with the compensating source measure $E_{\mathcal K(S)}(\lambda)$, which is concentrated on a disjoint destination, the set $(S \cap \partial (\mathcal K(S)))\setminus O(\mathcal K(S))$. Not all measures in $\mathcal M(A)$ are allowed to start the scheme: the admissible ones are those and only those $\lambda$ for which their exit law is concentrated on the intended destination. Or equivalently, the restriction of the exit law to $F=\partial (\mathcal K(S)) \setminus S$ is zero.
\end{rem}

\begin{thm}\label{thm:complete-iden-chacterization}
 Let $S\subset\mathbb R\times (0,\infty)$ be nonempty and compact.
 \begin{itemize}
 \item[(i)] If 
$S\cap O(\mathcal K(S))=\emptyset$, then $\mathcal P(S)$ is identifiable.
\item[(ii)]  If $\partial (\mathcal K(S))\subset S$, then
$$
\mathcal P(S)\text{ is identifiable}
\quad\Longleftrightarrow\quad
S\cap O(\mathcal K(S))=\emptyset.
$$
\end{itemize}
\end{thm}
\begin{proof}
  If $S\cap O(\mathcal K(S))=\emptyset$, by Theorem \ref{thm-null-space-characterization} leads to that $\mathcal N(S)=\{0\}$, which proves identifiability.
  To prove part (ii), it is noted that in this case, $\partial (\mathcal K(S))\setminus S=\emptyset$.
  Therefore, if $S\cap O(\mathcal K(S))\neq\emptyset$, then $\mathcal N(S)=\{\lambda-E_{\mathcal K(S)}\lambda: \lambda\in\mathcal M(A)\}$, which always contains non-zero elements. 
\end{proof}

\begin{rem}
(i) Notice that the necessary and sufficient conditions are characterized directly in terms of the topological and geometric structure of $S$, via its $\mu$-direction hull $\mathcal K(S)$, the boundary $\partial \mathcal K(S)$ and the backward Brownian interior $O(\mathcal K(S))$.
(ii) It is remarkable that Theorem \ref{thm:thin-scale-proj} follows immediately from Theorem \ref{thm:complete-iden-chacterization}. Indeed, if $\operatorname{Int} \pi_\sigma(S) = \emptyset$, the presence of scale gaps means that every point in $S$ exits immediately from $K(S)$, so $O(\mathcal K(S)) = \emptyset$. Thus $\mathcal P(S)$ is identifiable. It is striking how the conclusion is achieved swiftly here in contrast to the original proof. (iii) Thanks to Proposition \ref{prop:emptyinterior} below, in a similar fashion Theorem \ref{thm-sufficient-iden-K(S)} also follows quickly by applying part (i) of Theorem \ref{thm:complete-iden-chacterization}. 
(iv) Part (ii) of Theorem \ref{thm:complete-iden-chacterization} provides an useful tool for us to provide a complete characterization of identifiability for a variety of cases where the support $S$ is a parameterized curve. 

\end{rem}
The following proposition provides an interesting fact about these geometric objects involved.

\begin{prop}
\label{prop:emptyinterior}
Let $K\subset\mathbb R\times(0,\infty)$ be compact. 
	Suppose that every nonempty horizontal fiber
	$$
	K_\sigma
	:=
	\{\mu\in\mathbb R:(\mu,\sigma)\in K\}
	$$
	is an interval.
Then,
    \begin{equation}
    \operatorname{Int}K=\emptyset
		\quad\Longleftrightarrow\quad
		O(K)=\emptyset.
    \end{equation}
\end{prop}
\begin{proof}
    If $O(K) = \emptyset$ then the inclusion property $O(K) \supset \operatorname{Int} K$ from Lemma \ref{lem:exit-law}(i) entails that $\operatorname{Int} K = \emptyset$.
%
    For the converse, assume that $O(K) \neq \emptyset$. Then we will show that $\operatorname{int } K \neq \emptyset$.
    Take a $q=(\mu_0,\sigma_0)\in O(K)$.
    By Lemma \ref{lem:exit-law}(ii), 
    $\PP(\tau_K^q>0)=1.$
    This means
    $$
	\lim_{h\downarrow0}\mathbb P(\tau_K^q>h)=1.
	$$
	But
	$
	\{\tau_K^q>h\}
	\subset
	\{Y_h^q\in K\}
	=
	\{\mu_0+B_h\in K_{\sigma_0-h}\},
	$
	so	$
	\mathbb P(\mu_0+B_h\in K_{\sigma_0-h})
	\longrightarrow1$ as $h\downarrow 0$.

	Let $\Phi$ denote the standard normal distribution
	function. There exists $\delta\in(0,\sigma_0)$ such that
	\begin{equation}\label{eq-Phi-1}
	\mathbb P(\mu_0+B_h\in K_{\sigma_0-h})>\Phi(1)
	\qquad\text{for every }0<h<\delta.
	\end{equation}
	For each such $h$, the fiber $K_{\sigma_0-h}$ is nonempty.
	By compactness and the interval-fiber assumption,
	we may write
	$
	K_{\sigma_0-h}=[\ell_h,r_h]$. Thus, the probability bound \eqref{eq-Phi-1} may be written as
    \begin{equation}
        \mathbb P(B_h\in [l_h-\mu_0,r_h-\mu_0]) > \Phi(1).
    \end{equation}
    Since $(1/\sqrt{h}) B_h \sim N(0,1)$, $\mathbb P(B_h \leq \sqrt{h})= \mathbb P(B_h\geq -\sqrt{h}) = \Phi(1)$.
    Combining this identity with the previous bound yields $r_h -\mu_0 > \sqrt{h}$ and $l_h-\mu_0 < -\sqrt{h}$, for all $h\in (0,\delta)$.


Since $[\mu_0-\sqrt{h},\mu_{0}+\sqrt{h}] \subset (l_h,r_h) = \operatorname{Int} K_{\sigma_0-h}$ for all $h\in (0,\delta)$,
it follows that 
the nonempty open rectangle
	$$
	R
	=
	\big(
	\mu_0-\frac{\sqrt\delta}{2},
	\mu_0+\frac{\sqrt\delta}{2}
	\big)
	\times
	\big(
	\sigma_0-\frac{3\delta}{4},
	\sigma_0-\frac{\delta}{4}
	\big).
	$$
    is contained in $K$. Thus, $\operatorname{Int} K \neq \emptyset$. 
    This concludes the proof.
\end{proof}

\section{Identifiability and canonical mixing representation on curved supports}\label{sec:5}

The purpose of this section is multifold. First, Section \ref{sec:parameterized-curve} serves as an illustration and application of the general identifiability theorems obtained in the previous sections. We give a catalog of concrete examples of identifiability, or the violation thereof, for a variety of support set $S$ parameterized by a compact curve in $\mathbb R\times (0,\infty)$. In many examples the application boils down to analyzing the behavior of the backward space-time Brownian motion departing from $S$. The exit operator also suggests a special role for such classes of mixing measures, as it will be shown that every location-scale Gaussian mixture with compact mixing support can be represented uniquely by a mixing measure supported on a fixed boundary curve of a rectangular or circular enclosure. This representation result may be potentially useful for modeling applications, and is presented in Section \ref{sec:canonical-rep}. Finally, in Section \ref{sec:transformation} we explore the role of parameterization, which determines the shape of the curved support for the mixing measure, on the Gaussian mixture identifiability.

\subsection{Examples for curved supports}\label{sec:parameterized-curve}
Although Theorems \ref{thm-sufficient-iden-K(S)} and \ref{thm-sufficient-iden} provide only sufficient conditions for identifiability, they are almost straightforward to verify if applicable.
Theorems \ref{thm-null-space-characterization} and \ref{thm:complete-iden-chacterization} provides general and sharper characterization. Verification of sufficient or necessary conditions from such theorems requires an analysis of the exit behavior of the space-time backward Brownian motion. In the examples below such verification is mostly quite simple.

\subsubsection{Single graphs}
\label{sec:single-graph}


First, we consider $S$ being the single graph of some function of $\sigma$ or some function of $\mu$.

\begin{cor}[Graph over the scale variable]
	\label{cor-continuous-scale-graph}
	Let $J\subset(0,\infty)$ be a compact interval and let
	$
	r
	$
	be any continuous function on $J$.
	Set
$$
	S=\{(r(\sigma),\sigma):\sigma\in J\}.
$$
Then $\mathcal P(S)$ is
	identifiable.
\end{cor}

\begin{cor}[Characteristic-fiber condition]
	\label{thm-fiber-condition}
	Let $I\subset \mathbb R$ be compact and let
	$
	s
    \in C(I;(0,\infty)).
	$
	Set
	$$
	S:=\{(\mu,s(\mu)):\mu\in I\}\subset \mathbb R\times(0,\infty).
$$
	Assume that
	\begin{equation}
		\label{eq-bad-levels}
		\operatorname{int}
		\left\{
		c\in\mathbb R:
		s^{-1}(c)
		\text{ is disconnected as a subset of }\mathbb R
		\right\}
		=\emptyset.
	\end{equation}
	Then
	$\mathcal P(S)$ is identifiable.
\end{cor}

\begin{figure}[H]
	\centering
	\captionsetup[subfigure]{justification=centering}
	
	\begin{subfigure}[t]{0.32\textwidth}
		\centering
		
		\begin{tikzpicture}[scale=0.7]
			\draw[->, thick] (0,0) -- (5.3,0)
			node[right] {$\mu$};
			\draw[->, thick] (0,0) -- (0,4.6)
			node[above] {$\sigma$};
			
			\draw[
			blue,
			very thick,
			smooth
			]
			plot coordinates {
				(1.0,0.9)
				(1.55,1.15)
				(1.95,1.45)
				(1.75,1.75)
				(1.25,2.05)
				(1.45,2.35)
				(2.05,2.65)
				(2.45,2.95)
				(2.10,3.25)
				(2.35,3.55)
				(2.95,3.85)
				(3.35,4.15)
			};
			
			\fill[blue] (1.0,0.9) circle (2pt);
			\fill[blue] (3.35,4.15) circle (2pt);
		\end{tikzpicture}
		\caption{Identifiable, by Corollary \ref{cor-continuous-scale-graph}:
			$S=\{(r(\sigma),\sigma):\sigma\in J\},
			$
			where $J$ is compact and $r\in C(J)$.}
		
			
	\end{subfigure}
	\hfill
	\begin{subfigure}[t]{0.32\textwidth}
		\centering
		
		\begin{tikzpicture}[scale=0.7]
			\draw[->, thick] (0,0) -- (5.3,0)
			node[right] {$\mu$};
			\draw[->, thick] (0,0) -- (0,4.6)
			node[above] {$\sigma$};
			
			\draw[
			blue,
			very thick,
			smooth,
			samples=120,
			domain=1.75:3.55
			]
			plot (
			{\x},
			{1.0 + 1.8*(\x-2.4)^2}
			);
			
			\fill[blue]
			(1.75,{1.0 + 1.8*(1.75-2.4)^2})
			circle (2pt);
			
			\fill[blue]
			(3.55,{1.0 + 1.8*(3.55-2.4)^2})
			circle (2pt);
		\end{tikzpicture}
		\caption{
        Identifiable, by Corollary \ref{prop:iden-graph-over-location}:
         $S=\{(\mu,s(\mu)):\mu\in I\},
			$
			where $I$ is an interval and $s$ satisfies conditions in Corollary \ref{prop:iden-graph-over-location} (i).}
		
			
	\end{subfigure}
	\hfill
\hfill
\begin{subfigure}[t]{0.32\textwidth}
    \centering
    \begin{tikzpicture}[scale=0.7]
        \draw[->, thick] (0,0) -- (5.3,0)
            node[right] {$\mu$};
        \draw[->, thick] (0,0) -- (0,4.6)
            node[above] {$\sigma$};

        \draw[dashed, gray] (2.6,0) -- (2.6,3.7);
        \node[below] at (2.6,0) {$\mu^*$};

        \draw[
            blue,
            very thick,
            smooth,
            samples=80,
            domain=1.0:2.6
        ]
            plot (\x,{1.0 + 2.7*exp(\x-2.6)});

        \draw[
            blue,
            very thick,
            smooth,
            samples=80,
            domain=2.6:4.2
        ]
            plot (\x,{1.0 + 2.7*exp(2.6-\x)});

        \fill[blue]
            (1.0,{1.0 + 2.7*exp(-1.6)})
            circle (2pt);
        \fill[blue]
            (4.2,{1.0 + 2.7*exp(-1.6)})
            circle (2pt);
    \end{tikzpicture}

    \caption{
        Identifiable by Corollary~
        \ref{prop:iden-graph-over-location}\textup{(ii)}:
        $s$ has two $C^1$ branches, with
        $s'>0$ to the left of $\mu^*$ and
        $s'<0$ to its right.
    }
\end{subfigure}
		
			
			
			
			
			
			
	\caption{Examples of support classes as single graphs.}
	\label{fig:single graph}
\end{figure}
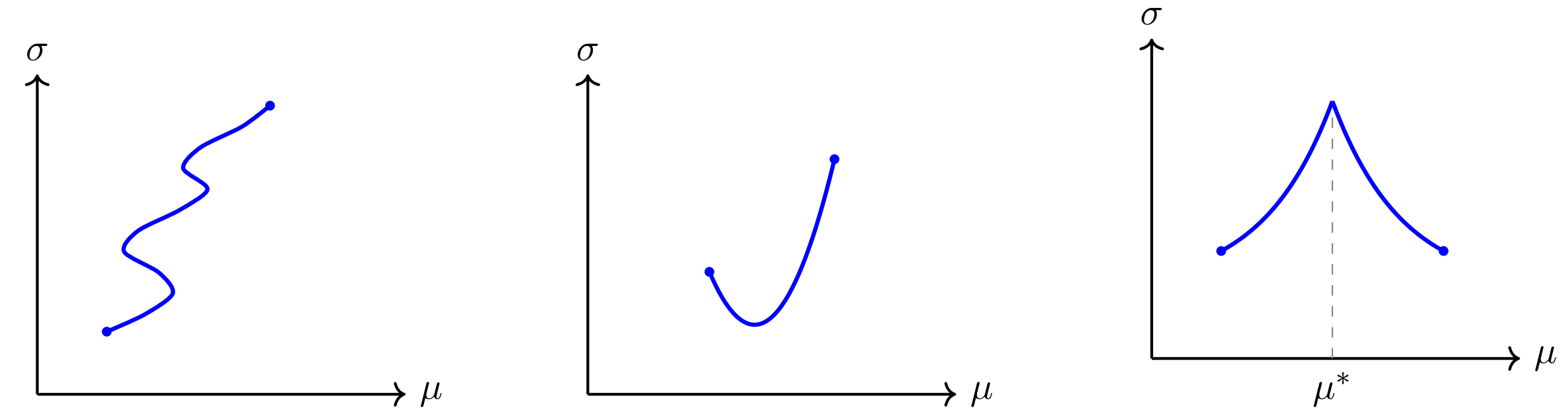

\begin{cor}[Graph over the location variable]
	\label{prop:iden-graph-over-location}
	Let $I=[a,b]$, with $a<b$, and let
	$s:I\to(0,\infty)$ be continuous. Set
	$$
	S=\{(\mu,s(\mu)):\mu\in I\}.
	$$
	Then 
	$\mathcal P(S)$ is identifiable under any of the following conditions:
	\begin{enumerate}
		
		\item[(i)]
		$s$ is quasiconvex, i.e., the level sets $\{\mu: s(\mu) \leq c\}$ are connected intervals for all $c > 0$.
		In particular, this includes the case $s$ is convex.
		
		
		\item[(ii)]
		There exists $\mu^*\in(a,b)$ such that
		$$
		s\in C^1((a,\mu^*))\cap C^1((\mu^*,b)),
		\quad
		s'(\mu)>0\text{ for } a<\mu<\mu^*,
		\quad
		s'(\mu)<0\text{ for } \mu^*<\mu<b.
		$$
		
		\item[(iii)]
		$s$ is concave and has a unique maximizer on $I$.
	\end{enumerate}
\end{cor}
\begin{rem}
The proof of Corollary \ref{cor-continuous-scale-graph} follows immediately from either Theorem \ref{thm-sufficient-iden-K(S)} or Theorem \ref{thm-sufficient-iden}. Corollary \ref{thm-fiber-condition} follows directly from Theorem \ref{thm-sufficient-iden}. Corollary \ref{prop:iden-graph-over-location} is proved by verifying the condition in Theorem \ref{thm:complete-iden-chacterization}, which is postponed to Appendix \ref{sec:app-section5}. The cases given in Corollary \ref{prop:iden-graph-over-location} also show that the sufficient conditions given in Theorem \ref{thm:thin-scale-proj}, Theorem \ref{thm-sufficient-iden-K(S)} and Theorem \ref{thm-sufficient-iden} are not necessary. 
Some examples are also illustrated in Figure \ref{fig:single graph}.
\end{rem}
\begin{rem}\label{rem:concrete-example-location-graph}
Corollary \ref{prop:iden-graph-over-location} applies to every positive continuous monotone function on $I$
(including affine, exponential, and logistic functions), as well as every
polynomial of degree at most two that is positive on $I$.
Further examples include:
$
s(\mu)=c\pm A|\mu-m|^p,
$
$
s(\mu)=c+A e^{-\lambda|\mu-m|^p},
$
$
s(\mu)=c+\frac{A}{1+\lambda(\mu-m)^2},
$
where $A,\lambda,p>0$, $m\in(a,b)$, and $c$ is chosen so that
$s>0$ on $I$. The plus-power family satisfies \textup{(i)}, while the
minus-power and bell-shaped families satisfy \textup{(ii)}.
On $I=[0,1]$, the beta-type function
$
s(\mu)=c+A\mu^\alpha(1-\mu)^\beta,
\; c,A,\alpha,\beta>0,
$
also satisfy \textup{(ii)}, with unique maximizer
$\mu^*=\alpha/(\alpha+\beta)$.
\end{rem}

\subsubsection{U-shape curves}
As an application of Theorem \ref{thm:complete-iden-chacterization}, we obtain the identifiability result for U-shaped supports. The proof is given in Appendix \ref{sec:app-section5}.
\begin{cor}\label{cor:U-shape}
	Let $0<\sigma_0<\sigma_1$, and let
	$a,b\in C([\sigma_0,\sigma_1])$ satisfy
	$$
	a(\sigma_0)\le b(\sigma_0),
	\qquad
	a(\sigma)<b(\sigma),
	\qquad \forall \sigma\in(\sigma_0,\sigma_1].
	$$
	Suppose that, for every $\sigma\in(\sigma_0,\sigma_1]$,
	there exist constants $C_\sigma\ge0$ and
	$\delta_\sigma\in(0,\sigma-\sigma_0)$ such that
	\begin{equation}\label{eq-U-shape-regular-branch}
	a(\sigma-h)-a(\sigma)\ge-C_\sigma\sqrt h,
	\qquad
	b(\sigma-h)-b(\sigma)\le C_\sigma\sqrt h,
	\qquad
	\forall\; 0<h<\delta_\sigma.
	\end{equation}
    	Define
	$$
	\Gamma_a
	=
	\{(a(\sigma),\sigma):
	\sigma_0\le\sigma\le\sigma_1\},
	\quad
	\Gamma_b
	=
	\{(b(\sigma),\sigma):
	\sigma_0\le\sigma\le\sigma_1\},\quad
    D_0
	=
	[a(\sigma_0),b(\sigma_0)]\times\{\sigma_0\},
	$$
	and set
	$$
	S=\Gamma_a\cup\Gamma_b\cup D_0.
	$$
    Then $\mathcal P(S)$ is identifiable.

    	Furthermore, the conclusion holds if
	$a,b\in C^{0,\alpha}([\sigma_0,\sigma_1])$
	for some $\alpha\geq \frac 12$.
\end{cor}

In Corollary \ref{cor:U-shape}, absence of the upper cap $(a(\sigma_1),b(\sigma_1))\times\{\sigma_1\}$ is necessary for identifiability of U-shape supports.  
    In fact, non-identifiability will follow if one adds just one more point of the cap $(a(\sigma_1),b(\sigma_1))\times\{\sigma_1\}$ into the support. Indeed, applying Lemma \ref{lem:Brownian-stop-rep} (as well as Remark \ref{rem:non-identiable}), we obtain the following corollary.
\begin{cor}
[U-shape plus one source point]\label{cor:domain-exit-balayage}
   Let
	$
	\Omega
	$
	be a bounded open set having closure in $\mathbb R\times (0,\infty)$, and $a=(\mu_0,\sigma_0)\in\Omega$.
    Let $U=\partial \Omega\cap\{\sigma\leq\sigma_0\}$, and $S=U\cup\{a\}$. Then $\mathcal P(S)$ is not identifiable.
\end{cor}    

\medskip

The following result is interesting because it exhibits a sharp support threshold for identifiability with respect to the parabolic U-shape. The proof is postponed to Appendix \ref{sec:app-section5}.

 

\begin{cor}[Sharp threshold for the parabolic U-shape]
	\label{cor:U-shape-sharp-threshold}
Fix $\mu_0\in\mathbb R$, $\sigma_0>0$, $\varepsilon>0$,
and $d\in(-\varepsilon,\varepsilon)$. Define
$$
C_\varepsilon
=
\{(\mu_0+y,\sigma_0+y^2):|y|\le\varepsilon\},
\quad
a_\varepsilon
=
(\mu_0+d,\sigma_0+\varepsilon^2),
\quad\text{and}\quad S=C_\eps\cup\{a_\eps\}.
$$
Then:
\begin{enumerate}
\item [(i)] $\mathcal P(S)$ is not identifiable.
\item[(ii)]
For every proper compact subset
$S' \subsetneq S$,
$\mathcal P(S')$ is identifiable.

\end{enumerate}
\end{cor}

Corollary \ref{cor:U-shape-sharp-threshold} can be generalized as following results.
\begin{cor}[One-point source and a full-support target]
	\label{cor:one-point-full-support-target}
	Let
	$
	a\in\mathbb R\times(0,\infty),
	\;
	\Gamma\subset\mathbb R\times (0,\infty)
	$
	be compact with
	$
	a\notin\Gamma.
	$
	Suppose that $\mathcal P(\Gamma)$ is identifiable and that there exists
	a probability measure
	$
	G_a\in\mathcal P(\Gamma)
	$
	such that
	\begin{equation}
		\label{eq:point-rep-measure}
		p_{G_a}
		=
		p_{\delta_a}
\quad
\text{and} 
\quad 
		\operatorname{supp}G_a
		=
		\Gamma.
	\end{equation}
	Then
	\begin{itemize}
\item[(i)]	
$
		\mathcal P(\{a\}\cup\Gamma)$
is not identifiable.
\item[(ii)]
For every proper compact subset
	$
	\Gamma'\subsetneq\Gamma,
	$
	one has
$			\mathcal P(\{a\}\cup\Gamma')$
 is identifiable.
\end{itemize}
\end{cor}

\medskip

\subsubsection{Convex curves}
The U-shaped examples show that adding a single point in
the open upper cap can destroy identifiability. Completing
the curve, however, need not have the same effect: a lower
semicircle can be completed to an identifiable full circle
by adjoining the upper semicircle, which avoids the open
horizontal cap. This contrast leads us to examine how
the geometry of a closed boundary governs cancellation.

For closed convex curves, every exit location remains on the support,
so by Theorem \ref{thm:complete-iden-chacterization} (ii), identifiability is determined by whether the backward
space-time Brownian motion exits immediately from every boundary
point. The following family exhibits a sharp transition:
even among analytic, strictly convex curves, sufficient
flatness at the uppermost point permits nontrivial cancellation.

\begin{cor}[A sharp parabolic threshold]
	\label{thm:pq-threshold}
	Let
	$
	p,q\ge1,
	\;
	s_0>1,
	$
	and define
	$$
	\Omega_{p,q}
	:=
	\left\{
	(\mu,\sigma)\in\mathbb R\times(0,\infty):
	|\mu|^p+|\sigma-s_0|^q<1
	\right\}.
	$$
	Set
	$$
	S_{p,q}
	:=
	\partial\Omega_{p,q}
	=
	\left\{
	(\mu,\sigma)\in\mathbb R\times(0,\infty):
	|\mu|^p+|\sigma-s_0|^q=1
	\right\}.
	$$
	Then
	\begin{equation}
		\label{eq:pq-threshold}
			\mathcal P(S_{p,q})
			\text{ is identifiable}
			\quad\Longleftrightarrow\quad
			1\le p\le2.
	\end{equation}
\end{cor}
The proof is given in Appendix \ref{sec:app-section5}.

\subsubsection{Non-identifiable three-branch structure} 

The previous examples of non-identifiability share a common feature: the support set $S$ contains a particular point or points which act as a source for which one may find an alternative compensating measure of sources which produce the same Gaussian mixture output
(cf. Corollary \ref{cor:domain-exit-balayage}, Corollary \ref{cor:U-shape-sharp-threshold} (i), Corollary \ref{cor:one-point-full-support-target} (i) and Corollary \ref{thm:pq-threshold}). Theorem \ref{thm-null-space-characterization} reveals a general description for these ``culprits": all culprit points must be in the Brownian interior $O(\mathcal K(S))$ while \emph{also} belonging to $S$. The redistributing sources are sent by the backward Brownian motion to land on the portion of boundary of the support hull that belongs to $S$ \emph{and} away from the Brownian interior. If the set $S$ is also contained on the boundary of the support hull, for the culprit points to be in the Brownian interior of the support hull while remaining on $S$ means that these points must carry some sort of geometric irregularity in their position relative to $S$, e.g., isolatedness or flatness, so that the simulated backward Brownian process initiated from such points must remain locally inside the support hull $\mathcal K(S)$ within a positive time. 

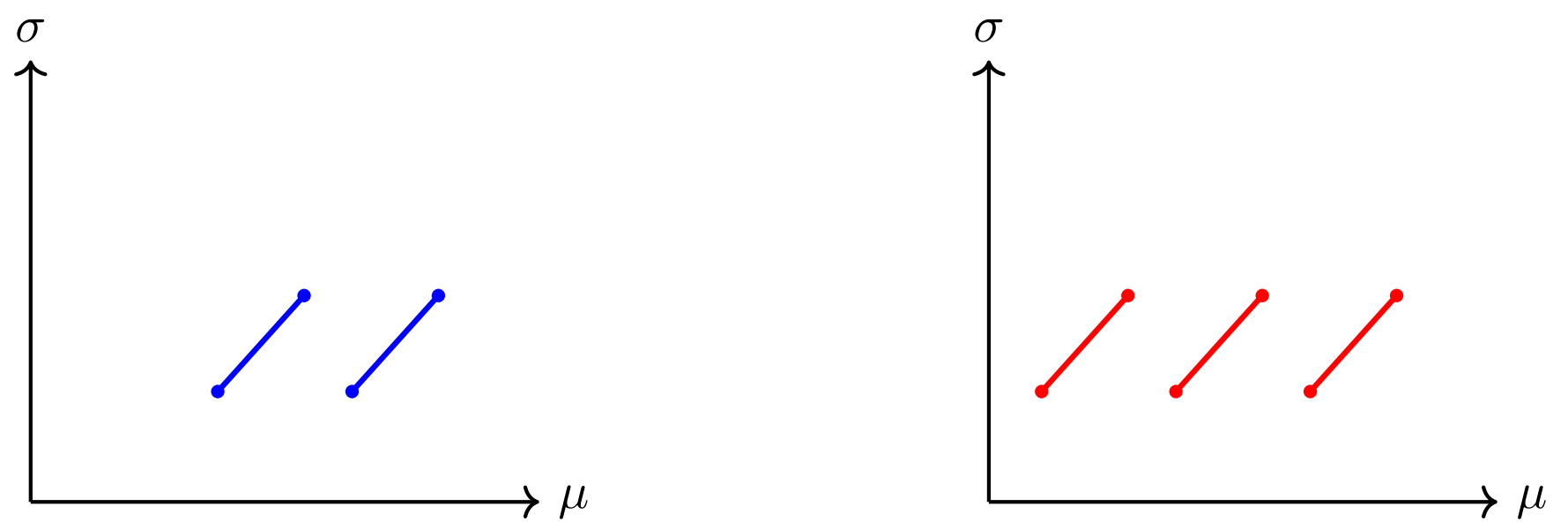
\begin{figure}[H]
	\centering
	\captionsetup[subfigure]{justification=centering}
	\begin{subfigure}[t]{0.32\textwidth}
		\centering
		\begin{tikzpicture}[scale=0.72]
			\draw[->, thick] (0,0) -- (5.3,0)
			node[right] {$\mu$};
			\draw[->, thick] (0,0) -- (0,4.6)
			node[above] {$\sigma$};
			
			
			\draw[blue, very thick]
			(1.95,1.15) -- (2.85,2.15);
			
			\draw[blue, very thick]
			(3.35,1.15) -- (4.25,2.15);
			
			
			\fill[blue] (1.95,1.15) circle (2pt);
			\fill[blue] (2.85,2.15) circle (2pt);
			
			\fill[blue] (3.35,1.15) circle (2pt);
			\fill[blue] (4.25,2.15) circle (2pt);
		\end{tikzpicture}
		\caption{Identifiable, by Corollary \ref{cor:U-shape}:
        A support consisting of two parallel line segments. }
	\end{subfigure}
\hspace{0.1\textwidth}
	\begin{subfigure}[t]{0.32\textwidth}
    
		\centering
		
		\begin{tikzpicture}[scale=0.72]
			\draw[->, thick] (0,0) -- (5.3,0)
			node[right] {$\mu$};
			\draw[->, thick] (0,0) -- (0,4.6)
			node[above] {$\sigma$};
			
			\draw[red, very thick]
			(0.55,1.15) -- (1.45,2.15);
			
			\draw[red, very thick]
			(1.95,1.15) -- (2.85,2.15);
			
			\draw[red, very thick]
			(3.35,1.15) -- (4.25,2.15);
			
			\fill[red] (0.55,1.15) circle (2pt);
			\fill[red] (1.45,2.15) circle (2pt);
			
			\fill[red] (1.95,1.15) circle (2pt);
			\fill[red] (2.85,2.15) circle (2pt);
			
			\fill[red] (3.35,1.15) circle (2pt);
			\fill[red] (4.25,2.15) circle (2pt);
		\end{tikzpicture}
		\caption{Non-identifiable, by Proposition \ref{prop:three}:
        A support consisting of three parallel line segments.
     }
	\end{subfigure}
	\caption{Two branch vs three branch structure.}
	\label{fig:three-branch}
\end{figure}

The following non-identifiability example is different,
because its source of non-identifiability can be attributed to a global geometric feature of the support set. In fact, we shall construct $S$ as smooth curve indexed by the location variable, namely $S=(\mu, s(\mu))$, where $\mu$ is restricted in a compact set $I$, and $s\in C^\infty(I)$, and yet $\mathcal P(S)$ is non-identifiable. By contrast, we already showed in Corollary \ref{cor-continuous-scale-graph} that $S=(r(\sigma),\sigma)$, with $\sigma$ being restricted in a compact set $J$, is identifiable for any continuous $r$. Moreover $S=(\mu,s(\mu))$, $\mu$ being a compact set $I$ is identifiable if $s$ satisfies certain conditions on its graph. Such conditions are of a global type, as given Corollary \ref{prop:iden-graph-over-location} and \ref{cor:U-shape}, and notably share a common feature: the horizontal fiber $S_c$ contains at most two elements for every $c > 0$. Roughly speaking, the support structure contains at most two branches, a crucial feature enabling identifiability. In the example that we will present next, the horizontal fibers $S_c$ may contain up to three elements, resulting in a three-branch structure for the set $S$. In particular, a non-negligible portion of $S$ lies in the Euclidean interior of the support hull $\mathcal K(S)$. This is a crucial global feature that induces non-identifiability, allowing us to obtain an example of non-identifiable $\mathcal{P}(S)$, even when $S$ is a smooth graph over the location variable. See Figure \ref{fig:three-branch} for an illustration.

The basic structure of such type of non-identifiability is revealed by the following proposition. The proof is given in Appendix \ref{sec:app-section5}.

\begin{prop}
	\label{prop:three}
	
	Let
	$
	0<T_0<T_1<\infty,
	\;
	c\in\mathbb R\setminus\{0\},
	$
	and 
	$
	a_1<a_2<a_3
	$
	be three distinct real numbers. For $j=1,2,3$, define
	$$
	\Gamma_j
	:=
	\left\{
	(a_j+ct,t):t\in[T_0,T_1]
	\right\}
	\subset
	\mathbb R\times(0,\infty),
	$$
	and set
	$
	S:=\Gamma_1\cup\Gamma_2\cup\Gamma_3.
	$
	Then, $\mathcal{P}(\Gamma_1 \cup \Gamma_2)$ is identifiable. However,
	$\mathcal P(S)$ is not identifiable.
\end{prop}

\begin{rem}[A gap in the argument of Bruni--Koch \cite{bruni1985identifiability}]\label{rem-Bruni}
Proposition \ref{prop:three} provides a counterexample to a main result in \cite{bruni1985identifiability}. In 
\cite[Theorem 1]{bruni1985identifiability},
they claimed an identifiability result which implies that Gaussian mixtures are identifiable when $S$ is represented by a compact $C^1$ curve. This is not correct. Indeed, by choosing $a_1<a_2<a_3$, $0<c<\frac{\min\{a_2-a_1,a_3-a_2\}}{T_1-T_0}$, the set $S$ given in the above Proposition can be extended to belong to a smooth graph $S=\{(\mu, s(\mu)) | \mu\in I\}$, where $s\in C^\infty(I)$ for $I = [a_1 + c T_0, a_3 + c T_1]$. The class of Gaussian mixtures whose mixing measure is supported on such a set $S$ is therefore not identifiable.


    The induction step in the preliminary lemma of Bruni and Koch [pg. 1344]\cite{bruni1985identifiability}  requires a boundedness estimate for a quotient of normalized differences of weighted Gaussian integrals. When both the dominant scale and the dominant location coincide, the leading exponential factors cancel, and the numerator and denominator may both tend to zero. The required bound therefore depends on their relative rates of vanishing. The proof does not
    establish this comparison: convergence of the individual truncated-integral ratios alone does not control the quotient of the corresponding vanishing differences.
	Instead, it states that the relevant ratios of truncated Gaussian
	integrals are ``essentially independent'' of $x$ and converge, as
	the truncation parameter tends to zero, to a bounded quantity
	independent of the indices. 
	Equality of the
	dominant location and scale determines only the leading exponential
	rate of the corresponding Laplace-type integrals; it does not, for
	general mixing measures, control their subexponential factors or the
	rate at which their ratio approaches \(1\). Hence this is not merely
	an omitted routine estimate, but a substantive step on which the
	induction depends. 
	
	The counterexample above shows, moreover, that
	under the assumptions stated in \cite{bruni1985identifiability} this      step
	cannot in general be repaired: the claimed identifiability conclusion
	itself fails even when the parameter curve is fixed.
\end{rem}

\subsection{Canonical and identifiable representation for Gaussian mixtures}
\label{sec:canonical-rep}

The exit construction also yields a canonical representation
principle for Gaussian mixtures. Every Gaussian mixture with compact parameter support for its mixing measure can be represented by a Gaussian mixture whose mixing measure is supported on a fixed boundary curve of a rectangular or disk-shaped enclosure.
For a rectangle, the bottom and two vertical sides suffice;
for a disk, the entire circle serves this role.
In both cases, the representation is unique, and the curve
is the smallest compact subset of the enclosure generating
its full mixture class.
This result is stated in the following theorem, 
whose proof is given in Appendix \ref{sec:app-section5}.

For a compact parameter set $A \in \mathbb R \times (0,1)$, denote by
$$
\mathcal  G(A):=\biggr \{p_G \;\bigr |\; p_G(x) = \int \phi(x;\mu,\sigma) G(d\mu, d\sigma),\; G\in\mathcal P(A) \biggr \}
$$
the class of all Gaussian mixture density functions associated with some mixing probability measure whose support is contained in $A$.

\begin{thm}[Canonical and minimal boundary representations]\label{prop:rep-mixing-to-boundary}
Let $K\subset\mathbb R\times(0,\infty)$ be either a rectangle
$[a,b]\times[s_0,s_1]$, where $a<b$ and $0<s_0<s_1$,
or a closed Euclidean disk of positive radius. Define
$$
\Gamma=
\begin{cases}
([a,b]\times\{s_0\})\cup
(\{a,b\}\times[s_0,s_1]),
& \text{if $K$ is the rectangle},\\[2mm]
\partial K,
& \text{if $K$ is a disk}.
\end{cases}
$$
Then the following statements hold.
\begin{enumerate}
\item[(i)]
Every $G\in\mathcal P(K)$ admits a unique
$G_\Gamma\in\mathcal P(\Gamma)$ such that
$
p_{G_\Gamma}=p_G.
$
This representative actually is $G_\Gamma=E_KG$.
Consequently,
$$
\mathcal G(K)=\mathcal G(\Gamma).
$$

\item[(ii)]
For every compact $A\subseteq K$,
$$
\mathcal G(A)=\mathcal G(K)
\quad\Longleftrightarrow\quad
\Gamma\subseteq A.
$$
Thus $\Gamma$ is the smallest compact parameter subset
of $K$ generating the class of Gaussian mixture densities 
$\mathcal G(K)$.
\end{enumerate}
\end{thm}

The two geometries show how the direction of the backward Brownian diffusion
and the shape of the boundary determine the canonical
support. The open top edge of a rectangle consists of
active sources whose Gaussian outputs can be reproduced
on the bottom and sides, leaving a U-shaped curve as
the smallest compact generating set. For a disk, every
boundary point is an immediate-exit point, including
the uppermost point, so the entire circle is indispensable.
Thus exact reduction from a two-dimensional parameter
region to an identifiable curve can produce either a
curve with endpoints or a smooth closed curve.
In both cases, uniqueness concerns the representative
on the chosen curve; the original mixing measure
on the full region need not be identifiable.

\begin{cor}\label{cor:circle-representation}
    For every compact
$S\subset\mathbb R\times(0,\infty)$, there exists a circle
$\Gamma\subset\mathbb R\times(0,\infty)$ enclosing $S$
such that every Gaussian mixture generated by
$G\in\mathcal P(S)$ admits a unique mixing representation
on $\Gamma$. The curve $\Gamma$ depends only on $S$.
\end{cor}

The representation theorem gives curved supports a direct
modeling role. A suitable enclosing curve can preserve
the entire class of Gaussian mixture densities generated
by a two-dimensional parameter region, while providing
a unique mixing representation for each density.
This motivates a broader study of curved supports:
which geometric features ensure identifiability,
and how far does this property extend beyond the
enclosures considered above?

The uniqueness here is relative to the chosen curve.
It identifies a representative of the observed density
within that support class; the original mixing measure
in the enclosed region may remain non-identifiable.

\subsection{Identifiability under transformation}
\label{sec:transformation}
Section \ref{sec:parameterized-curve} provides several examples in which a seemingly innocent modification to the support set $S$ can drastically change its identifiability character.
This leads us to wonder how transformation of the support set affects identifiability.
Precisely, we want to find which transformation $\mathcal T$ preserves the identifiability character, namely $\mathcal P(\mathcal T(S))$ is  identifiable if and only if $\mathcal P(S)$ is.
This section is devoted to answering this question.

We first narrow the search for such transformations using the following necessary condition for identifiability invariance.



\begin{thm}
	\label{thm:uni-trans-pre}
	Let
	$$
	\mathcal T(\mu,\sigma)=(X(\mu,\sigma),Y(\mu,\sigma)):
	\mathbb R\times (0,\infty)
	\rightarrow
	\Omega
	$$
	be a $C^1$ diffeomorphism onto a domain
	$
	\Omega\subset\mathbb R\times (0,\infty).
	$
	Assume that:
	\begin{equation}
		\label{eq:univ-iden}
     \text{For every compact set }
	S\subset\mathbb R\times (0,\infty),\quad
		\mathcal P(S)\text{ identifiable}
		\quad\Longrightarrow\quad
		\mathcal P(\mathcal T(S))\text{ identifiable}.
	\end{equation}
	Then
	\begin{equation}
		\label{eq:V-independent-of-mu}
		\partial_\mu Y(\mu,\sigma)=0
		\qquad
		\text{for every }(\mu,\sigma)\in\mathbb R\times (0,\infty).
	\end{equation}
\end{thm}

\begin{proof}
	Suppose, to the contrary, that there exists
	$
	(\mu_0,\sigma_0)\in\mathbb R\times (0,\infty)
	$
	such that
$
		\partial_\mu Y(\mu_0,\sigma_0)\ne0.
$
	Write
	$$
	\mathcal T(\mu_0,\sigma_0)=(m_0,s_0)\in\Omega.
	$$
	Since $T$ is a $C^1$ diffeomorphism, write
	$$
	\mathcal T^{-1}(m,s)
	=
	\bigl(
	X^{-1}(m,s),Y^{-1}(m,s)
	\bigr).
	$$
	At $(m_0,s_0)$, the Jacobian of the inverse map $\mathcal T$ satisfies
	$
	D\mathcal T^{-1}(m_0,s_0)
	=
	[D\mathcal T(\mu_0,\sigma_0)]^{-1}.
	$
	Since
	$$
	D\mathcal T(\mu_0,\sigma_0)
	=
	\begin{pmatrix}
		X_\mu(\mu_0,\sigma_0) & X_\sigma(\mu_0,\sigma_0)
		\\
		Y_\mu(\mu_0,\sigma_0) & Y_\sigma(\mu_0,\sigma_0)
	\end{pmatrix},
	$$
	we obtain
	$$
	D\mathcal T^{-1}(m_0,s_0)
	=
	\frac{1}{\det D\mathcal T(\mu_0,\sigma_0))}
	\begin{pmatrix}
		Y_\sigma(\mu_0,\sigma_0)) & -X_\sigma(\mu_0,\sigma_0))
		\\
		-Y_\mu(\mu_0,\sigma_0)) & X_\mu(\mu_0,\sigma_0))
	\end{pmatrix}.
	$$
	Therefore
	\begin{equation}
		\label{eq:inverse-var}
		\partial_mY^{-1}(m_0,s_0)
		=
		-\frac{Y_\mu(\mu_0,\sigma_0))}{\det D\mathcal T(\mu_0,\sigma_0))}
		\ne0.
	\end{equation}
	
	Since $\Omega$ is open and $(m_0,s_0)\in\Omega$, we may choose
	$\varepsilon>0$ sufficiently small so that
	$$
	D_\varepsilon
	:=
	\left\{
	(m_0+y,s_0+y^2):
	|y|\le\varepsilon
	\right\}
	\cup
	\left\{
	(m_0,s_0+\varepsilon^2)
	\right\}
	$$
	is contained in $\Omega$.
	By Corollary \ref{cor:U-shape-sharp-threshold},
	\begin{equation}
		\label{eq:L-eps-noniden}
		\mathcal P(D_\varepsilon)
		\text{ is not identifiable}.
	\end{equation}
	
	Consider the preimage $\mathcal T^{-1}(C_\eps)$ of the parabolic arc
	$$
	C_\varepsilon
	=
	\left\{
	(m_0+y,s_0+y^2):
	|y|\le\varepsilon
	\right\}.
	$$
	It is parametrized by
	$$
	\mathcal T^{-1}(m_0+y,s_0+y^2),
	\qquad
	|y|\le\varepsilon.
	$$
	Let
	$$
	w(y)
	:=
	Y^{-1}(m_0+y,s_0+y^2)
	$$
	be its scale coordinate. Then
	$$
	w'(y)
	=
	\partial_mY^{-1}(m_0+y,s_0+y^2)
	+
	2y\partial_sY^{-1}(m_0+y,s_0+y^2).
	$$
	In particular, by \eqref{eq:inverse-var}
	$$
	w'(0)=\partial_mY^{-1}(m_0,s_0)\ne0.
	$$
	By continuity, after decreasing $\varepsilon$ if necessary,
	$
	w'(y)\ne0
	\;
	\text{for every }|y|\le\varepsilon,
	$
	and $w'$ has constant sign on this interval. Hence
	$
	w:[-\varepsilon,\varepsilon]\longrightarrow\mathbb R
	$
	is strictly monotone.
	It follows that the compact curve
	$
	\mathcal T^{-1}(C_\varepsilon)
	$
	intersects every horizontal line
	$
	\mathbb R\times\{c\}
	$
	in at most one point.
	
	Let
	$$
	S_\varepsilon
	:=
	\mathcal T^{-1}(D_\varepsilon)
	=
	\mathcal T^{-1}(C_\varepsilon)
	\cup
	\mathcal T^{-1}(m_0,s_0+\varepsilon^2).
	$$
	Since $\mathcal T^{-1}$ is continuous and $D_\varepsilon$ is compact,
	$S_\varepsilon$ is compact.
	
	For every
	$
	c\ne \pi_\sigma(\mathcal T^{-1}(m_0,s_0+\varepsilon^2)),
	$
	the fiber
	$
	(S_\varepsilon)_c
	$
	contains at most one point. Hence
	$
	\mathbb R\setminus(S_\varepsilon)_c
	$
	has no bounded connected component.
	Therefore condition \eqref{eq-sufficient-cond} holds for $S_\varepsilon$.
By Theorem \ref{thm-sufficient-iden},
		$\mathcal P(S_\varepsilon)$ is identifiable.
	But
	$
	\mathcal T(S_\varepsilon)=D_\varepsilon,
	$
	so \eqref{eq:L-eps-noniden} contradicts the invariance condition
	\eqref{eq:univ-iden}.
	Hence
	$
	\partial_\mu Y\equiv0
	$
	on $\mathbb R\times (0,\infty)$.
    \end{proof}

We next describe several classes of transformations of the
location-scale parameter space under which identifiability is invariant.

\begin{prop}[Affine-parabolic transformations]
	\label{prop:affine}
	Let
	$
	S\subset\mathbb R\times(0,\infty)
	$
	be compact, and let
	$$
	\mathcal A_{a,b,c}(\mu,\sigma)
	=
	(a\mu+b,a^2\sigma+c),\;
\text{where }
	a\ne0,\;
	a^2\sigma+c>0
	\;
	\text{for every }(\mu,\sigma)\in S.
$$
	Then
	$$
	\mathcal P(S)\text{ is identifiable}
	\quad\Longleftrightarrow\quad
	\mathcal P(\mathcal A_{a,b,c}(S))\text{ is identifiable}.
	$$
\end{prop}

\begin{proof}
	The map $\mathcal A_{a,b,c}$ is a homeomorphism from $S$ onto its image.
	Hence
	$
	G\longmapsto (\mathcal A_{a,b,c})_\#G
	$
	is a bijection from $\mathcal P(S)$ onto
	$\mathcal P(\mathcal A_{a,b,c}(S))$.
Consider Fourier transformation
	$$
	\widehat p_G(t)
	=
	\int_{\mathbb R}e^{itx}p_G(x)\,dx
	=
	\int_S
	e^{it\mu-\sigma t^2/2}
	\,dG(\mu,\sigma).
	$$
	Therefore
	\begin{align}\label{eq:affine-fourier}
		\widehat p_{(\mathcal A_{a,b,c})_\#G}(t)
		&=
		\int_S
		e^{it(a\mu+b)-(a^2\sigma+c)t^2/2}
		\,dG(\mu,\sigma)
		=
		e^{ibt-ct^2/2}
		\widehat p_G(at).
	\end{align}
	
	Suppose first that
	$
	p_{G_1}=p_{G_2}.
	$
	Then
	$
	\widehat p_{G_1}=\widehat p_{G_2},
	$
	and \eqref{eq:affine-fourier} yields
	$
	\widehat p_{(\mathcal A_{a,b,c})_\#G_1}
	=
	\widehat p_{(\mathcal A_{a,b,c})_\#G_2}.
	$
	By uniqueness of the Fourier transform,
	the corresponding mixture distributions are equal. Since
	$\mathcal A_{a,b,c}(S)$ is compact in
	$\mathbb R\times(0,\infty)$, its scales are bounded away from zero,
	so the two mixture densities are continuous. Hence they are equal
	everywhere.
	
	Conversely, suppose
	$
	p_{(\mathcal A_{a,b,c})_\#G_1}
	=
	p_{(\mathcal A_{a,b,c})_\#G_2}.
	$
	Then, by \eqref{eq:affine-fourier},
	$$
	e^{ibt-ct^2/2}
	\widehat p_{G_1}(at)
	=
	e^{ibt-ct^2/2}
	\widehat p_{G_2}(at),\text{ for every }t\in\mathbb R.
	$$
	Since
	$
	e^{ibt-ct^2/2}\ne0
	$
	and $a\ne0$, we obtain
	$
	\widehat p_{G_1}(s)
	=
	\widehat p_{G_2}(s)
	\;
	\text{for every }s\in\mathbb R.
	$
	Again by Fourier uniqueness and continuity of the mixture densities,
	$
	p_{G_1}=p_{G_2}.
	$
	
	Thus, the equivalence is proved.
\end{proof}

\begin{prop}
	\label{prop:quadratic}
	Let
	$
	S\subset\mathbb R\times(0,\infty)
	$
	be compact. Let $r,q\in\mathbb R$ satisfy
	$
	1+q\sigma>0
	\;
	\text{for every }(\mu,\sigma)\in S,
	$
	and define
	$$
	Q_{r,q}(\mu,\sigma)
	=
	\left(
	\frac{\mu+r\sigma}{1+q\sigma},
	\frac{\sigma}{1+q\sigma}
	\right).
	$$
	Then
	$$
	\mathcal P(S)\text{ is identifiable}
	\quad\Longleftrightarrow\quad
	\mathcal P(Q_{r,q}(S))\text{ is identifiable}.
	$$
\end{prop}

\begin{proof}
	A direct completion of the square gives
	\begin{equation}
		\label{eq:quadratic-kernel}
		h_{r,q}(x)\phi(x;\mu,\sigma)
		=
		Z_{r,q}(\mu,\sigma)
		\phi\left(
		x;
		Q_{r,q}(\mu,\sigma)
		\right),
	\end{equation}
where
	$$
	h_{r,q}(x)
	=
	\exp\left(
	rx-\frac q2x^2
	\right)\quad
	\text{and}\quad
	Z_{r,q}(\mu,\sigma)
	=
	\frac{1}{\sqrt{1+q\sigma}}
	\exp\left(
	\frac{2r\mu+r^2\sigma-q\mu^2}{2(1+q\sigma)}
	\right).
	$$    
	Since $S$ is compact and
	$
	1+q\sigma>0
	$
	on $S$, the function $Z_{r,q}$ is continuous and $
	0<
	\min_S Z_{r,q}
	\le
	\max_S Z_{r,q}
	<
	\infty.
	$
	
	For
	$
	G\in\mathcal P(S),
	$
	define
	$$
	m_G
	:=
	\int_S
	Z_{r,q}(\mu,\sigma)
	\,dG(\mu,\sigma)\quad
\text{and}\quad
	dG^Z(\mu,\sigma)
	=
	\frac{Z_{r,q}(\mu,\sigma)}{m_G}
	\,dG(\mu,\sigma).
	$$
	Then $G^Z\in\mathcal P(S)$. Define
	$$
	\Psi(G)
	:=
	(Q_{r,q})_\#G^Z.
	$$ Then, $\Psi(G) \in \mathcal P(Q_{r,q}(S))$. Moreover, notice that
    $$p_{\Psi(G)}(x) = \int \phi(x; Q_{r,q}(\mu,\sigma) d G^Z(\mu,\sigma)
    = \frac{1}{m_G} \int \phi(x;Q_{r,q}(\mu,\sigma) Z_{r,q}(\mu,\sigma) dG(\mu,\sigma).$$
Hence, integrating \eqref{eq:quadratic-kernel} against $G$ and $x$ yields
	\begin{equation}
		\label{eq:quadratic-mix}
		p_{\Psi(G)}(x)
		=
		\frac{h_{r,q}(x)}{m_G}
		p_G(x),\qquad\text{and}\qquad		m_G
		=
		\int_{\mathbb R}
		h_{r,q}(x)p_G(x)\,dx.
	\end{equation}
	
	We next verify that $\Psi$ is a bijection. The inverse of $Q_{r,q}$
	on its image is
	$
	Q_{r,q}^{-1}=Q_{-r,-q}.
	$
	Let
	$
	G_0\in\mathcal P(Q_{r,q}(S))
	$
	and set
	$
	\widetilde G_0
	:=
	(Q_{r,q}^{-1})_\#G_0.
	$
	Define
	$
	C_{G_0}
	:=
	\int_S
	(1/Z_{r,q})\,d\widetilde G_0.
	$
	Since $Z_{r,q}$ is bounded above and below by positive constants on
	$S$,
	$
	0<C_{G_0}<\infty.
	$
	Define
	$
	dG
	=
	\frac{(1/Z_{r,q})}{C_{G_0}}
	\,d\widetilde G_0.
	$
	Then
	$
	G\in\mathcal P(S)
	$
	and
	$
	m_G
	=
	\int_S Z_{r,q}\,dG
	=
	\frac1{C_{G_0}}.
	$
	Hence
	$
	G^Z=\widetilde G_0,
	$
	and therefore
	$
	\Psi(G)=G_0.
	$
	Thus $\Psi$ is surjective. The same calculation shows uniqueness of
	the preimage, so $\Psi$ is bijective.

	Suppose
	$
	p_{G_1}=p_{G_2}.
	$
	Then 
     \eqref{eq:quadratic-mix}
    implies
	$
	m_{G_1}=m_{G_2},
	$
    and then,
	$
	p_{\Psi(G_1)}
	=
	p_{\Psi(G_2)}.
	$
	
	Conversely, suppose
	$
	p_{\Psi(G_1)}
	=
	p_{\Psi(G_2)}.
	$
	Since
	$
	h_{r,q}(x)>0
	\;
	\text{for every }x,
	$
	equation \eqref{eq:quadratic-mix} gives
	$
	\frac{p_{G_1}(x)}{m_{G_1}}
	=
	\frac{p_{G_2}(x)}{m_{G_2}}
	\;
	\text{for every }x.
	$
	Integrating over $\mathbb R$ and using
	$
	\int_{\mathbb R}p_{G_j}(x)\,dx=1
	$
	gives
	$
	m_{G_1}=m_{G_2}.
	$
	Therefore
	$
	p_{G_1}=p_{G_2}.
	$
	
	Thus $\Psi$ is a bijection between the two mixing classes and preserves
	the equivalence relation of producing the same Gaussian mixture in
	both directions. This proves the result.
\end{proof}

\begin{prop}
	\label{thm:gauss}
	Let
	$
	S\subset\mathbb R\times(0,\infty)
	$
	be compact. Consider
	\begin{equation*}
		\mathcal T(\mu,\sigma)
		=
		\left(
		a\frac{\mu+r\sigma}{1+q\sigma}+b,
		\,
		a^2\frac{\sigma}{1+q\sigma}+c
		\right),\quad
\text{where }
	a,b,c,r,q\in\mathbb R,
	\;
	a\ne0.
\end{equation*}
	Assume that
	$$
	1+q\sigma>0
\quad
\text{and}
\quad
	a^2\frac{\sigma}{1+q\sigma}+c>0,\qquad
\text{for every }
	(\mu,\sigma)\in S.
	$$
	Then
	$$
	\mathcal P(S)\text{ is identifiable}
	\quad\Longleftrightarrow\quad
	\mathcal P(\mathcal T(S))\text{ is identifiable}.
	$$
\end{prop}

\begin{proof}
	Write
	$$
	\mathcal T=\mathcal A_{a,b,c}\circ Q_{r,q},
	$$
	where
	$$
	Q_{r,q}(\mu,\sigma)
	=
	\left(
	\frac{\mu+r\sigma}{1+q\sigma},
	\frac{\sigma}{1+q\sigma}
	\right)
\quad
\text{and}
	\quad
	\mathcal A_{a,b,c}(u,v)
	=
	(au+b,a^2v+c).
	$$
	The assumptions ensure that both transformations are admissible on
	the relevant compact sets. The conclusion follows successively from
	Propositions \ref{prop:quadratic}
	and \ref{prop:affine}.
\end{proof}

\begin{thm}[Fractional-linear Gaussian--heat transformations]
	\label{thm:mobius}
	Let
	$
	S\subset\mathbb R\times(0,\infty)
	$
	be compact. Consider
	\begin{equation}
		\label{eq:mobius-gaussian-heat-map}
		\mathcal T(\mu,\sigma)
		=
		\left(
		\frac{\lambda\mu+\eta\sigma+\xi}
		{\gamma\sigma+\delta},
		\,
		\frac{\alpha\sigma+\beta}
		{\gamma\sigma+\delta}
		\right),
	\end{equation}
	where
	$
	\alpha,\beta,\gamma,\delta,
	\lambda,\eta,\xi\in\mathbb R
	$
	satisfy
	\begin{equation}
		\label{eq:mobius-det}
		\alpha\delta-\beta\gamma
		=
		\lambda^2
		>
		0.
	\end{equation}
	Assume that
	\begin{equation}
		\label{eq:mobius-admis}
		\gamma\sigma+\delta\ne0,
		\qquad
		\frac{\alpha\sigma+\beta}
		{\gamma\sigma+\delta}>0,\quad
			\text{for every }
		\sigma\in\pi_\sigma(S).
	\end{equation}

	Then
	\begin{equation}
		\label{eq:mobius-iden}
		\mathcal P(S)\text{ is identifiable}
		\quad\Longleftrightarrow\quad
		\mathcal P(\mathcal T(S))\text{ is identifiable}.
	\end{equation}
\end{thm}

\begin{proof}
	Set
	$$
	d(\sigma):=\gamma\sigma+\delta,
	\quad
	\text{and}
	\quad
	v(\sigma)
	:=
	\frac{\alpha\sigma+\beta}
	{\gamma\sigma+\delta}.
	$$
	
	We first consider the case where $d(\sigma)$ has a constant sign on $\pi_\sigma(S)$.
	Choose
	$
	s_0<\min \pi_\sigma(S)
	$
	such that
	$
	d_0:=\gamma s_0+\delta
	$
	has the same sign as $d(\sigma)$ for $\sigma\in E$.
	Define
	$
	\mathcal A_{s_0}(\mu,\sigma)
	=
	(\mu,\sigma-s_0).
	$
	By Proposition \ref{prop:affine},
	$$
	\mathcal P(S)\text{ is identifiable}
	\quad 
\Longleftrightarrow
	\quad
	\mathcal P(\mathcal A_{s_0}(S))\text{ is identifiable}.
	$$
Let $\tau=\sigma-s_0$, direct calculations show
	\begin{equation*}
		T(\mu,\tau+s_0)
		=
		\left(
		a\frac{\mu+r\tau}{1+q\tau}+b,
		\,
		a^2\frac{\tau}{1+q\tau}+c
		\right),
	\end{equation*}
    for
	$$
	a:=\frac{\lambda}{d_0},
	\;
	q:=\frac{\gamma}{d_0},
	\;
	b:=\frac{\eta s_0+\xi}{d_0},
	\;
	c:=\frac{\alpha s_0+\beta}{d_0},
	\;
	r
	:=
	\frac{\eta/d_0-bq}{a}.
	$$
It is easily seen that $a\neq 0$,
and $1+q\tau
	=
	\frac{\gamma\sigma+\delta}{d_0}
	>0
	$
	on $S$, since numerator and denominator have the same sign, and $a^2\frac{\tau}{1+q\tau}+c=v(\sigma)>0$.
    
    	Therefore
	Proposition \ref{thm:gauss}
	applied to $\mathcal A_{s_0}(S)$ yields
	$$
	\mathcal P(\mathcal A_{s_0}(S))\text{ is identifiable}
\quad
\Longleftrightarrow
\quad
	\mathcal P(T(S))\text{ is identifiable}.
	$$
	Combining this with the invariance under $\mathcal A_{s_0}$ proves
	\eqref{eq:mobius-iden} whenever $d$ has constant
	sign on $\pi_\sigma(S)$.

    	We now consider the case where $d$ changes sign on $\pi_\sigma(S)$. Then
	$
	\gamma\ne0,
	$
	and the pole
	$
	\sigma_0:=-\frac{\delta}{\gamma}
	$
	does not belong to $\pi_\sigma(S)$ by
	\eqref{eq:mobius-admis}.
	Define
	$$
	S_-
	:=
	S\cap\{\sigma<\sigma_0\},
	\quad
\text{and}
\quad
	S_+
	:=
	S\cap\{\sigma>\sigma_0\}.
	$$
	If either $S_-$ or $S_+$ is empty, then $d$ has constant sign on $E$
	and the preceding argument applies. Thus, suppose both are nonempty.
    	Since $\pi_\sigma(S)$ is compact and
	$
	\sigma_0\notin \pi_\sigma(S),
	$
	there is a positive scale gap between $S_-$ and $S_+$, i.e.,
	$
	\max\pi_\sigma(S_-)
	<
	\min\pi_\sigma(S_+).
	$
	Corollary \ref{cor:var-gap} therefore gives
	$$
	\mathcal P(S)\text{ is identifiable}
	\quad
	\Longleftrightarrow
    \quad
    \text{both }
	\mathcal P(S_-)
	\text{ and }
	\mathcal P(S_+)
    \text{ are identifiable.}
	$$
	On each of $S_-$ and $S_+$, the denominator
	$
	\gamma\sigma+\delta
	$
	has constant sign. Therefore the first part of the proof gives
	$$
	\mathcal P(S_-)\text{ is identifiable}
	\quad\Longleftrightarrow\quad
	\mathcal P(\mathcal T(S_-))\text{ is identifiable},
$$
and
$$
	\mathcal P(S_+)\text{ is identifiable}
	\quad\Longleftrightarrow\quad
	\mathcal P(\mathcal T(S_+))\text{ is identifiable}.
	$$
Therefore, it remains to verify that the two image pieces are separated in the
	scale coordinate. A direct computation using
	\eqref{eq:mobius-det} gives
	\begin{align}\label{eq:mobius-image-gap}
		v(\sigma)-\frac{\alpha}{\gamma}
		=
		\frac{
			\gamma(\alpha\sigma+\beta)
			-
			\alpha(\gamma\sigma+\delta)
		}{
			\gamma(\gamma\sigma+\delta)
		}
		=
		-\frac{\lambda^2}
		{\gamma(\gamma\sigma+\delta)}.
	\end{align}
	Since
	$
	\gamma(\gamma\sigma+\delta)
	=
	\gamma^2(\sigma-\sigma_0),
	$
	we obtain
	$
	v(\sigma)>\frac{\alpha}{\gamma}
	\;
	\text{when }\sigma<\sigma_0,
	$
	and
	$
	v(\sigma)<\frac{\alpha}{\gamma}
	\;
	\text{when }\sigma>\sigma_0.
	$
	Thus
	$
	\pi_\sigma(\mathcal T(S_-))
	$
	lies strictly above $\alpha/\gamma$, whereas
	$
	\pi_\sigma(\mathcal T(S_+))
	$
	lies strictly below $\alpha/\gamma$.
	These are disjoint compact subsets of $(0,\infty)$ and are therefore
	separated by a positive distance. In particular, 
	Proposition \ref{cor:var-gap} yields
	$$
	\mathcal P(\mathcal T(S))\text{ is identifiable}
	\quad
\Longleftrightarrow
	\quad
    \text{both }
	\mathcal P(\mathcal T(S_-))
	\text{ and }
	\mathcal P(\mathcal T(S_+))
    \text{ are identifiable}.
	$$
	Combining with the preceding equivalences proves
	\eqref{eq:mobius-iden}.
\end{proof}

\begin{rem}
The transformation \eqref{eq:mobius-gaussian-heat-map}
can 
change the geometric
form of $S$. 
For example, if
$$
S=\{(\mu,s(\mu)):\mu\in I\},
\qquad
s(\mu)=a\mu^2+b\mu+c,
\qquad a\neq0,
$$
then $\mathcal T(S)$ is contained in a nondegenerate conic. If
$\gamma=0$, this conic is again a parabola. If $\gamma\neq0$, its type
is determined by
$
\Delta_\infty
=
b^2-4a\left(c+\frac{\delta}{\gamma}\right).
$
More precisely,
$
\Delta_\infty<0,\;
\Delta_\infty=0,\;
\Delta_\infty>0
$
correspond respectively to an ellipse, a parabola, and a hyperbola.
Thus, although $\mathcal T$ does not enlarge the class of sets satisfying
\eqref{eq-sufficient-cond}, it can turn quadratic graphs into geometrically quite
different conic arcs.
\end{rem}

\appendix

\section{Proofs of Some Lemmas and Propositions}
\label{app:proofs}
\subsection{Proof of Proposition \ref{thm-ness-cond}}
\begin{prop}
\label{thm-ness-cond}
	Let
	$
	S\subset \mathbb{R}\times(0,\infty)
	$
	be compact. Suppose that the class $\mathcal P(S)$ is identifiable.
	Then
	$
	\operatorname{Int}(S)=\emptyset.
	$
	Equivalently, if
	$
	\operatorname{Int}(S)\neq\emptyset,
	$
	then $\mathcal P(S)$ is not identifiable.
\end{prop}
\begin{proof}[Proof of Proposition \ref{thm-ness-cond}]
The proof is similar to \cite{le2026partial}.
	Assume, toward a contradiction, that
	$
	\operatorname{Int}(S)\neq\emptyset.
	$
	Then there exists a nonzero function
	$
	q\in C_c^\infty(\operatorname{int}(S))
	$
	such that
	$$
	h(\mu,\sigma)
	:=
	\frac{\partial q}{\partial \sigma}(\mu,\sigma)
	+
	\frac12
	\frac{\partial^2 q}{\partial\mu^2}(\mu,\sigma)
	\not\equiv0.
	$$
	Define a finite signed measure $\nu$ on $S$ by
	$$
	\nu(d\mu,d\sigma)
	=
	h(\mu,\sigma)\,d\mu\,d\sigma.
	$$
	
	We first show that $\nu$ has total mass zero. Since $q$ is compactly
	supported in $\operatorname{int}(S)$,
	$$
	\nu(S)
	=
	\int_S
	\left(
	\frac{\partial q}{\partial \sigma}
	+
	\frac12
	\frac{\partial^2 q}{\partial\mu^2}
	\right)
	\,d\mu\,d\sigma
	=
	0,
	$$
	because both terms are total derivatives and all boundary terms vanish.
	
	Next, we show that $\nu$ annihilates every Gaussian observation
	function. 
	For every fixed $x\in\mathbb R$, integration by parts gives
	\begin{align*}
		\int_S
		\phi(x;\mu,\sigma)\,\nu(d\mu,d\sigma)
		&=
		\int_S
		\phi(x;\mu,\sigma)
		\left(
		q_\sigma(\mu,\sigma)
		+
		\frac12 q_{\mu\mu}(\mu,\sigma)
		\right)
		\,d\mu\,d\sigma
		\\
		&=
		-
		\int_S
		q(\mu,\sigma)
		\phi_\sigma(x;\mu,\sigma)
		\,d\mu\,d\sigma
		+
		\frac12
		\int_S
		q(\mu,\sigma)
		\phi_{\mu\mu}(x;\mu,\sigma)
		\,d\mu\,d\sigma.
	\end{align*}
	Again, no boundary terms appear because
	$
	q\in C_c^\infty(\operatorname{int}(S)).
	$
	Using the heat equation, we obtain
	\begin{equation}\label{eq-phinu}
		\int_S
		\phi(x;\mu,\sigma)\,\nu(d\mu,d\sigma)
		=
		0,
		\qquad
		\forall x\in\mathbb R.
	\end{equation}
	
	Since $h\not\equiv0$, the signed measure $\nu$ is nonzero.
	Let
	$
	\nu=\nu^+-\nu^-
	$
	be its Jordan decomposition. Since $\nu(S)=0$,
	$
	\nu^+(S)=\nu^-(S).
	$
	Since $\nu\neq0$, this common mass is strictly positive. Set
	$
	M:=\nu^+(S)=\nu^-(S)>0
	$
	and define
	$
	G_1:=\frac{\nu^+}{M},
	\;
	G_2:=\frac{\nu^-}{M}.
	$
	Then
	$
	G_1,G_2\in\mathcal P(S),
	\;
	G_1\neq G_2.
	$
	
	On the other hand, for every $x\in\mathbb R$, using \eqref{eq-phinu},
	\begin{align*}
		p_{G_1}(x)-p_{G_2}(x)
		&=
		\int_S
		\phi(x;\mu,\sigma)\,G_1(d\mu,d\sigma)
		-
		\int_S
		\phi(x;\mu,\sigma)\,G_2(d\mu,d\sigma)
		\\
		&=
		\frac1M
		\int_S
		\phi(x;\mu,\sigma)\,\nu(d\mu,d\sigma)
		=
		0.
	\end{align*}
	Therefore,
	$
	p_{G_1}=p_{G_2}
	$
	while
	$	G_1\neq G_2.$
	Hence $\mathcal P(S)$ is not identifiable, contradicting the assumption.
	Therefore,
	$
	\operatorname{Int}(S)=\emptyset.
	$
\end{proof}


\subsection{Proof for Section \ref{sec:prelim}}\label{sec:proof-Sec2-Appendix}
\begin{proof}[Proof of Proposition \ref{prop:equi-Gauss-Hermite}]
        We first prove that
\begin{equation}\label{eq-phi}
	(T\nu)(x):=\int_S\phi_x\,d\nu=0
	\quad\text{for every }x\in\mathbb R
	\end{equation}
	if and only if
	\begin{equation}\label{eq-Hn}
	\int_SH_n\,d\nu=0
	\quad\text{for every }n\geq0.
	\end{equation}
	
	Assume that \eqref{eq-phi} holds. 
	Since
$
	\int_{\mathbb R}
	\phi(x;\mu,\sigma)\,dx
	=
	1,
$
	Tonelli's theorem gives
	\begin{align*}
		\int_{\mathbb R}|(T\nu)(x)|\,dx
		\leq
		\int_{\mathbb R}
		\int_S
		\phi(x;\mu,\sigma)\,
		|\nu|(d\mu,d\sigma)\,dx=
		\int_S
		\left(
		\int_{\mathbb R}
		\phi(x;\mu,\sigma)\,dx
		\right)
		|\nu|(d\mu,d\sigma)=
		|\nu|(S)<\infty.
	\end{align*}
	Thus
$
	T\nu\in L^1(\mathbb R).
$
	By assumption,
$
	T\nu(x)=0,\;
	\forall x\in\mathbb R,
$
taking the Fourier transform in $x$, we obtain
$$
	0
	=
	\widehat {T\nu}(t)
	=
	\int_{\mathbb R}
	e^{itx}(T\nu)(x)\,dx.
$$
	Using Fubini's theorem,
	\begin{align*}
		\widehat {T\nu}(t)
		=
		\int_S
		\left[
		\int_{\mathbb R}
		e^{itx}
		\phi(x;\mu,\sigma)\,dx
		\right]
		\nu(d\mu,d\sigma)=
		\int_S
		e^{it\mu-\frac12\sigma t^2}
		\,\nu(d\mu,d\sigma).
	\end{align*}
	Hence
\begin{equation*}
	F_\nu(t)
	:=
	\int_S
	e^{it\mu-\frac12\sigma t^2}
	\,d\nu
	=
	0,
	\qquad
	\forall t\in\mathbb R.
\end{equation*}
	
	Since $S$ is compact, 
	differentiation under
	the integral sign is valid in a neighborhood of $t=0$. Therefore,
$$
	F_\nu^{(n)}(0)
	=
	\int_S
	\left.
	\frac{d^n}{dt^n}
	e^{it\mu-\frac12\sigma t^2}
	\right|_{t=0}
	\,d\nu.
$$
	By the generating function, 
$
	\left.
	\frac{d^n}{dt^n}
	e^{it\mu-\frac12\sigma t^2}
	\right|_{t=0}
	=
	i^nH_n(\mu,\sigma).
$
	Hence,
$
	F_\nu^{(n)}(0)
	=
	i^n
	\int_SH_n\,d\nu.
$
	Since $F_\nu\equiv0$, all of its derivatives at zero vanish, and thus
$$
	\int_SH_n\,d\nu=0,
	\qquad
	\forall n\geq0.
$$
	This proves \eqref{eq-phi}$\Rightarrow$\eqref{eq-Hn}.
	
	Conversely, assume that \eqref{eq-Hn} holds. Define again
$
	F_\nu(t)
	=
	\int_S
	e^{it\mu-\frac12\sigma t^2}
	\,d\nu.
$
	For every $R>0$, the expansion
$$
	e^{it\mu-\frac12\sigma t^2}
	=
	\sum_{n=0}^{\infty}
	\frac{i^n}{n!}
	H_n(\mu,\sigma)t^n
$$
	converges uniformly on
$
	[-R,R]\times S.
$
	We may therefore integrate term by term to get
	\begin{align*}
		F_\nu(t)
		=
		\int_S
		\sum_{n=0}^{\infty}
		\frac{i^n}{n!}
		H_n(\mu,\sigma)t^n
		\,d\nu=
		\sum_{n=0}^{\infty}
		\frac{i^nt^n}{n!}
		\int_SH_n\,d\nu=
		0,\quad\forall t\in[-R,R].
	\end{align*}
	Since $R>0$ was arbitrary,
$
	F_\nu(t)=0,
	\;
	\forall t\in\mathbb R.
$
	
	As shown above,
$
	T\nu\in L^1(\mathbb R),
$
	and
$
	\widehat {T\nu}(t)=F_\nu(t)=0.
$
	By injectivity of the Fourier transform on $L^1(\mathbb R)$,
$
	T\nu=0
	\text{ a.e. on }\mathbb R.
$
	Moreover, $g_\nu$ is continuous. Indeed, if $x_j\to x$, then
$
	\phi(x_j;\mu,\sigma)
	\rightarrow
	\phi(x;\mu,\sigma)
$
	pointwise on $S$, and compactness of $S$ yields a uniform bound on the
	kernel for $x_j$ in a fixed bounded neighborhood of $x$. Dominated
	convergence therefore gives
$
	(T\nu)(x_j)\to (T\nu)(x).
$
	Thus a continuous function which vanishes almost everywhere must vanish
	everywhere, and hence
	$
    (T\nu)(x)=0,\; \forall x\in\mathbb R.
	$
	This proves \eqref{eq-Hn}$\Rightarrow$\eqref{eq-phi}.
	We have therefore established
\begin{equation}\label{eq-5}
		\int_S\phi_x\,d\nu=0
		\ \forall x
		\quad\Longleftrightarrow\quad
		\int_SH_n\,d\nu=0
		\ \forall n.
\end{equation}

    Define
$$
	\mathcal V_\phi^\perp
	=
	\left\{
	\nu\in\mathcal M(S):
	\int_Sf\,d\nu=0,
	\quad\forall f\in\mathcal V_\phi
	\right\}
$$
	and define $\mathcal V_H^\perp$ analogously. Since
	$\mathcal V_\phi$ and $\mathcal V_H$ are the closed linear spans of
	their respective generating families, \eqref{eq-5} implies
	\begin{equation}\label{eq-6}
	\mathcal V_\phi^\perp
	=
	\mathcal V_H^\perp.
	\end{equation}
	
	Now, suppose for contradiction that
$
	\mathcal V_\phi\neq\mathcal V_H.
$
	Without loss of generality, assume that there exists
$
	f\in\mathcal V_\phi\setminus\mathcal V_H.
$
	Since $\mathcal V_H$ is a closed linear subspace of $C(S)$,
	the Hahn--Banach theorem gives a continuous linear functional
	$L\in C(S)^*$ such that
$$
	L|_{\mathcal V_H}=0
	\qquad\text{and}\qquad
	L(f)\neq0.
$$
	By the Riesz-Markov-Kakutani representation theorem, there exists
	$\nu\in\mathcal M(S)$ such that
$
	L(g)=\int_Sg\,d\nu.
$
	Thus
$
	\nu\in\mathcal V_H^\perp.
$
	By \eqref{eq-6},
$
	\nu\in\mathcal V_\phi^\perp.
$
	Since $f\in\mathcal V_\phi$, this implies
$
	L(f)=0,
$
	contradicting the choice of $L$.
	Therefore,
$
	\mathcal V_\phi=\mathcal V_H.
$
\end{proof}

\begin{proof}[Proof of Theorem \ref{thm:iden-characterization}]
%
The equivalence of (i) and (ii) is straightforward.
    We next prove equivalence of (i) and (iii).
	Assume that
	$	\mathcal V_\phi=C(S).$
	Let $G_1,G_2\in\mathcal P(S)$ satisfy
	$
	p_{G_1}(x)=p_{G_2}(x),
	\;
\forall	x\in\mathbb{R}.
	$
	Define the finite signed measure
	$
	\nu:=G_1-G_2.
	$
	Then, for every $x\in\mathbb{R}$,
	$
	\int_S \phi_x\,d\nu
	=
	0.
	$
	Therefore,
	$$
	\int_S f\,d\nu=0,\quad
	\text{ for every }
	f\in
	\operatorname{span}
	\{\phi_x:x\in\mathbb{R}\}.
$$
	
	Now consider the linear functional
	$$
	L_\nu:C(S)\to\mathbb{R},
	\qquad
	L_\nu(f)=\int_S f\,d\nu.
	$$
	Since $\nu$ is a finite signed measure,
	$
	|L_\nu(f)|
	\leq
	\|f\|\,\|\nu\|_{\mathrm{TV}}(S),
	\;
\forall	f\in C(S),
	$
	where $\|\nu\|_{\mathrm{TV}}$ denotes the total variation measure of $\nu$.
	Hence $L_\nu$ is continuous with respect to the uniform norm.

    Given $\mathcal V_\phi=C(S)$,
	$
	\operatorname{span}
\{\phi_x:x\in\mathbb{R}\}
	$
	is dense in $C(S)$.
    So, for every $f\in C(S)$ there exists a sequence
	$\{f_n\}$ in this span such that
	$
	\|f_n-f\|\to0.
	$
	By continuity of $L_\nu$,
	$
	L_\nu(f)
	=
	\lim_{n\to\infty}L_\nu(f_n)
	=
	0.
	$
	Thus
	$
	\int_S f\,d\nu=0,
	\;
	\forall f\in C(S).
	$
	By the uniqueness in the Riesz representation theorem, we conclude that
$
	\nu=0.
$
	Consequently,
$
	G_1=G_2.
$
	Hence, $\mathcal P(S)$ is identifiable.

	Assume now that $\mathcal P(S)$ is identifiable. We prove that
	$
	\mathcal V_\phi=C(S).
	$
	Suppose, toward a contradiction, that
	$
	\mathcal V_\phi\neq C(S).
	$
	Since $\mathcal V_\phi$ is a proper closed linear subspace of the Banach
	space $C(S)$, the Hahn--Banach theorem implies that there exists a
	nonzero continuous linear functional
$$
	L\in C(S)^*,
	\;
	L\neq0,
\text{ such that }
	L(f)=0,
	\;
	\forall f\in\mathcal V_\phi.
$$
	
	By the Riesz-Markov-Kakutani representation theorem, there exists a nonzero finite signed
    measure $\nu$ on $S$ such that
	$$
	L(f)=\int_S f\,d\nu,
	\qquad
	f\in C(S).
	$$
	Since $\phi_x\in\mathcal V$ for every $x\in\mathbb{R}$, we have
	$$
    \int_S \phi_x\,d\nu
	=
	0,
	\;
	\forall x\in\mathbb{R}
    \qquad\text{and}\qquad
    \nu(S)
	=
	\int_S 1\,d\nu
	=
	\int_S\int\phi_xdxd\nu
	=
	0
	$$
    Therefore
	$
	\mathcal V_\phi=C(S).
	$
	The proof of equivalence of (i) and (iii) is complete.


Equivalence of (iii) and (iv) follows Proposition \ref{prop:equi-Gauss-Hermite}.
The theorem is therefore complete.
\end{proof}
\subsection{Proof for Section \ref{sec:topo}}\label{sec:proof34-Appendix}

\begin{proof}[Proof of Lemma \ref{lem:K(S)-compact}]
If $S_c=\emptyset$, then $\mathbb R\setminus S_c=\mathbb R$
has no bounded connected components, so $(\mathcal K(S))_c=\emptyset$.
Also, if $S_c=[\min S_c,\max S_c]$ then $(\mathcal K(S))_c=[\min S_c,\max S_c].$

Suppose now that $S_c\neq\emptyset$ and $S_c\neq [\min S_c,\max S_c]$. Since $S$ is compact,
$S_c$ is compact. 
As $\mathbb R\setminus S_c$ is open, each of its bounded connected
components is an open interval $(\alpha,\beta)$.
Both endpoints belong to $S_c$: otherwise, openness of
$\mathbb R\setminus S_c$ would allow the interval to be extended,
contradicting maximality of the component.
Thus $\min S_c\leq\alpha<\beta\leq \max S_c$, proving
$$
S_c\cup\!\!\bigcup_{\substack{J\text{ bounded connected}\\
    \text{component of }\mathbb R\setminus S_c}}J
\subseteq[\min S_c,\max S_c].
$$

Conversely, let $x\in[\min S_c,\max S_c]\setminus S_c$. The sets
$S_c\cap(-\infty,x]$ and $S_c\cap[x,\infty)$ are nonempty and compact.
Hence the numbers
$
\alpha:=\max\bigl(S_c\cap(-\infty,x]\bigr),
\;
\beta:=\min\bigl(S_c\cap[x,\infty)\bigr)
$
are well defined. Since $x\notin S_c$, we have
$\alpha<x<\beta$, and the definitions imply
$(\alpha,\beta)\cap S_c=\emptyset$.
Moreover, $\alpha,\beta\in S_c$, so no connected subset of
$\mathbb R\setminus S_c$ containing $(\alpha,\beta)$ can extend
beyond either endpoint. Thus $(\alpha,\beta)$ is a bounded
connected component containing $x$ and then, is contained in $(\mathcal K(S))_c$. These prove the claimed
description of $(\mathcal K(S))_c$.

Consequently, we can parameterize $\mathcal K(S)$ as the follows.
$$
\mathcal K(S)=
\{(1-t)p+tq:\ p,q\in S,\ \pi_\sigma(p)=\pi_\sigma(q),\ t\in[0,1]\}.
$$
The set
$$
D:=\{(p,q,t)\in S\times S\times[0,1]:\pi_\sigma(p)=\pi_\sigma(q)\}
$$
is closed in the compact space $S\times S\times[0,1]$,
and hence compact. The map
$$
F:D\longrightarrow\mathbb R^2,
\qquad F(p,q,t):=(1-t)p+tq,
$$
is continuous and satisfies $F(D)=\mathcal K(S)$.
Therefore $\mathcal K(S)$ is compact.
\end{proof}

\subsection{Proof for Section \ref{sec:5}}\label{sec:app-section5}

\begin{proof}[Proof of Collorary \ref{prop:iden-graph-over-location}]
We first prove (i).
	Quasiconvexity gives
	$$
    K:=\mathcal K(S)\subset\{(\mu,\sigma):\mu\in I,\ \sigma\ge s(\mu)\}.
	$$
	Now, fix $q=(\mu,s(\mu))\in S$.
	Almost surely, there are positive times $t_n\downarrow0$
	such that $B_{t_n}=0$. At these times,
	$$
	Y_{t_n}^q=(\mu,s(\mu)-t_n)\notin K.
	$$
	Consequently $\tau_{K}^q=0$ almost surely.
	Thus $S\cap O(K)=\emptyset$.

Now, we move to part (ii).
Let $\sigma^*=s(\mu^*)$ and $q^*=(\mu^*,\sigma^*)$.
Fix $a<\mu<\mu^*$. Since $s'(\mu)>0$, the inverse of the left branch is locally Lipschitz near $s(\mu)$. Thus, for some $C>0$, decreasing the scale by $h$ moves the left endpoint of the $\mu$-directional section of $K$ leftward by at most $Ch$. Consequently, for all sufficiently small $h>0$,
$$
\big\{\tau_{K}^{(\mu,s(\mu))}\geq h\big\}
\subset
\big\{Y_h^{(\mu,s(\mu))}=(\mu+B_h,s(\mu)-h)\in K\big\}
\subset
\{\mu+B_h\ge\mu-Ch\}
\subset\{B_h\ge-Ch\}.
$$
On the other hand, $\lim_{h\to 0}\PP\{B_h\ge-Ch\}=\frac 12$.
Therefore, by Lemma \ref{lem:exit-law}(ii), $(\mu,s(\mu))\notin O(K)$ for all $\mu\in(a,\mu^*)$. Similar conclusion is applied for $\mu\in(\mu^*,b)$.

It follows that $O(K)\subset\{q^*\}$.
If $q^*\notin O(K)$, then $S\cap O(K)=\emptyset$. Identifiability follows directly from Theorem \ref{thm:complete-iden-chacterization} (i).
Now, suppose that $q^*\in O(K)$. 
Let $\sigma_0=\max\{s(a),s(b)\}$.
For $\sigma\in(\sigma_0,\sigma^*)$, we will write $K_\sigma$ as 
$$
K_\sigma=[\ell(\sigma),r(\sigma)].
$$
Let
$$
F_0:=(\ell(\sigma_0),r(\sigma_0))\times\{\sigma_0\}\subset F:=\partial K\setminus S.
$$
We will claim that the exit law from $q^*$ satisfies
\begin{equation}\label{eq:apen-exit-law-1}
    L_{q^*}^{K}(F_0)>0.
\end{equation}


    
Indeed,  starting from
$q^*=(\mu^*,\sigma^*)$, the process is
$
Y_t^{q^*}=(\mu^*+B_t,\sigma^*-t).
$
Write $T=\sigma^*-\sigma_0$, and
$\tau^*=\tau_K^{q^*}$.
First, note that the fibers $K_c$ for $c< \sigma_0$ are
empty or singletons. Since $\mu^*+B_t$ has a density,
for every $t>T$ sufficiently close to $T$,
$$
\mathbb P(\tau^*>t)
\le
\mathbb P\bigl(\mu^*+B_t\in K_{\sigma^*-t}\bigr)
=0.
$$
Taking a sequence $t\downarrow T$ gives $\tau^*\le T$
almost surely.

Next, choose $\delta\in(0,T)$ such that
$\mathbb P(\tau^*>\delta)>0$.
On $\{\tau^*>\delta\}$, the point $Y_\delta^{q^*}$
belongs to $K$. Since $K_{\sigma^*-\delta}$ is an
interval and $\mu^*+B_\delta$ has a density, its
spatial coordinate lies in the interior of this
fiber almost surely on this event.
We may therefore choose $c\in\mathbb R$ and $r>0$
such that
$$
J:=[c-r,c+r]
\subset \operatorname{Int}K_{\sigma^*-\delta},
\qquad
\mathbb P\bigl(
\tau^*>\delta,\ |\mu^*+B_\delta-c|<r/2
\bigr)>0.
$$
The condition on the graph $S$ implies that the
fiber $K_{c}$ expands as $c$ decreases in $[\sigma_0,\sigma^*]$. 
The definition of
$F_0$ gives
$$
J\subseteq K_s
\quad\text{for }\sigma_0\le s\le\sigma^*-\delta,
\qquad
J\times\{\sigma_0\}\subseteq F_0.
$$
By independence of Brownian increments, the event
$$
E=
\Big\{
\tau^*>\delta,\quad
|\mu^*+B_\delta-c|<r/2,\quad
\sup_{\delta\le u\le T}
|B_{u}-B_\delta|<r/2
\Big\}
$$
has probability
$$
\mathbb P(E)
=
\mathbb P\Big(
\tau^*>\delta,\ |\mu^*+B_\delta-c|<r/2
\Big)
\,
\mathbb P\Big(
\sup_{\delta\le u\le T}|B_u - B_\delta|<r/2
\Big)
>0.
$$
On $E$, the spatial coordinate of $Y_t^{q^*}$
remains in $J$ for $\delta\le t\le T$.
Thus $Y_t^{q^*}\in K$ throughout this time interval,
and $Y_T^{q^*}\in F_0$.
Since $\tau^*>\delta$ on $E$, such constraints on $Y_t^{q^*}$ yield that $\tau^*\ge T$.
Together with $\tau^*\le T$ almost surely, we obtain
$\tau^*=T$ on $E$, almost surely. Consequently,
$$
L_{q^*}^{K}(F_0)
=
\mathbb P\bigl(Y_{\tau^*}^{q^*}\in F_0\bigr)
\ge \mathbb P(E)>0,
$$
proving the claim
\eqref{eq:apen-exit-law-1}.

Now, every $\lambda\in\mathcal M(S\cap O(K))$ has the form $\lambda=c\delta_{q^*}$.
By Theorem \ref{thm-null-space-characterization}, 
$$
0=(E_{K}\lambda)\vert_F=cL_{q^*}^{K}\vert_F.
$$
This together with \eqref{eq:apen-exit-law-1} forces $c=0$. So, $\mathcal N(S)=\{0\}.$
Therefore, part (ii) is proved.

The proof of part (iii) is similar to that of part (ii).
\end{proof}

\begin{proof}[Proof of Corollary \ref{cor:U-shape}]
	By continuity of $a$ and $b$, the set $S$ is compact.
	Its support hull is
	$$
	K:=\mathcal K(S)
	=
	\left\{
	(\mu,\sigma):
	\sigma_0\le\sigma\le\sigma_1,\ 
	a(\sigma)\le\mu\le b(\sigma)
	\right\}.
	$$
	 By Theorem \ref{thm:complete-iden-chacterization}, the identifiability of $\mathcal P(S)$ follows immediately once we verify that 
     $S \cap O(K) = \emptyset$. That is, 
	$$
	\mathbb P(\tau_K^q=0)=1
	\quad\text{for every }q\in S.
    $$

    	If $q=(\mu,\sigma_0)\in D_0$, then, for every $t>0$,
	the scale coordinate of $Y_t^q$ is
	$\sigma_0-t<\sigma_0$. Thus $Y_t^q\notin K$, and
	$
	\tau_K^q=0
	$
	for every continuous sample path.
	This also covers the 
    case
	$a(\sigma_0)=b(\sigma_0)$.

    	Next, fix $\sigma\in(\sigma_0,\sigma_1]$ and consider
	$q=(a(\sigma),\sigma)\in\Gamma_a$.
	The event $\{\tau_K^q>h\}$ implies $Y_h^q\in K$ and hence
	$
	a(\sigma)+B_h\ge a(\sigma-h).
	$ 
    For all $0<h<\delta_\sigma$,  
    using \eqref{eq-U-shape-regular-branch}, we obtain that 
    	$$
	\{\tau_K^q>h\}
	\subset
	\{B_h\ge-C_\sigma\sqrt h\}.
	$$
    	Let $\Phi$ denote the distribution function of a standard
	normal random variable. Since $B_h/\sqrt h$ is a standard
	normal variable,
	$$
	\mathbb P(\tau_K^q>h)
	\le
	\mathbb P(B_h\ge-C_\sigma\sqrt h)
	=
	\Phi(C_\sigma).
	$$
	Letting $h\downarrow0$ gives
	$
	\mathbb P(\tau_K^q>0)
	\le
	\Phi(C_\sigma)
	<1.
	$
    By Lemma \ref{lem:exit-law} (ii), $q \notin O(K)$, which means
	$
	\mathbb P(\tau_K^q>0)=0.
	$
If $q=(b(\sigma),\sigma)\in\Gamma_b$, the argument is similar.
The proof is therefore complete.
\end{proof}

\begin{proof}[Proof of Corollary \ref{cor:U-shape-sharp-threshold}]
    Set $S=\{a_\varepsilon\}\cup C_\varepsilon$. Its $\mu$-directional 
hull is
$$
K=\mathcal K(C_\varepsilon)
=
\{(\mu,\sigma):
(\mu-\mu_0)^2\le \sigma-\sigma_0\le\varepsilon^2\}.
$$
Every point of $C_\varepsilon$ is an immediate-exit point
of $K$. Indeed, starting at
$q=(\mu_0+y,\sigma_0+y^2)$, the 
the backward
space-time Brownian motion
$
Y_t^q=(\mu_0+y+B_t,\sigma_0+y^2-t)$
lies outside $K$
whenever
$$
(y+B_t)^2>y^2-t,
\qquad\text{equivalently}\qquad
2yB_t+B_t^2+t>0.
$$
For $y=0$, it holds for every
$t>0$.
For $y\ne0$, 
Brownian motion almost surely takes values
of the same sign as $y$ at arbitrarily small positive
times. At each such time,
$
2yB_t+B_t^2+t>0,
$
so $Y_t^q\notin K$. Thus the exit time from $K$ is zero
almost surely.

Starting instead at $a_\varepsilon=(\mu_0+d,\sigma_0+\eps^2)$, define
$$
\tau_\varepsilon
=
\inf\{t\ge0:(d+B_t)^2+t>\varepsilon^2\},
\qquad
G_\varepsilon
=
\operatorname{Law}
\bigl(\mu_0+d+B_{\tau_\varepsilon},
      \sigma_0+\varepsilon^2-\tau_\varepsilon\bigr).
$$
By a.s. continuity of the Brownian sample paths, we have almost surely
$
(d+B_{\tau_\varepsilon})^2+\tau_\varepsilon
=\varepsilon^2,
$ and $0<\tau_\varepsilon\le\varepsilon^2$.
Thus $G_\varepsilon=E_K\delta_{a_\varepsilon}$ belongs
to $\mathcal P(C_\varepsilon)$. Since $a_\varepsilon$
is the only active point of $S$ and its exit law is
concentrated on $C_\varepsilon\subset S$, the null-space
characterization theorem (Theorem \ref{thm-null-space-characterization})
gives
$$
\mathcal N(S)
=
\{c(\delta_{a_\varepsilon}-G_\varepsilon):
c\in\mathbb R\}.
$$

We now claim that $\operatorname{supp}G_\varepsilon=C_\varepsilon$.
Fix an arbitrary $y \in (-\varepsilon, \varepsilon)$, set $T=\varepsilon^2-y^2$, and
consider the linear path
$$
f(t)=\frac{y-d}{T}t,
\qquad 0\le t\le2\varepsilon^2,
$$
which has the property that $f(0) = 0$ and $f(T) = y-d$. In particular, at $t=T$, the stopping identity $(d+f(t))^2+t = \epsilon^2$ holds, and the
coordinate $(\mu_0+d+f(T), \sigma_0 + \varepsilon^2-T) = (\mu_0+y,\sigma_0+y^2)$ corresponds to a point in $C_\epsilon$. Notice that the Wiener measure has full support on $C_0([0,2\varepsilon^2]):=\{\omega\in C[0,2\varepsilon^2]: \omega(0)=0\}$, so every neighborhood under the uniform norm of this linear path $f$ has a positive $G_\varepsilon$-probability.

Since the function
$$
h(t)=(d+f(t))^2+t-\varepsilon^2
$$
is convex, satisfies $h(0)=d^2-\varepsilon^2<0$ and
$h(T)=0$, and has
$$
h'(T)
=
\frac{(y-d)^2+\varepsilon^2-d^2}{T}>0.
$$
we obtain that $h(t)<0$ for $0\le t<T$ and $h(t)>0$
for $t>T$. These strict inequalities imply that
sufficiently small uniform perturbations of $f$ have
exit times close to $T$ and exit locations close to
$(\mu_0+y,\sigma_0+y^2)$.
Thus $(\mu_0+y,\sigma_0+y^2)$ must be contained in the support of 
the exit law $G_\epsilon$.
The two endpoints follow by closedness of support,
proving the claim.

Since $\operatorname{supp}G_\varepsilon=C_\varepsilon$
and $a_\varepsilon\notin C_\varepsilon$, the measure
$\delta_{a_\varepsilon}-G_\varepsilon$ is nonzero and
is supported exactly on $S$. In particular,
$\mathcal N(S)\ne\{0\}$, proving \textup{(i)}.

Now let $S'\subsetneq S$ be compact and take
$\nu\in\mathcal N(S')$. Viewing $\nu$ as a measure on
$S$, the above characterization gives
$$
\nu=c(\delta_{a_\varepsilon}-G_\varepsilon)
$$
for some $c\in\mathbb R$. If $c\ne0$, then
$\operatorname{supp}\nu=S$, contradicting
$\operatorname{supp}\nu\subseteq S'$.
Thus $\nu=0$, and $\mathcal P(S')$ is identifiable.
This proves \textup{(ii)}.
\end{proof}
\begin{proof}[Proof Corollary \ref{cor:one-point-full-support-target}]
We observe that
	$
	\mathcal M(\{a\})
	=
	\{c\delta_a:c\in\mathbb R\}
	$
    and that
	\begin{equation}
		\label{eq:one-dim-null}
		\mathcal N(\{a\}\cup\Gamma)
		=
		\left\{
		c(\delta_a-G_a):
		c\in\mathbb R
		\right\}.
	\end{equation}
	Since
	$
	\delta_a-G_a\ne0,
	$
	the set
	$
	\{a\}\cup\Gamma
	$
	is non-identifiable.
	
	Now let
	$
	\Gamma'\subsetneq\Gamma
	$
	be compact. Suppose
	$
	\nu\in\mathcal N(\{a\}\cup\Gamma').
	$
	By
	\eqref{eq:one-dim-null},
	$
	\nu
	=
	c(\delta_a-G_a)
	$
	for some $c\in\mathbb R$.
	If $c\ne0$, then the restriction of $\nu$ to $\Gamma$ has support
	$
	\operatorname{supp}(cG_a)
	=
	\Gamma,
	$
	by
	\eqref{eq:point-rep-measure}. This is impossible
	because $\nu$ is supported on
	$
	\{a\}\cup\Gamma',
	\;
	\Gamma'\subsetneq\Gamma.
	$
	Hence
	$
	c=0,
	$
	and therefore
	$
	\nu=0.
	$
	Thus
	$
	\mathcal P(\{a\}\cup\Gamma')
	$
	is identifiable.
\end{proof}

\begin{proof}[Proof of Corollary \ref{thm:pq-threshold}]
Write $S=S_{p,q}$. Its support hull is
$$
K=\mathcal K(S)=\overline{\Omega_{p,q}},
\qquad \partial K=S.
$$
Since the entire boundary belongs to $S$, by Theorem \ref{thm:complete-iden-chacterization} (ii),
the question of identifiability is reduced to determining 
the active points of $S$, namely the set $S \cap O(K)$.


We first show that every boundary point other than
$a_+=(0,s_0+1)$ is an immediate-exit point.
The bottom point $a_-=(0,s_0-1)$ exits immediately
because its variance coordinate decreases below $s_0-1$.
At any other point $q=(\mu,\sigma)\in S$ with $\mu\ne0$,
convexity of $x\mapsto |x|^p$ (for $p\geq 1$) and local Lipschitz
continuity of $v\mapsto |v|^q$ give
$$
|\mu+B_t|^p+|\sigma-s_0-t|^q
\ge
1+p\,\operatorname{sgn}(\mu)|\mu|^{p-1}B_t-Ct
$$
for all sufficiently small $t>0$ and some $C<\infty$.
The law of the iterated logarithm implies
$$
\limsup_{t\downarrow0}
\frac{\operatorname{sgn}(\mu)B_t}{t}
=\infty
\qquad\text{almost surely}.
$$
Thus $Y_t^q\notin K$ at arbitrarily small positive
times, proving that $\tau_K^q=0$ almost surely.

It remains to examine $a_+$. For $0<t<1$,
$$
|B_t|^p+|1-t|^q-1
=
|B_t|^p-qt+o(t).
$$
If $1\le p\le2$, the law of the iterated logarithm gives
$$
\limsup_{t\downarrow0}\frac{|B_t|^p}{t}
=\infty
\qquad\text{almost surely}.
$$
Consequently, the preceding expression is positive
at arbitrarily small times, so $\tau_K^{a_+}=0$
almost surely. Since $S\cap O(K) = \emptyset$, Theorem \ref{thm:complete-iden-chacterization}
yields
$\mathcal N(S)=\{0\}$.

If $p>2$, choose $\alpha$ with $1/p<\alpha<1/2$.
Almost surely, Brownian paths are locally
$\alpha$-H\"older continuous, and hence
$
|B_t|^p=O(t^{\alpha p})=o(t).
$
It follows that
$$
|B_t|^p+|1-t|^q-1=-qt+o(t)<0
$$
for all sufficiently small positive $t$, almost surely.
Thus $\tau_K^{a_+}>0$ almost surely, and $a_+$ is the
unique element of $O(K)$.

Let
$
G=E_K\delta_{a_+}$
Theorem \ref{thm-null-space-characterization} gives
$
\mathcal N(S)
=
\{c(\delta_{a_+}-G):c\in\mathbb R\}.
$
Moreover, $\tau_K^{a_+}>0$ implies $G(\{a_+\})=0$,
so $\delta_{a_+}-G\ne0$. Hence $\mathcal P(S)$
is not identifiable when $p>2$. That completes the proof.
\end{proof}

\begin{proof}[Proof of Proposition \ref{prop:three}]
The identifiability of $\mathcal P(\Gamma_1\cup\Gamma_2)$ follows from Corollary \ref{cor:U-shape}.
To establish non-identifiability of $\mathcal P(S)$, it suffices to construct explicitly a nonzero finite signed measure $\nu$ supported
	on $S$ such that
	\begin{equation}\label{eq:counter}
		\int_S
		e^{z\mu+\frac12z^2\sigma}
		\,\nu(d\mu,d\sigma)
		=
		0,
		\qquad
		\forall z\in\mathbb C.
	\end{equation}
	Once $\nu$ is constructed, we differentiate
		\eqref{eq:counter} at $z=0$, and get that
	for every $n\geq0$,
	$$
	0
	=
	\left.
	\frac{d^n}{dz^n}
	\int_S
	e^{z\mu+\frac12z^2\sigma}\,d\nu
	\right|_{z=0}
	=
	\int_S
	H_n(\mu,\sigma)\,d\nu.
	$$
        (This is admissible because
		$S$ is compact, and thus, the generating series
		$$
		e^{z\mu+\frac12z^2\sigma}
		=
		\sum_{n=0}^\infty
		H_n(\mu,\sigma)\frac{z^n}{n!}
		$$
		converges uniformly on $S$ for $z$ in every compact subset of
		$\mathbb C$). 
	Hence
	$
	\int_S H_n\,d\nu=0,
	\;\forall
	n=0,1,2,\ldots.
	$
	Since $\nu\neq0$, the Riesz characterization of density
	immediately gives
	$$
	\overline{
		\operatorname{span}
		\{H_n|_S:n\geq0\}
	}
	\neq C(S).
	$$
	Therefore, by Theorem \ref{thm:iden-characterization}, $\mathcal P(S)$ is not identifiable.
    
	Now, we construct nonzero signed measure $\nu$ satisfying \eqref{eq:counter}.

    \bigskip
	
	\noindent
	\textbf{Step 1: Reduction to three one-dimensional Laplace transforms.}
	For $j=1,2,3$, let
	$$
	\gamma_j:[T_0,T_1]\longrightarrow\Gamma_j,
	\qquad
	\gamma_j(t)=(a_j+ct,t).
	$$
	We shall construct signed measures on $\Gamma_j$ via the pushforward measures
	$$
	\nu_j
	=
	(\gamma_j)_\#
	\bigl(g_j(t)\,dt\bigr),$$
	for some suitable $g_j\in C_c^\infty((T_0,T_1))
	$ to be chosen, and then set 
	$\nu:=\nu_1+\nu_2+\nu_3$.
	
	Define the Laplace transforms
	$$
	L[g_j](\lambda)
	:=
	\int_{T_0}^{T_1}
	e^{\lambda t}g_j(t)\,dt,
	\qquad
	\lambda\in\mathbb C.
	$$
	Since $g_j$ has compact support, $L[g_j]$ is entire.
	For $z\in\mathbb C$,
	\begin{align*}
		\int_{\Gamma_j}
		e^{z\mu+\frac12z^2\sigma}
		\,d\nu_j(\mu,\sigma)
		=
		\int_{T_0}^{T_1}
		e^{z(a_j+ct)+\frac12z^2t}
		g_j(t)\,dt
		=
		e^{a_jz}
		L[g_j]\left(cz+\frac12z^2\right).
	\end{align*}
	Thus \eqref{eq:counter} will follow if
	\begin{equation}\label{eq:three-L-condition}
		\sum_{j=1}^3
		e^{a_jz}
		L[g_j]\left(cz+\frac12z^2\right)
		=
		0,
		\qquad
		\forall z\in\mathbb C.
	\end{equation}
	
   Let
    $\lambda := cz+\frac{1}{2}z^2 = \frac{1}{2}(w^2-c^2)$,
    where
	$
	w:=z+c.
	$
	Then,
	$
	e^{a_jz}
	=
	e^{-ca_j}e^{a_jw}.
	$
	Let
	$
	A_j:=e^{-ca_j}>0$
    and $L_j$ stand for $L[g_j]$, for $j=1,2,3$.
	Then \eqref{eq:three-L-condition} becomes
	$$
	\sum_{j=1}^3
	A_je^{a_jw}L_j(\lambda)
	=
	0.
	$$
	Since the same value of $\lambda$ corresponds to both $w$ and $-w$,
	it is natural to impose simultaneously
	$$
	\sum_{j=1}^3
	A_je^{a_jw}L_j(\lambda)=0,
	\text{ and }
	\sum_{j=1}^3
	A_je^{-a_jw}L_j(\lambda)=0.
	$$
	These constitute two homogeneous linear equations for the three
	quantities $L_1(\lambda),L_2(\lambda),L_3(\lambda)$.

\medskip
    
	\noindent
	\textbf{Step 2: An explicit null vector.}	
	Consider the two row vectors
	$$
	\bigl(
	A_1e^{a_1w},
	A_2e^{a_2w},
	A_3e^{a_3w}
	\bigr)
	\text{ and }
	\bigl(
	A_1e^{-a_1w},
	A_2e^{-a_2w},
	A_3e^{-a_3w}
	\bigr),
	$$
	and	their cross product, which gives a vector
	$
	\bigl(
	N_1(w),N_2(w),N_3(w)
	\bigr)
	$
	orthogonal to both of them.
	A direct computation gives
    \begin{eqnarray*}
	N_1(w)
	& = &
	-2A_2A_3
	\sinh\bigl((a_3-a_2)w\bigr),\\
	N_2(w)
	& = &
	2A_1A_3
	\sinh\bigl((a_3-a_1)w\bigr),\\
	N_3(w)
	& = &
	-2A_1A_2
	\sinh\bigl((a_2-a_1)w\bigr).
    \end{eqnarray*}
	Therefore
	\begin{equation}\label{eq:null-plus}
		\sum_{j=1}^3
		A_je^{a_jw}N_j(w)
		=
		0,
		\text{ and }	
		\sum_{j=1}^3
		A_je^{-a_jw}N_j(w)
		=
		0.
	\end{equation}
	
	\medskip
	
	\noindent
	\textbf{Step 3: Compactly supported realization of Laplace multipliers.}
	In this step, we shall construct the functions $g_j \in C_c^\infty((T_0,T_1))$ whose Laplace transform $L[g_j](\lambda)$ is proportional to the functions $N_j(w)$ just defined.
    This is the most technical part of the proof.
	The main idea is to exploit Taylor series of hyperbolic sine function:
	$$
	\sinh z=\sum_{n=0}^{\infty}\frac{z^{2n+1}}{(2n+1)!}.
	$$
	Since the null vector constructed above naturally produces hyperbolic sine factors, we seek a function whose Laplace transform generates such a multiplier.  

    Take a function $\varphi \in C_c^\infty((T_0,T_1))$, integration by parts gives
	$$\int_{T_0}^{T_1} e^{\lambda t}\varphi'(t)\,dt = -\lambda \int_{T_0}^{T_1} e^{\lambda t}\varphi(t)\,dt = -\lambda L[\varphi](\lambda). $$
    Note the relation between the two variables $w$ and $\lambda$ defined earlier: $w^2 = c^2 + 2\lambda$.

    Thus, applying the differential operator $c^2-2\partial_t$ to $\varphi$ results in the Laplace transform 
	$$
	\int
	e^{\lambda t}
	(c^2-2\partial_t)\varphi(t)\,dt
	=
	(c^2+2\lambda)L[\varphi](\lambda) = w^2 L[\varphi(\lambda)],
	$$
	and, by iterating, we obtain for any $n \geq 0$, 
	$$
	\int
	e^{\lambda t}
	(c^2-2\partial_t)^n\varphi(t)\,dt
	=
	w^{2n} L[\varphi](\lambda).
	$$

    For $d>0$, define formally
	\begin{equation}\label{eq:Gd-def}
		G_d(t)
		:=
		\sum_{n=0}^\infty
		\frac{d^{2n+1}}{(2n+1)!}
		(c^2-2\partial_t)^n \varphi(t).
	\end{equation}

    We will verify in the sequel that for a suitably chosen $\varphi \in C_c^\infty((T_0,T_1))$, the infinite series that defines $G_d$ converges uniformly in $(T_0,T_1)$, and in fact  
	\begin{equation}\label{eq:Gd-uniconv}
    G_d \in C_c^\infty((T_0,T_1)).
    \end{equation}
    It follows that the Laplace transform of $G_d$ may be obtained by applying the transform to the 
	the series term by term:
	\begin{align*}
		\int
		e^{\lambda t}G_d(t)\,dt
		&=
		L[\varphi](\lambda)
		\sum_{n=0}^\infty
		\frac{d^{2n+1}w^{2n}}{(2n+1)!} = L[\varphi](\lambda)\frac{\sinh(dw)}{w}.
	\end{align*}
	Hence
	\begin{equation}\label{eq:Gd-Laplace}
		L[G_d](\lambda)
		=
		L[\varphi](\lambda)
		\frac{\sinh(d\sqrt{c^2+2\lambda})}
		{\sqrt{c^2+2\lambda}},
	\end{equation}
	where the quotient is understood through its power-series
	expression. In particular, the right-hand side is a single-valued
	entire function of $\lambda$.

	%

    It remains in this proof step to verify \ref{eq:Gd-uniconv}. Choose a number
	$
	\rho\in(1,2).
	$
	There exists a nonzero compactly supported function
	$
	\varphi\in C_c^\infty((T_0,T_1))
	$ for which there are constants $C,R>0$ such that
	\begin{equation}\label{eq:Gevrey-bound}
		\|\varphi^{(m)}\|_\infty
		\leq
		CR^m(m!)^\rho,
		\qquad
		m=0,1,2,\ldots.
	\end{equation}
    The class of functions with this property are known as Gevrey functions of order $\rho$. Since $\varphi$ is nonzero, its Laplace transform is nonzero as well. 
	
	
    For any fixed $m\geq0$,
	$$
	\partial_t^m
	(c^2-2\partial_t)^n\varphi
	=
	\sum_{\ell=0}^n
	\binom{n}{\ell}
	c^{2(n-\ell)}
	(-2)^\ell
	\varphi^{(m+\ell)}.
	$$
	Using \eqref{eq:Gevrey-bound}, there exists a constant $C_m>0$
	such that
	\begin{equation}\label{eq:operator-growth}
		\left\|
		\partial_t^m
		(c^2-2\partial_t)^n\varphi
		\right\|_\infty
		\leq
		C_m^{\,n+1}
		((m+n)!)^\rho.
	\end{equation}
	For fixed $m$,
	$$
	((m+n)!)^\rho
	\leq
	C_m'(n!)^\rho n^{m\rho}
	$$
	for a suitable constant $C_m'>0$. On the other hand, Stirling's
	formula gives
	$$
	(2n+1)!
	\asymp
	4^n(n!)^2
	\sqrt{n}
	$$
	up to polynomial factors. Since
	$\rho<2,$
	the factorial growth in $(2n+1)!$ dominates the bound in
	\eqref{eq:operator-growth}. Hence, for every $m\geq0$, the series
	$$
	\sum_{n=0}^\infty
	\frac{d^{2n+1}}{(2n+1)!}
	\partial_t^m
	(c^2-2\partial_t)^n\varphi
	$$
	converges uniformly.
	Thus \eqref{eq:Gd-def} defines
	$
	G_d\in C^\infty(\mathbb R).
	$
	Furthermore, differentiation does not enlarge support, and every term
	in \eqref{eq:Gd-def} is supported in $\operatorname{supp}\varphi$.
	Uniform convergence therefore yields
	$$
	\operatorname{supp}G_d
	\subseteq
	\operatorname{supp}\varphi
	\subset
	(T_0,T_1).
	$$
    Thus \ref{eq:Gd-uniconv} holds.
	
	

\medskip
    
	\noindent
	\textbf{Step 4: Definition of the three densities.} With the function $G_d$ obtained, 
	set
	$$
	d_{23}:=a_3-a_2>0,
	\qquad
	d_{13}:=a_3-a_1>0,
	\qquad
	d_{12}:=a_2-a_1>0,
	$$
	and
	\begin{equation}\label{eq-define-gj}
		g_1
		:=
		-2A_2A_3\,G_{d_{23}},
		\qquad
		g_2
		:=
		2A_1A_3\,G_{d_{13}},
		\quad 
		\text	{and }
		g_3
		:=
		-2A_1A_2\,G_{d_{12}}.
	\end{equation}
	By \eqref{eq:Gd-Laplace}, their Laplace transforms are
	$$
	L[g_1](\lambda)
	=
	-2A_2A_3 L[\varphi](\lambda)
	\frac{\sinh((a_3-a_2)w)}{w},
	$$
	$$
	L[g_2](\lambda)
	=
	2A_1A_3 L[\varphi](\lambda)
	\frac{\sinh((a_3-a_1)w)}{w},
	$$
	and
	$$
	L[g_3](\lambda)
	=
	-2A_1A_2 L[\varphi](\lambda)
	\frac{\sinh((a_2-a_1)w)}{w},
	$$
	where
	$	
	w^2=c^2+2\lambda.
	$
	That says
	\begin{equation}\label{eq:Lj-nullvector}
		L[g_j](\lambda)
		=
		L[\varphi](\lambda)\frac{N_j(w)}{w},
		\qquad
		j=1,2,3.
	\end{equation}
	Using \eqref{eq:Lj-nullvector} and \eqref{eq:null-plus},
	\begin{align*}
		\sum_{j=1}^3
		e^{a_jz}L[g_j](\lambda)
		&= L[\varphi](\lambda)
		\sum_{j=1}^3
		A_je^{a_jw}
		\frac{N_j(w)}{w}
		=
		\frac{L[\varphi](\lambda)}{w}
		\sum_{j=1}^3
		A_je^{a_jw}N_j(w)
		=0.
	\end{align*}
	At $w=0$ the expression is interpreted by analytic continuation;
	all functions involved are entire, so the identity remains valid.
	
	Therefore
	$$
	\sum_{j=1}^3
	e^{a_jz}
	L_j\left(cz+\frac12z^2\right)
	=
	0,
	\qquad
	\forall z\in\mathbb C.
	$$
	As a result,
	$$
	\int_S
	e^{z\mu+\frac12z^2\sigma}\,d\nu(\mu,\sigma)
	=
	0,
	\qquad
	\forall z\in\mathbb C.
	$$
		
	\noindent
	\textbf{Step 5: Verifying the signed measure $\nu$ is nonzero.}
	We verify that the preceding construction is nontrivial.
	For every $d>0$, the entire multiplier
	$$
	M_d(\lambda)
	:=
	\frac{\sinh(d\sqrt{c^2+2\lambda})}
	{\sqrt{c^2+2\lambda}}
	$$
	is non-zero. Since
	$
	L[\varphi]\not\equiv0,
	$
	the product
	$
	L[\varphi](\lambda)M_d(\lambda)
	$
	is not identically zero.
	Therefore
	$
	L[G_d]\not\equiv0.
	$
	By uniqueness of the Laplace transform of compactly supported finite
	measures,
	$
	G_d\not\equiv0.
	$
	In particular, each $g_j$ in \eqref{eq-define-gj} is nonzero.
	Moreover, the three curves $\Gamma_j$ are pairwise disjoint. 
	Consequently the nonzero component measures $\nu_j$ have mutually
	disjoint supports, and hence
	$
	\nu=\nu_1+\nu_2+\nu_3\neq0.
	$
\end{proof}

\begin{rem}[Nature of the three-branch structure]
	The mechanism behind Proposition \ref{prop:three}
	is the two-sheeted nature of the quadratic spectral map
	$
	q(z)
	=
	cz+\frac12z^2.
	$
	For a generic value of
	$
	\lambda=q(z),
	$
	there are exactly two corresponding spectral parameters,
	$
	z_\pm
	=
	-c\pm\sqrt{c^2+2\lambda}.
	$
	For a union of $N$ parallel affine branches
	$
	\mu=a_j+c\sigma,
	\; j=1,\ldots,N,
	$
	the annihilation condition therefore produces, at each generic
	$\lambda$, two linear equations
	$$
	\sum_{j=1}^N
	e^{a_jz_+}L_j(\lambda)=0,
	\qquad
	\sum_{j=1}^N
	e^{a_jz_-}L_j(\lambda)=0
	$$
	for the $N$ unknown Laplace transforms $L_j(\lambda)$.
	
	When $N=2$, the resulting $2\times2$ system is generically invertible
	provided the two branches are distinct, which forces both Laplace
	transforms to vanish. One may draw a parallel to the numerous examples given in Section \ref{sec:5} where the the support set $S$ can be seen as possessing a two-branch structure for which identifiability holds.
    In contrast, when $N=3$, the system has a
	nontrivial pointwise nullspace.  The construction above shows that a
	nonzero vector in this nullspace can in fact be realized by Laplace
	transforms of smooth compactly supported signed densities.  
    Thus, the analysis of three parallel branches reflects a 
    crucial feature of the global geometric structure of the support set $S$
    which leads to non-identifiability.
\end{rem}

\begin{proof}[Proof of Theorem \ref{prop:rep-mixing-to-boundary}]
In both cases, the previously established exit properties
give that $\mathcal K(\Gamma)=K$, and that  every point of
$\Gamma$ immediately exit $K$, and that an exit law is concentrated on $\Gamma$ from every starting point in $K$. The null-space characterization
therefore implies that $\mathcal P(\Gamma)$ is identifiable.
Since the exit operator preserves Gaussian output,
$E_KG$ provides the unique representative in (i).

To prove (ii), fix $q\in\Gamma$. For every $z\in K$,
$$
(E_K\delta_z)(\{q\})=\mathbf 1_{\{z=q\}}.
$$
Indeed, immediate exit gives the identity when $z=q$.
For $z\ne q$, an exit at $q$ would require Brownian motion
to take a prescribed value at a fixed positive time,
determined by the difference in scale coordinates;
this event has probability zero. Hence
$$
(E_KG)(\{q\})=G(\{q\})
\qquad\text{for every }G\in\mathcal P(K).
$$
If $\mathcal G(A)=\mathcal G(K)$, choose
$G\in\mathcal P(A)$ with $p_G=p_{\delta_q}$.
The uniqueness in (i) gives that $E_KG=\delta_q$,
so $G(\{q\})=1$ and therefore $q\in A$.
Thus $\Gamma\subseteq A$. The converse follows from
$\Gamma\subseteq A\subseteq K$ and (i).
\end{proof}

\bigskip

\noindent{\bf Acknowledgements.} 
Part of this research was performed during both authors' visits to 
Vietnam Institute for Advanced Study in Mathematics (VIASM) in summer 2024 and
summer 2026. We wish to thank VIASM for a stimulating environment and the kind hospitality. ChatGPT 5.6-Sol was used to assist with literature searches
and for some mathematical discussions during this work.
The authors take full responsibility for the content
of the manuscript.

\bibliographystyle{amsalpha}
\bibliography{iden}
\end{document}